\documentstyle{article}
\def\noheaderplainsetup{%

     \topmargin=50pt \headheight=0pt \headsep=0pt  
     \oddsidemargin=45pt \evensidemargin=45pt       
     \textheight=7.48truein \textwidth=5.2truein}  

\noheaderplainsetup

\begin{document}


\newcommand{\inp}[1]{{\mathbf Int}^{#1}}
\newcommand{\lit}{\mbox{\bf Literal}}
\newcommand{\por}{\mbox{\bf Port}}
\newcommand{\fn}[1]{\mbox{\bf Tfn}^{#1}}
\newcommand{\str}[1]{\mbox{\bf St}^{#1}}
\newcommand{\rsp}[1]{\mbox{\bf Pl}^{#1}} 

\newcommand{\seq}[1]{\langle #1 \rangle}           


\newcommand{\mla}{\mbox{{\Large $\wedge$}}}
\newcommand{\mle}{\mbox{{\Large $\vee$}}}

\newcommand{\pst}{\mbox{\raisebox{-0.01cm}{\scriptsize $\wedge$}\hspace{-4pt}\raisebox{0.16cm}{\tiny $\mid$}\hspace{2pt}}}
\newcommand{\gneg}{\neg}                  
\newcommand{\mli}{\rightarrow}                     
\newcommand{\cla}{\mbox{\large $\forall$}}      
\newcommand{\cle}{\mbox{\large $\exists$}}        
\newcommand{\mld}{\vee}    
\newcommand{\mlc}{\wedge}  
\newcommand{\ade}{\mbox{\Large $\sqcup$}}      
\newcommand{\ada}{\mbox{\Large $\sqcap$}}      
\newcommand{\add}{\sqcup}                      
\newcommand{\adc}{\sqcap}                      

\newcommand{\tlg}{\bot}               
\newcommand{\twg}{\top}               
\newcommand{\st}{\mbox{\raisebox{-0.05cm}{$\circ$}\hspace{-0.13cm}\raisebox{0.16cm}{\tiny $\mid$}\hspace{2pt}}}
\newcommand{\cost}{\mbox{\raisebox{0.12cm}{$\circ$}\hspace{-0.13cm}\raisebox{0.02cm}{\tiny $\mid$}\hspace{2pt}}}
\newcommand{\pcost}{\mbox{\raisebox{0.12cm}{\scriptsize $\vee$}\hspace{-4pt}\raisebox{0.02cm}{\tiny $\mid$}\hspace{2pt}}}


\newtheorem{theoremm}{Theorem}[section]
\newtheorem{thesiss}[theoremm]{Thesis}
\newtheorem{definitionn}[theoremm]{Definition}
\newtheorem{lemmaa}[theoremm]{Lemma}
\newtheorem{propositionn}[theoremm]{Proposition}
\newtheorem{conventionn}[theoremm]{Convention}
\newtheorem{examplee}[theoremm]{Example}
\newtheorem{remarkk}[theoremm]{Remark}
\newtheorem{factt}[theoremm]{Fact}
\newtheorem{exercisee}[theoremm]{Exercise}

\newenvironment{exercise}{\begin{exercisee} \em}{ \end{exercisee}}
\newenvironment{definition}{\begin{definitionn} \em}{ \end{definitionn}}
\newenvironment{theorem}{\begin{theoremm}}{\end{theoremm}}
\newenvironment{lemma}{\begin{lemmaa}}{\end{lemmaa}}
\newenvironment{proposition}{\begin{propositionn} }{\end{propositionn}}
\newenvironment{convention}{\begin{conventionn} \em}{\end{conventionn}}
\newenvironment{remark}{\begin{remarkk} \em}{\end{remarkk}}
\newenvironment{proof}{ {\bf Proof.} }{\  $\Box$ \vspace{.1in} }
\newenvironment{example}{\begin{examplee} \em}{\end{examplee}}
\newenvironment{fact}{\begin{factt}}{\end{factt}}

\title{Introduction to cirquent calculus and abstract resource semantics}
\author{Giorgi Japaridze\thanks{This material is based upon work supported by the National Science Foundation under Grant No. 0208816, and 2005 Summer Research Grant from Villanova University.} 
\\  \\ Department of Computing Sciences, Villanova University, 
\\ 800 Lancaster Avenue, Villanova, PA 19085, USA\\
Email: giorgi.japaridze@villanova.edu\\
URL: http://www.csc.villanova.edu/$^\sim$japaridz/
}
\date{}
\maketitle

\begin{abstract} This paper introduces a refinement of the sequent calculus approach called {\em cirquent calculus}. Roughly speaking, the difference between the two is that, while in Gentzen-style proof trees sibling (or cousin, etc.) sequents are disjoint and independent sequences of formulas, in cirquent calculus they are permitted to share elements. Explicitly allowing or disallowing shared resources and thus taking to a more subtle level the resource-awareness intuitions underlying substructural logics, cirquent calculus offers much greater flexibility and power than sequent calculus does. A need for substantially new deductive tools came with the birth of {\em computability logic}  --- the semantically constructed formal theory of computational resources, which has stubbornly resisted all axiomatization attempts within the framework of traditional syntactic approaches. Cirquent calculus breaks the ice. Removing contraction from the full (``classical'') collection of its rules yields a sound and complete system for the basic fragment {\bf CL5} of computability logic, previously thought to be ``most unaxiomatizable''. Deleting the offending rule of contraction in ordinary sequent calculus, on the other hand, throws out the baby with the bath water, resulting in the strictly 
weaker affine logic. An implied claim of computability logic is that 
 it is {\bf CL5} rather than affine logic that adequately materializes the resource philosophy traditionally associated with the latter. To strengthen this claim, 
the paper further introduces an abstract resource semantics and shows the soundness and completeness of {\bf CL5} with respect to it. Unlike the semantics of computability logic, which understands resources in a special --- computational --- sense, abstract resource semantics can be seen as a direct formalization of the more general yet naive intuitions in the ``can you get both a candy and an apple for one dollar?'' style. The inherent incompleteness of affine or linear logics, resulting from the fundamental limitations of the underlying sequent-calculus approach, is apparently the reason why such intuitions
and examples, while so heavily relied on in the popular linear-logic literature, have never really found a good explication in the form of a mathematically strict and intuitively convincing semantics. 

 The paper is written in a style accessible to a wide range of readers. Some basic familiarity with computability logic, sequent calculus or linear logic is desirable 
only as much as to be able to duly appreciate the import of the present contribution. 

\end{abstract}


\  


\section{Introduction}\label{intr}

This paper introduces a refinement of the sequent calculus approach called {\em cirquent calculus}. Roughly speaking, the difference between the two is that, while in Gentzen-style proof trees sibling (or cousin, etc.) sequents are disjoint and independent sequences of formulas, in cirquent calculus they are permitted to share elements. Explicitly allowing or disallowing shared resources and thus taking to a more subtle level the resource-awareness intuitions underlying substructural logics, cirquent calculus offers much greater flexibility and power than sequent calculus does. 

A need for substantially new deductive tools came with the recent (2003) birth of {\em computability logic} (CL), characterized in  \cite{CL1} as ``a formal theory of computability in the same sense as classical logic is a formal theory of truth''. Indeed, formulas in CL are seen as computational problems rather than propositions or predicates, and their ``truth'' seen as algorithmic solvability. In turn, computational problems, understood in their most general --- interactive --- sense, are defined as games played by an interactive Turing machine against its environment, with ``algorithmic solvability'' meaning existence of a machine that wins the game against any possible (behavior of the) environment. A core collection of the most basic and natural operations 
on computational problems forms the logical vocabulary of the theory, with some of those operations, as logical operators, resembling those of linear logic. With this semantics, CL provides a systematic answer to the fundamental question ``{\em what can be computed?}\hspace{1pt}'', just as classical logic is a systematic tool for telling what is true. Furthermore, as it turns out, in positive cases ``{\em what} can be computed'' always allows itself to be replaced by ``{\em how} can be computed'', which makes CL of interest in not only theoretical computer science, but some more applied areas as well, including interactive knowledge base systems and resource oriented systems for planning and action. On the logical side, CL can serve as a basis for constructive   
applied theories. This is a very brief summary. See  
\cite{Jap03}, \cite{Japic} or \cite{Japfin} for elaborated expositions of the philosophy, motivations and techniques of computability logic.\footnote{A comprehensive online source on CL can be found at http://www.cis.upenn.edu/$\sim$giorgi/cl.html}

The above-mentioned fact of resemblance between computability-logic and linear-logic operators is no accident. Both logics claim to be ``logics of resources'', with their philosophies and ambitions thus having a significant overlap. The ways this common philosophy is materialized, however, are rather different. Computability logic directly captures resource intuitions through its semantics. Resources, understood in the specific sense of  {\em computational resources}, are dual/symmetric to computational problems: what is a problem for the machine, is a resource for the environment (=user), and vice versa. So, as a logic of computational problems, CL also automatically is a logic of computational resources. The scheme that CL follows can be characterized
as ``{\em from semantics to syntax}'': it starts with  a clear concept of resources (=computational problems) and resource-semantical validity (=universal algorithmic solvability), and only after that, as a natural second step, asks 
what the corresponding syntax is, i.e. how the set of valid formulas can be axiomatized. On the other hand, it would be accurate to say that linear logic, as a logic of resources (rather than that of phases or coherence spaces), has started directly from the second step, essentially by taking classical sequent calculus and deleting the structural rules unsound from a naive, purely intuitive resource point of view. For simplicity, in this discussion we narrow linear logic down to its multiplicative fragment; 
furthermore, taking some terminological liberty,  by ``linear logic'' we mean the version of it more commonly known as {\em affine logic}, which is classical sequent calculus without the contraction rule (Girard's \cite{Gir87} canonical
 system 
for linear logic further deletes the rule of weakening as well). Even the most naive and vague resource intuitions are sufficient to see that the deleted rule of contraction, responsible for the principle $P\mli P\mlc P$, was indeed wrong: having \$1 does not imply having \$1 {\em and} \$1, i.e. \$2. Such intuitions can also be safely relied upon in deeming all the other rules of classical sequent calculus ``right''.  To summarize, linear logic is undoubtedly sound as a logic of resources. But even more so is ... the empty logic. Completeness is thus a crucial issue. This is where the need for a mathematically strict and intuitively convincing resource semantics becomes critical, without which a question on completeness cannot even be meaningfully asked. 
Despite intensive efforts, however, such a semantics has never really been found for linear logic. And apparently the reason for this failure is as straightforward as it could possibly be: linear logic, as a resource logic, is simply incomplete. At least, this is what CL believes, for it has been shown (\cite{Japfin}) that the semantics of the latter, with its well-justified claims to be a semantics of resources, validates a strictly bigger class of formulas than linear (=affine) logic does.

Taking pride in the meaningfulness of its semantics, computability logic, at the same time, has been suffering from one apparent disadvantage: the absence of a good syntax, i.e. proof-theoretically ``reasonable'' and nice deductive systems, as opposed to the beauty and harmony of the Gentzen-style axiomatizations for linear logic and its variations, let alone proof nets. True, certain sound and complete systems, named {\bf CL1}, {\bf CL2}, {\bf CL3} and {\bf CL4}, have been constructed for incrementally expressive (and rather expressive) fragments of CL in \cite{CL1,CL2,CL3,CL4}, and probably more results in the same style are still to come.  Yet, hardly many would perceive those systems as ``logical calculi'', and perhaps not everyone would even call them ``deductive systems''. Rather, those somewhat bizarre constructions --- one of which ({\bf CL2}) will be reproduced later in Section \ref{slast} --- might be seen as just ad hoc syntactic characterizations, offering sound and complete decision or enumeration procedures for the corresponding sets of valid formulas of CL, but otherwise providing no real proof-theoretic insights into this new logic. Repeated attempts to find Gentzen- or Hilbert-style equivalents of those systems have hopelessly failed even at the most basic, $\gneg,\mlc,\mld,\mli$ (``multiplicative'') level. And probably this failure, just like the failure to find a good resource semantics for linear logic, is no accident. The traditional deductive methods have been originally developed with traditional logics in mind. There are no reasons to expect for those  
methods to be general and flexible enough to just as successfully accommodate the needs of finer-level semantic approaches, such as the computational semantics of CL, or resource semantics in general. Switching to a novel vision in semantics may require doing the same in syntax. 

This is where cirquent calculus as a nontraditional syntax comes in, breaking the stubborn resistance of CL to axiomatization attempts.
   While the full collection of its rules just offers an alternative axiomatization for 
the kind old classical logic,  removing (the cirquent-calculus version of) contraction from that collection --- we call the resulting system {\bf CL5} --- yields a sound and complete system for the $(\gneg,\mlc,\mld,\mli)$-fragment of CL, the very core of the logic previously appearing to be ``most unaxiomatizable''. Being complete, {\bf CL5} is thus strictly stronger than 
the incomplete affine logic. The latter, by merely deleting
the offending rule of contraction without otherwise trying to first appropriately re(de)fine ordinary sequent calculus, has thrown out the baby with the bath water. Among the innocent victims expelled together with contraction 
is {\em Blass's} \cite{Bla92} {\em principle} 

\[(P\mlc Q)\mld(R\mlc S)\mli (P\mld R)\mlc(Q\mld S),\]
provable in {\bf CL5} but not in affine logic, which, in fact, even fails to prove the less general formula  
\[(P\mlc P)\mld( P\mlc P)\mli (P\mld P)\mlc(P\mld P).\]

To strengthen the implied claim of computability logic that it is {\bf CL5} rather than affine logic that adequately materializes the resource philosophy traditionally associated with the latter,
the present paper further introduces an abstract resource semantics and shows that {\bf CL5} is sound and complete with respect to that semantics as well. Unlike the semantics of computability logic, which understands resources in the special --- computational --- sense, abstract resource semantics can be seen as a direct formalization of the more general intuitions in the style ``having \$1 does not imply having  \$1 and \$1'' or ``one cannot get both a candy and an apple for a dollar even if one dollar can buy either''. As noted earlier, the inherent incompleteness of linear logic, resulting from the fundamental limitations of the underlying sequent-calculus approach, is the reason why such intuitions
and examples, while so heavily relied on in the popular linear-logic literature, have never really found a good explication in the form of a mathematically well-defined semantics. 

The set of theorems of {\bf CL5} admits an alternative, simple yet non-deductive characterization, according to which this is the set of all binary tautologies and their substitutional instances.  Here {\em binary tautologies} mean tautologies of classical propositional logic in which no propositional letter occurs more than twice. The class of such formulas has naturally emerged in the past in several unrelated contexts. The earliest relevant piece  of literature  of which the author is aware is \cite{Jas63}, dating back to 1963, where Ja\'{s}kowski studied binary tautologies as the solution to the problem of characterizing the provable formulas of a certain deductive system. Andreas Blass came across the same class of formulas twice. In \cite{Bla92} he introduced a game semantics for linear-logic connectives and found that the multiplicative fragment of the corresponding logic was exactly the class of the substitutional instances of binary tautologies. In the same paper he argued that this class was inherently unaxiomatizable --- using his words, ``entirely foreign to proof theory''. Such an assessment was both right and wrong, depending on whether proof theory is understood in the strictly traditional (sequent calculus) or a more generous (cirquent calculus) sense. 11 years later, in \cite{Bla03}, using Herbrand's Theorem, Blass introduced the concept of universal simple Herbrand validity, a natural sort of resource consciousness that makes sense in classical logic. Blass found in \cite{Bla03} that this  (non-game) semantics validates exactly the same class of propositional formulas as his unrelated game semantics for the multiplicative fragment of the language of linear logic does. While independently experimenting with various semantical approaches prior to the invention of computability logic, the author of the present paper, too, had found game-semantical soundness and completeness of the class of binary tautologies and their substitutional instances. Once this happened in \cite{Jap97} and then, again, in \cite{Jap00,Jap02}. The underlying semantics in those two cases 
were rather different from each other, as well as different from that of CL or Blass's game semantics. The fact that 
the set of the theorems of {\bf CL5} arises in different approaches by different authors with various motivations and traditions, serves as additional empirical  evidence for the naturalness of {\bf CL5}. This is somewhat in the same sense as the existence of various models of computation that eventually yield the same class of computable functions speaks in favor of the Church-Turing Thesis. 

The version of cirquent calculus presented in this paper captures the most basic yet only a modest fraction of the otherwise very expressive language of computability logic. For instance, the
 formalism of the earlier-mentioned system {\bf CL4}, in addition to $\gneg,\mlc,\mld,\mli$ called {\em parallel connectives}, 
has the connectives $\adc,\add$ ({\em choice connectives}, resembling the additives of linear logic), 
 and the two groups $\ada,\ade$ ({\em choice}) and $\cla,\cle$ ({\em blind}) of quantifiers. Among the other operators officially introduced within the framework of CL so far are the 
{\em parallel} (``multiplicative'') {\em quantifiers}  $\mla,\mle$, and the two groups $\st,\cost$ ({\em branching}) and $\pst,\pcost$ ({\em parallel}) of {\em recurrence} (``exponential'') {\em operators}. Extending the cirquent-calculus approach so as to accommodate incrementally expressive fragments of CL is a task for the future. The results of the present paper could be seen just as first steps on that long road. What  is important is that the syntactic ice cover of computability logic, previously having seemed to be unbreakable, is now cracked. 
 
It should be noted that, even though computability logic has been at the center of discussion almost throughout this introductory section, essentially its relevance to the present paper is  limited to being the primary 
 source of motivation and inspiration. 
Cirquent calculus (the very idea of it), abstract resource semantics and all related technical results presented in this paper are new and, as the author wishes to hope, valuable in their own rights.\footnote{For the exception of the soundness and completeness theorem for {\bf CL5} with respect to the CL semantics, of course.} There is no overlap with any prior work on CL, and 
familiarity with the latter, while desirable, is not at all necessary for understanding the present material.

\section{Cirquents}\label{mul}

Throughout the rest of this paper, unless otherwise specified, by a {\bf formula} we mean one of the language of classical propositional logic. We consider the version of this language that has infinitely many non-logical {\bf atoms}
(also called {\bf propositional letters}), for which we use the metavariables $P,Q,R,S$, and no logical atoms such as $\twg$ or $\tlg$. The propositional connectives are limited to the unary $\gneg$ and the binary $\mlc,\mld$. If we write $F\mli G$, it is to be understood as an abbreviation of $\gneg F\mld G$. Furthermore, we officially allow $\gneg$ to be applied only to atoms. $\gneg\gneg F$ is to be understood as $F$, $\gneg (F\mlc G)$ as $\gneg F\mld \gneg G$, and $\gneg(F\mld G)$ as $\gneg F\mlc\gneg G$. 
  
Where $k$ is a natural number, by  a $k$-ary  {\bf pool} we mean a sequence $\seq{F_1,\ldots,F_k}$ of formulas. Such a sequence may have repetitions, and we refer 
to a  particular occurrence of a formula in a pool as an {\bf oformula}, with the prefix ``o'' derived from ``occurrence''. This prefix will as well be used with a similar meaning in other words and contexts where same objects --- such as, say, atoms or subformulas --- may have several occurrences.  Thus, the pool $\seq{F,G,F}$ has two formulas but three oformulas; and the formula $P\mld(P\mlc \gneg P)$ has one  
atom but three oatoms. For readability, we usually refer to oformulas (or oatoms, etc.) by the name of the corresponding formula (atom, etc.), as in the phrase ``the oformula $F$'', 
assuming that it is clear from the context which of the possibly many occurrences of $F$ we have in mind.

A $k$-ary {\bf cirquentstructure}, or simply {\bf structure}, is a finite sequence ${\bf St}\ =\ \langle\Gamma_1,$
$\ldots,\Gamma_m\rangle$ ($m\geq 0$), where each $\Gamma_i$, called a {\bf group} of ${\bf St}$, is a subset of $\{1,\ldots,k\}$. As in pools, here we may have $\Gamma_i=\Gamma_j$ for some $i\not=j$. Again, to differentiate between a particular occurrence of a group in the structure from the group as such, we use the term {\bf ogroup}. The structure 
$\seq{\{1,2\},\emptyset,\{1,3\},\{1,2\}}$ thus has three groups but four ogroups. Yet, as in the case of oformulas or oatoms, we may just say ``the ogroup $\{1,2\}$'' if it is clear  which of the two occurrences of the group $\{1,2\}$ is meant.  

\begin{definition}\label{first}
A $k$-ary ($k\geq 0$) 
{\bf cirquent}  
is a pair $C=(\str{C},\rsp{C})$, where $\str{C}$, called the {\bf structure} of $C$, is a $k$-ary cirquentstructure, and  $\rsp{C}$, called the {\bf pool} of $C$, is a $k$-ary pool. 

An (o)group of such a $C$ will mean an (o)group of $\str{C}$, and an (o)formula of $C$ mean an (o)formula of $\rsp{C}$. Also, we often prefer to think of the groups of $C$ as sets of its oformulas rather than sets of the corresponding ordinal numbers. For example, if $\rsp{C}=\seq{F,G,H}$ and $\Gamma=\{1,3\}$, we can understand the same $\Gamma$ as the set $\{F,H\}$ of oformulas. In this case we say that $\Gamma$ {\bf contains} $F$ and $H$. When $\Gamma$ is seen as an ogroup rather than a group, we may also refer to such a set $\{F,H\}$  as the {\bf content} of $\Gamma$.  

An $1$-ary cirquent whose only ogroup is $\{1\}$ is said to be a {\bf singleton}.
\end{definition}

We represent cirquents using {\bf diagrams}, such as the one shown below: 

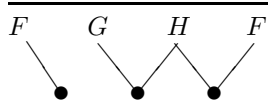
\begin{center} \begin{picture}(98,54)
\put(0,44){\line(1,0){98}}
\put(0,32){$F$}
\put(30,32){$G$}
\put(60,32){$H$}
\put(90,32){$F$}

\put(20,10){\line(-2,3){13}}
\put(49,10){\line(-4,5){15}}
\put(49,10){\line(4,5){15}}
\put(78,10){\line(-4,5){15}}
\put(78,10){\line(4,5){15}}
\put(20,10){\circle*{5}}
\put(49,10){\circle*{5}}
\put(78,10){\circle*{5}}
\end{picture}
\end{center}

This diagram represents the cirquent whose pool is $\seq{F,G,H,F}$ and whose structure is $\seq{\{1\},\{2,3\},\{3,4\}}$. We typically do not terminologically distinguish between cirquents and diagrams: for us, a diagram {\em is} (rather than {\em represents}) a cirquent, and a cirquent {\em is} a diagram. 
The top level of a diagram thus lists the oformulas of the cirquent, and the bottom level lists its ogroups, with each ogroup represented by (and identified with) a {\large $\bullet$}, where
the {\bf arcs} (lines connecting the {\large $\bullet$} 
with oformulas) are {\bf pointing to} the oformulas that a given ogroup contains. The horizontal line at the top of the diagram is just to indicate that this is one cirquent rather than, say, two cirquents (one $1$-ary and one $3$-ary) put together. Our convention is that such a line should be present even if there is no potential ambiguity. It is required to be long enough --- and OK if longer than necessary --- 
to cover all of the oformulas and ogroups of the cirquent. 

The term ``cirquent" is a hybrid of ``circuit" and ``sequent". So is, in a sense, its meaning. Cirquents can be seen to generalize sequents by imposing circuit-style structures on their oformulas. 
In a preliminary attempt to see some familiar meaning in cirquents, it might be helpful to think of them as Boolean circuits of depth 2, with oformulas serving as inputs, all first-level gates --- representing ogroups --- being $\mld$-gates, and the only second-level gate, connected to each first-level gate, being an $\mlc$-gate.
This is illustrated in Figure 1:
  
\begin{center} \begin{picture}(213,114)
\put(34,104){\em cirquent}
\put(0,92){\line(1,0){98}}
\put(0,80){$F$}
\put(30,80){$G$}
\put(60,80){$H$}
\put(90,80){$F$}

\put(20,58){\line(-2,3){13}}
\put(49,58){\line(-4,5){15}}
\put(49,58){\line(4,5){15}}
\put(78,58){\line(-4,5){15}}
\put(78,58){\line(4,5){15}}
\put(20,58){\circle*{5}}
\put(49,58){\circle*{5}}
\put(78,58){\circle*{5}}

\put(155,104){\em circuit}
\put(119,80){$F$}
\put(149,80){$G$}
\put(179,80){$H$}
\put(209,80){$F$}
\put(123,83){\circle{15}}
\put(153,83){\circle{15}}
\put(183,83){\circle{15}}
\put(213,83){\circle{15}}

\put(140,64){\line(-5,4){15}}
\put(169,64){\line(-5,4){15}}
\put(169,64){\line(5,4){15}}
\put(198,64){\line(-5,4){15}}
\put(198,64){\line(5,4){15}}
\put(140,56){\circle{15}}
\put(169,56){\circle{15}}
\put(198,56){\circle{15}}
\put(169,30){\circle{15}}
\put(137,53){$\mld$}
\put(166,53){$\mld$}
\put(195,53){$\mld$}
\put(166,27){$\mlc$}
\put(169,38){\line(0,1){10}}
\put(169,38){\line(5,2){27}}
\put(169,38){\line(-5,2){27}}
\put(90,10){\bf Figure 1}
\end{picture}
\end{center}

In traditional logic, circuits are interesting only in the context of representational complexity, and otherwise they do not offer any additional expressive power, for duplicating or merging identical nodes creates no difference when Boolean functions are all one sees in circuits. So, from the classical perspective,  the circuit of Figure 1 is equivalent to either circuit of Figure 2, with the tree-like circuit on the right being a direct reading of the formula $F\mlc(G\mld H)\mlc(H\mld F)$ expressing the Boolean function of the circuit of Figure 1, and the circuit on the left being a most economical representation of the same Boolean function: 

\begin{center} \begin{picture}(233,95)

\put(0,80){$G$}
\put(29,80){$F$}
\put(58,80){$H$}
\put(4,83){\circle{15}}
\put(33,83){\circle{15}}
\put(62,83){\circle{15}}

\put(4,64){\line(5,2){28}}
\put(33,64){\line(-5,2){28}}
\put(33,64){\line(5,2){29}}
\put(62,64){\line(-5,2){30}}
\put(62,64){\line(0,1){12}}
\put(4,56){\circle{15}}
\put(33,56){\circle{15}}
\put(62,56){\circle{15}}
\put(33,30){\circle{15}}
\put(1,53){$\mld$}
\put(30,53){$\mld$}
\put(59,53){$\mld$}
\put(30,27){$\mlc$}
\put(33,38){\line(0,1){10}}
\put(33,38){\line(5,2){27}}
\put(33,38){\line(-5,2){27}}
\put(96,10){\bf Figure 2}

\put(103,80){$F$}
\put(133,80){$G$}
\put(163,80){$H$}
\put(193,80){$H$}
\put(223,80){$F$}
\put(107,83){\circle{15}}
\put(137,83){\circle{15}}
\put(167,83){\circle{15}}
\put(197,83){\circle{15}}
\put(227,83){\circle{15}}

\put(124,64){\line(-5,4){15}}
\put(153,64){\line(-5,4){15}}
\put(153,64){\line(5,4){15}}
\put(212,64){\line(-5,4){15}}
\put(212,64){\line(5,4){15}}
\put(124,56){\circle{15}}
\put(153,56){\circle{15}}
\put(212,56){\circle{15}}
\put(166,30){\circle{15}}
\put(121,53){$\mld$}
\put(150,53){$\mld$}
\put(209,53){$\mld$}
\put(163,27){$\mlc$}
\put(166,38){\line(-1,1){11}}
\put(166,38){\line(4,1){43}}
\put(166,38){\line(-4,1){40}}
\end{picture}
\end{center}

Linear logic, understanding the nodes of the circuit as representing resources rather than 
just Boolean values, would not agree with such an equivalence though: the first and the fourth upper-level nodes of the circuit of Figure 1, even though having the same type, would be seen as two different individual resources. What linear logic generally fails to account for, however, is the possibility of resource sharing. $H$ is a resource shared by two different compound resources --- the resources represented by \#2 and \#3 $\mld$-gates of \mbox{Figure 1.} Allowing shared resources in cirquent calculus refines the otherwise crude approach of linear logic. And by no means does it mean departing from the idea that resources should  be accurately book-kept. Indeed, 
a shared resource does not mean a duplicated resource. Imagine Victor has \$10,000 on his bank account. One day he decides to give his wife access to the account. From now on the \$10,000 is shared. Two persons can 
use it, either at once, or portion by portion. Yet, this does not turn the \$10,000 into \$20,000, as the aggregate possible usage remains limited to \$10,000.

\section{Core cirquent calculus rules}\label{srules}

Different cirquent calculus systems will differ in what logical operators and atoms their underlying formal languages include, and what rules of inference they allow. The underlying language is fixed in this paper (it will only be slightly extended in the last paragraph of Section \ref{slast}). And all of the rules will come from the ones introduced in the present section. 
We explain those rules in a relaxed fashion, in terms of inserting arcs, swapping oformulas, etc. Such explanations are rather clear, and translating them into rigorous formulations in the style and terms of Definition \ref{first}, while possible, is hardly necessary.  

We need to agree on some additional terminology first. {\bf Adjacent oformulas} of a given cirquent are two oformulas $F$ and $G$ with $G$ appearing next to (to the right of) $F$ in the pool of the cirquent. We say that $F$ {\bf immediately precedes} $G$, and that $G$ {\bf immediately follows} $F$. Similarly for {\bf adjacent ogroups}. 

By {\bf merging two adjacent ogroups $\Gamma$ and $\Delta$} in a given cirquent $C$ we mean replacing in $C$ the two ogroups $\Gamma$ and $\Delta$ by the one ogroup $\Gamma\cup\Delta$, leaving the rest of the cirquent unchanged. The resulting cirquent will thus only differ from $C$ in that it will have one {\large $\bullet$} where $C$ had the two adjacent {\large $\bullet$}s,
with the arcs of this new {\large $\bullet$} pointing exactly to the oformulas to which the arcs of one or both of the old {\large $\bullet$}s were pointing. For example, the right cirquent of the following figure is the result of merging ogroups \#2 and \#3 in the left cirquent:

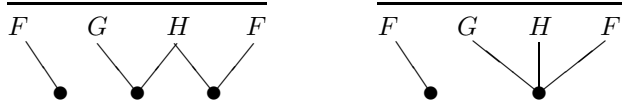
\begin{center} \begin{picture}(230,54)

\put(0,44){\line(1,0){98}}
\put(0,32){$F$}
\put(30,32){$G$}
\put(60,32){$H$}
\put(90,32){$F$}

\put(20,10){\line(-2,3){13}}
\put(49,10){\line(-4,5){15}}
\put(49,10){\line(4,5){15}}
\put(78,10){\line(-4,5){15}}
\put(78,10){\line(4,5){15}}
\put(20,10){\circle*{5}}
\put(49,10){\circle*{5}}
\put(78,10){\circle*{5}}

\put(140,44){\line(1,0){92}}
\put(140,32){$F$}
\put(170,32){$G$}
\put(197,32){$H$}
\put(224,32){$F$}

\put(160,10){\line(-2,3){13}}
\put(201,10){\line(-5,4){25}}
\put(201,10){\line(5,4){25}}
\put(201,10){\line(0,1){19}}
\put(160,10){\circle*{5}}
\put(201,10){\circle*{5}}
\end{picture}
\end{center}

{\bf Merging two adjacent oformulas $F$ and $G$ into $H$} means replacing those two oformulas by the one oformula $H$, and redirecting to it all arcs that were pointing to $F$ or $G$. For example,
the right cirquent of the following figure is the result of merging, in the left cirquent, (the first) $F$ and $G$ into $H$: 

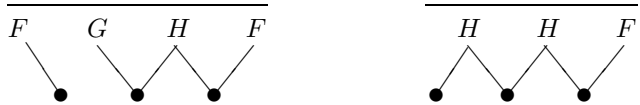
\begin{center} \begin{picture}(236,54)
\put(0,44){\line(1,0){98}}
\put(0,32){$F$}
\put(30,32){$G$}
\put(60,32){$H$}
\put(90,32){$F$}

\put(20,10){\line(-2,3){13}}
\put(49,10){\line(-4,5){15}}
\put(49,10){\line(4,5){15}}
\put(78,10){\line(-4,5){15}}
\put(78,10){\line(4,5){15}}
\put(20,10){\circle*{5}}
\put(49,10){\circle*{5}}
\put(78,10){\circle*{5}}

\put(158,44){\line(1,0){80}}
\put(170,32){$H$}
\put(200,32){$H$}
\put(230,32){$F$}
\put(162,10){\line(2,3){12}}
\put(189,10){\line(-4,5){15}}
\put(189,10){\line(4,5){15}}
\put(218,10){\line(-4,5){15}}
\put(218,10){\line(4,5){15}}
\put(162,10){\circle*{5}}
\put(189,10){\circle*{5}}
\put(218,10){\circle*{5}}

\end{picture}
\end{center}

Now we are ready to look at the rules. 

\subsection{Axioms ({\bf A})} Axioms are ``rules" with no premises. There are two axioms, called the {\bf empty cirquent axiom} and the {\bf identity axiom}. The first one introduces 
the {\bf empty cirquent} $(\seq{},\seq{})$ (both the pool and the structure are empty); the second one which --- just like the rest of the rules --- is, in fact, a scheme of rules because 
$F$ can be an arbitrary formula, 
introduces the cirquent $(\seq{\{1,2\}},\seq{\gneg F, F})$.

\begin{center} \begin{picture}(209,76)
\put(0,63){\em empty cirquent axiom}
\put(28,46){\line(1,0){36}}
\put(69,44){\scriptsize A}

\put(150,63){\em identity axiom}
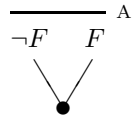
\put(160,46){\line(1,0){36}}
\put(200,44){\scriptsize A}
\put(160,33){$\gneg F$}
\put(188,33){$F$}
\put(180,10){\line(-3,5){11}}
\put(180,10){\line(3,5){11}}
\put(180,10){\circle*{5}}

\end{picture}
\end{center}

The letter ``A'' next to the horizontal line stands for the name of the rule by which the conclusion is obtained. We will follow the same notational practice for the other rules.

\subsection{Mix ({\bf M})} This rule takes two premises. The conclusion is obtained by simply putting one premise next to the other, thus creating one cirquent out of the two, as illustrated below: 

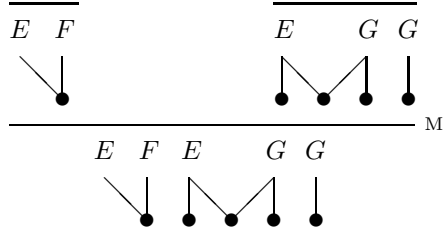
\begin{center} \begin{picture}(166,102)
\put(0,92){\line(1,0){26}}
\put(100,92){\line(1,0){54}}
\put(100,79){$E$}
\put(132,79){$G$}
\put(147,79){$G$}
\put(0,79){$E$}
\put(17,79){$F$}
\put(52,10){\circle*{5}}
\put(52,10){\line(0,1){16}}
\put(52,10){\line(-1,1){16}}
\put(68,10){\circle*{5}}
\put(84,10){\circle*{5}}
\put(100,10){\circle*{5}}
\put(116,10){\circle*{5}}
\put(68,10){\line(0,1){16}}
\put(84,10){\line(-1,1){16}}
\put(84,10){\line(1,1){16}}
\put(100,10){\line(0,1){16}}
\put(116,10){\line(0,1){16}}

\put(0,46){\line(1,0){153}}
\put(157,44){\scriptsize M}
\put(65,33){$E$}
\put(97,33){$G$}
\put(112,33){$G$}
\put(32,33){$E$}
\put(49,33){$F$}
\put(20,56){\circle*{5}}
\put(20,56){\line(0,1){16}}
\put(20,56){\line(-1,1){16}}
\put(103,56){\circle*{5}}
\put(119,56){\circle*{5}}
\put(135,56){\circle*{5}}
\put(151,56){\circle*{5}}
\put(103,56){\line(0,1){16}}
\put(119,56){\line(-1,1){16}}
\put(119,56){\line(1,1){16}}
\put(135,56){\line(0,1){16}}
\put(151,56){\line(0,1){16}}
\end{picture}
\end{center}

\subsection{Exchange (E)} This and all of the remaining rules take a single premise. The exchange rule comes in two flavors: {\bf oformula exchange} and {\bf ogroup exchange}. The conclusion of  oformula (resp. ogroup) exchange is the result of swapping in the premise  
two adjacent oformulas (resp. ogroups) and correspondingly redirecting all arcs. The following oformula exchange example swaps $F$ with $G$; and the ogroup exchange example swaps ogroup \#2 with ogroup \#3: 

\begin{center} \begin{picture}(182,124)
\put(-4,109){\em oformula exchange}
\put(10,92){\line(1,0){53}}
\put(10,79){$F$}
\put(33,79){$G$}
\put(55,79){$H$}
\put(14,56){\line(0,1){19}}
\put(14,56){\line(5,4){23}}
\put(14,56){\circle*{5}}
\put(37,56){\circle*{5}}
\put(59,56){\circle*{5}}
\put(37,56){\line(0,1){18}}
\put(37,56){\line(5,4){22}}
\put(59,56){\line(0,1){18}}

\put(10,46){\line(1,0){53}}
\put(67,44){\scriptsize E}
\put(10,33){$G$}
\put(33,33){$F$}
\put(55,33){$H$}
\put(14,10){\line(0,1){18}}
\put(14,10){\line(5,4){23}}
\put(14,10){\circle*{5}}
\put(37,10){\circle*{5}}
\put(59,10){\circle*{5}}
\put(37,10){\line(-5,4){22}}
\put(37,10){\line(5,4){22}}
\put(59,10){\line(0,1){18}}

\put(122,109){\em ogroup exchange}
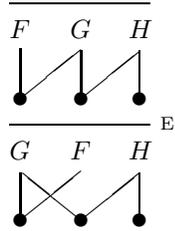
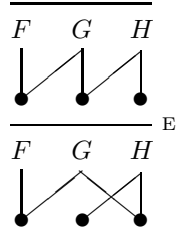
\put(130,92){\line(1,0){53}}
\put(130,79){$F$}
\put(153,79){$G$}
\put(175,79){$H$}
\put(134,56){\line(0,1){19}}
\put(134,56){\line(5,4){23}}
\put(134,56){\circle*{5}}
\put(157,56){\circle*{5}}
\put(179,56){\circle*{5}}
\put(157,56){\line(0,1){19}}
\put(157,56){\line(5,4){22}}
\put(179,56){\line(0,1){18}}

\put(130,46){\line(1,0){53}}
\put(187,44){\scriptsize E}
\put(130,33){$F$}
\put(153,33){$G$}
\put(175,33){$H$}
\put(134,10){\line(0,1){19}}
\put(134,10){\line(5,4){23}}
\put(134,10){\circle*{5}}
\put(157,10){\circle*{5}}
\put(179,10){\circle*{5}}
\put(179,10){\line(-5,4){23}}
\put(157,10){\line(5,4){22}}
\put(179,10){\line(0,1){18}}
\end{picture}
\end{center}

The presence of oformula exchange essentially allows us to treat the pool of a cirquent as a multiset rather than a sequence of formulas. Similarly, the presence of ogroup exchange makes it possible to see the structure of a cirquent as a multiset rather than a sequence of groups.  

\subsection{Weakening (W)} This rule, too, comes in two flavors: {\bf ogroup weakening} and {\bf pool weakening}. In the first case the conclusion is the result of inserting a new arc between an existing ogroup and an existing oformula of the premise. In the second case, the  conclusion is the result of inserting a new oformula anywhere in the pool of the premise. 

\begin{center} \begin{picture}(195,124)
\put(0,109){\em ogroup weakening}
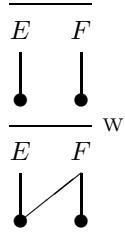
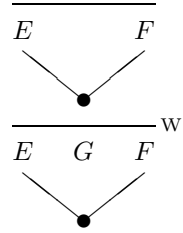
\put(20,92){\line(1,0){31}}
\put(20,79){$E$}
\put(43,79){$F$}
\put(24,56){\line(0,1){18}}
\put(24,56){\circle*{5}}
\put(47,56){\circle*{5}}
\put(47,56){\line(0,1){18}}

\put(20,46){\line(1,0){31}}
\put(55,44){\scriptsize W}
\put(20,33){$E$}
\put(43,33){$F$}
\put(24,10){\line(0,1){18}}
\put(24,10){\line(5,4){23}}
\put(24,10){\circle*{5}}
\put(47,10){\circle*{5}}
\put(47,10){\line(0,1){18}}

\put(140,109){\em pool weakening}
\put(142,92){\line(1,0){54}}
\put(142,79){$E$}
\put(188,79){$F$}
\put(169,56){\circle*{5}}
\put(169,56){\line(-5,4){23}}
\put(169,56){\line(5,4){23}}

\put(142,46){\line(1,0){54}}
\put(198,44){\scriptsize W}
\put(142,33){$E$}
\put(165,33){$G$}
\put(188,33){$F$}
\put(169,10){\line(-5,4){23}}
\put(169,10){\line(5,4){23}}
\put(169,10){\circle*{5}}

\end{picture}
\end{center}

\subsection{Duplication (D)} This rule comes in two versions as well: {\bf downward duplication} and {\bf upward duplication}. The conclusion (resp. premise) of downward (resp. upward) duplication is the result of replacing in the premise (resp. conclusion) some ogroup $\Gamma$ by two adjacent ogroups that, as groups, are identical with 
$\Gamma$.

\begin{center} \begin{picture}(238,122)
\put(0,107){\em downward duplication}
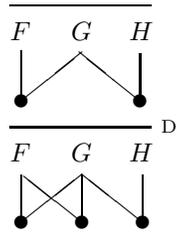
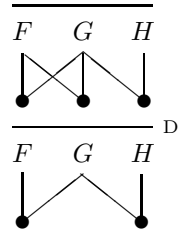
\put(19,92){\line(1,0){53}}
\put(19,79){$F$}
\put(42,79){$G$}
\put(64,79){$H$}
\put(23,56){\line(0,1){19}}
\put(23,56){\line(5,4){23}}
\put(23,56){\circle*{5}}
\put(68,56){\circle*{5}}
\put(68,56){\line(0,1){18}}
\put(68,56){\line(-5,4){23}}

\put(19,46){\line(1,0){53}}
\put(76,44){\scriptsize D}
\put(19,33){$F$}
\put(42,33){$G$}
\put(64,33){$H$}
\put(23,10){\line(0,1){18}}
\put(23,10){\line(5,4){23}}
\put(23,10){\circle*{5}}
\put(46,10){\circle*{5}}
\put(69,10){\circle*{5}}
\put(46,10){\line(-5,4){22}}
\put(46,10){\line(0,1){18}}
\put(69,10){\line(0,1){18}}
\put(69,10){\line(-5,4){23}}

\put(156,107){\em upward duplication}
\put(169,46){\line(1,0){53}}
\put(169,33){$F$}
\put(192,33){$G$}
\put(214,33){$H$}
\put(173,10){\line(0,1){19}}
\put(173,10){\line(5,4){23}}
\put(173,10){\circle*{5}}
\put(218,10){\circle*{5}}
\put(218,10){\line(0,1){18}}
\put(218,10){\line(-5,4){23}}

\put(169,92){\line(1,0){53}}
\put(226,44){\scriptsize D}
\put(169,79){$F$}
\put(192,79){$G$}
\put(214,79){$H$}
\put(173,56){\line(0,1){18}}
\put(173,56){\line(5,4){23}}
\put(173,56){\circle*{5}}
\put(196,56){\circle*{5}}
\put(219,56){\circle*{5}}
\put(196,56){\line(-5,4){22}}
\put(196,56){\line(0,1){18}}
\put(219,56){\line(0,1){18}}
\put(219,56){\line(-5,4){23}}

\end{picture}
\end{center}

Note that the presence of duplication together with ogroup exchange further allows us to think of the structure of a cirquent as a set rather than a sequence or multiset of groups.

\subsection{Contraction (C)} The premise of this rule is a cirquent with two adjacent oformulas $F,F$ that are identical as formulas. The conclusion is obtained from the premise by merging those two oformulas into $F$. 
The following two examples illustrate applications of contraction. 

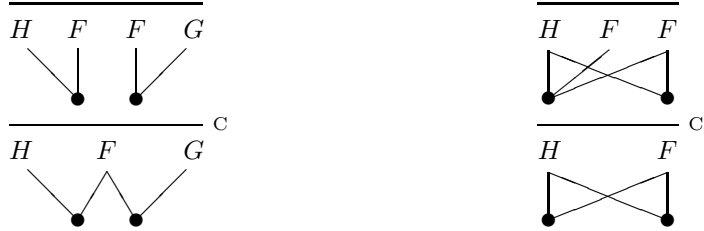
\begin{center} \begin{picture}(266,102)
\put(0,92){\line(1,0){73}}
\put(0,79){$H$}
\put(22,79){$F$}
\put(44,79){$F$}
\put(66,79){$G$}
\put(26,56){\line(0,1){19}}
\put(26,56){\circle*{5}}
\put(26,56){\line(-1,1){19}}
\put(48,56){\line(1,1){19}}
\put(48,56){\circle*{5}}
\put(48,56){\line(0,1){19}}

\put(0,46){\line(1,0){73}}
\put(77,44){\scriptsize C}
\put(0,33){$H$}
\put(33,33){$F$}
\put(66,33){$G$}
\put(26,10){\circle*{5}}
\put(48,10){\circle*{5}}
\put(26,10){\line(-1,1){19}}
\put(48,10){\line(1,1){19}}
\put(26,10){\line(3,5){11}}
\put(48,10){\line(-3,5){11}}

\put(200,92){\line(1,0){53}}
\put(200,79){$H$}
\put(223,79){$F$}
\put(245,79){$F$}
\put(204,56){\line(0,1){18}}
\put(204,56){\line(5,4){23}}
\put(204,56){\line(5,2){45}}
\put(249,56){\line(-5,2){45}}
\put(204,56){\circle*{5}}
\put(249,56){\circle*{5}}
\put(249,56){\line(0,1){18}}

\put(200,46){\line(1,0){53}}
\put(257,44){\scriptsize C}
\put(200,33){$H$}
\put(245,33){$F$}
\put(204,10){\line(0,1){18}}
\put(204,10){\line(5,2){45}}
\put(249,10){\line(-5,2){45}}
\put(204,10){\circle*{5}}
\put(249,10){\circle*{5}}
\put(249,10){\line(0,1){18}}
\end{picture}
\end{center}

\subsection{$\mld$-introduction ({\bf $\mld$})}  
The conclusion of this rule is obtained by merging in the premise some two adjacent oformulas $F$ and $G$ into $F\mld G$. We say that this application of the rule {\bf introduces} $F\mld G$. Below are three illustrations:

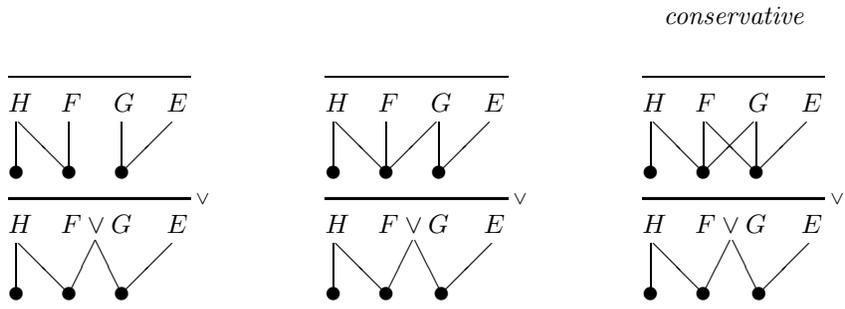
\begin{center} \begin{picture}(316,132)

\put(0,92){\line(1,0){69}}
\put(0,79){$H$}
\put(20,79){$F$}
\put(40,79){$G$}
\put(60,79){$E$}
\put(3,56){\circle*{5}}
\put(23,56){\circle*{5}}
\put(43,56){\circle*{5}}
\put(3,56){\line(0,1){19}}
\put(23,56){\line(0,1){19}}
\put(23,56){\line(-1,1){19}}
\put(43,56){\line(0,1){19}}
\put(43,56){\line(1,1){19}}

\put(0,46){\line(1,0){69}}
\put(71,44){\scriptsize $\mld$}

\put(0,33){$H$}
\put(20,33){$F\mld G$}
\put(60,33){$E$}
\put(3,10){\circle*{5}}
\put(23,10){\circle*{5}}
\put(43,10){\circle*{5}}
\put(3,10){\line(0,1){19}}
\put(23,10){\line(1,2){10}}
\put(23,10){\line(-1,1){19}}
\put(43,10){\line(-1,2){10}}
\put(43,10){\line(1,1){19}}

\put(120,92){\line(1,0){69}}
\put(120,79){$H$}
\put(140,79){$F$}
\put(160,79){$G$}
\put(180,79){$E$}
\put(123,56){\circle*{5}}
\put(143,56){\circle*{5}}
\put(163,56){\circle*{5}}
\put(123,56){\line(0,1){19}}
\put(143,56){\line(0,1){19}}
\put(143,56){\line(-1,1){19}}
\put(143,56){\line(1,1){19}}
\put(163,56){\line(0,1){19}}
\put(163,56){\line(1,1){19}}

\put(120,46){\line(1,0){69}}
\put(191,44){\scriptsize $\mld$}
\put(120,33){$H$}
\put(140,33){$F\mld G$}
\put(180,33){$E$}
\put(123,10){\circle*{5}}
\put(143,10){\circle*{5}}
\put(164,10){\circle*{5}}
\put(123,10){\line(0,1){19}}
\put(143,10){\line(1,2){10}}
\put(143,10){\line(-1,1){19}}
\put(164,10){\line(-1,2){10}}
\put(164,10){\line(1,1){19}}

\put(248,112){\em conservative}
\put(240,92){\line(1,0){69}}
\put(240,79){$H$}
\put(260,79){$F$}
\put(280,79){$G$}
\put(300,79){$E$}
\put(243,56){\circle*{5}}
\put(263,56){\circle*{5}}
\put(283,56){\circle*{5}}
\put(243,56){\line(0,1){19}}
\put(263,56){\line(0,1){19}}
\put(263,56){\line(-1,1){19}}
\put(263,56){\line(1,1){19}}
\put(283,56){\line(0,1){19}}
\put(283,56){\line(1,1){19}}
\put(283,56){\line(-1,1){19}}

\put(240,46){\line(1,0){69}}
\put(311,44){\scriptsize $\mld$}
\put(240,33){$H$}
\put(260,33){$F\mld G$}
\put(300,33){$E$}
\put(243,10){\circle*{5}}
\put(263,10){\circle*{5}}
\put(284,10){\circle*{5}}
\put(243,10){\line(0,1){19}}
\put(263,10){\line(1,2){10}}
\put(263,10){\line(-1,1){19}}
\put(284,10){\line(-1,2){10}}
\put(284,10){\line(1,1){19}}

\end{picture}
\end{center}

In what we call {\bf conservative $\mld$-introduction} (the rightmost example), which is a special case of $\mld$-introduction, the situation is that whenever an ogroup of the conclusion contains the introduced $F\mld G$, the corresponding ogroup of the premise contains {\em both} $F$ and $G$. In a general case (the first two examples), this is not necessary. What {\em is} always necessary, however, is that if an ogroup of the conclusion contains the introduced $F\mld G$, then the corresponding ogroup of the premise should contain {\em at least one} of the oformulas $F,G$.

We have just used and will continue to use the jargon ``the corresponding ogroup", whose meaning should be clear: the present rule does not change the number or order of ogroups, and it only modifies the contents of some of those ogroups. So, to ogroup \#$i$ of the conclusion {\bf corresponds} ogroup \#$i$ of the premise, and vice versa. The same applies to the rules of oformula exchange, weakening and contraction. In an application of ogroup exchange that swaps ogroups \#$i$ and \#$i+1$, to ogroup \#$i$ of the premise corresponds ogroup \#$i+1$ of the conclusion, and vice versa; to ogroup \#$i+1$ of the premise corresponds ogroup \#$i$ of the conclusion, and vice versa; and to any other ogroup \#$j$ of the premise corresponds ogroup \#$j$ of the conclusion and vice versa. Finally, in an application of mix, to ogroup \#$i$ of the first premise corresponds ogroup \#$i$ of the conclusion, and vice versa; and, where $n$ is the 
number of the ogroups of the first premise, to ogroup \#$i$ of the second premise corresponds ogroup \#$n+i$ of the conclusion, and vice versa.

\subsection{$\mlc$-introduction ($\mlc$)}
 
The premise of this rule is a cirquent with adjacent oformulas $F$ and $G$, such that the following two conditions are satisfied:
\begin{itemize}
\item No ogroup contains both $F$ and $G$.
\item Every ogroup containing $F$ is immediately followed by an ogroup containing $G$, and every ogroup containing $G$ is immediately preceeded by an ogroup containing $F$.  
\end{itemize}

The conclusion is obtained from the premise by merging each ogroup containing $F$ with the immediately following ogroup (containing $G$) and then, in the resulting cirquent, merging $F$ and $G$ into $F\mlc G$. In this case we say that the  the rule {\bf introduces} $F\mlc G$. 

Below are three examples for the simple case when there is only one ogroup in the conclusion that contains the introduced $F\mlc G$:

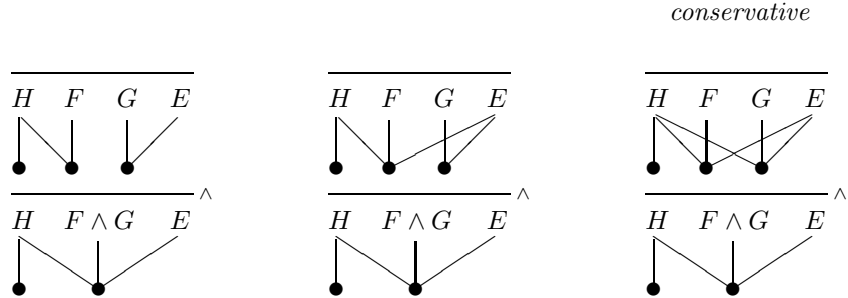
\begin{center} \begin{picture}(316,132)
\put(0,92){\line(1,0){69}}
\put(0,79){$H$}
\put(20,79){$F$}
\put(40,79){$G$}
\put(60,79){$E$}
\put(3,56){\circle*{5}}
\put(3,56){\line(0,1){19}}
\put(23,56){\circle*{5}}
\put(23,56){\line(0,1){18}}
\put(23,56){\line(-1,1){19}}
\put(44,56){\circle*{5}}
\put(44,56){\line(0,1){18}}
\put(44,56){\line(1,1){19}}

\put(0,46){\line(1,0){69}}
\put(71,44){\scriptsize $\mlc$}

\put(0,33){$H$}
\put(20,33){$F\mlc G$}
\put(60,33){$E$}
\put(3,10){\circle*{5}}
\put(3,10){\line(0,1){19}}
\put(33,10){\circle*{5}}
\put(33,10){\line(0,1){18}}
\put(33,10){\line(-3,2){30}}
\put(33,10){\line(3,2){30}}

\put(249,112){\em conservative}
\put(240,92){\line(1,0){69}}
\put(240,79){$H$}
\put(260,79){$F$}
\put(280,79){$G$}
\put(300,79){$E$}
\put(243,56){\circle*{5}}
\put(243,56){\line(0,1){19}}
\put(263,56){\circle*{5}}
\put(263,56){\line(0,1){18}}
\put(263,56){\line(-1,1){19}}
\put(263,56){\line(2,1){40}}
\put(284,56){\line(-2,1){40}}
\put(284,56){\circle*{5}}
\put(284,56){\line(0,1){18}}
\put(284,56){\line(1,1){19}}
\put(240,46){\line(1,0){69}}
\put(311,44){\scriptsize $\mlc$}
\put(240,33){$H$}
\put(260,33){$F\mlc G$}
\put(300,33){$E$}
\put(243,10){\circle*{5}}
\put(243,10){\line(0,1){19}}
\put(273,10){\circle*{5}}
\put(273,10){\line(0,1){18}}
\put(273,10){\line(-3,2){30}}
\put(273,10){\line(3,2){30}}

\put(120,92){\line(1,0){69}}
\put(120,79){$H$}
\put(120,33){$H$}
\put(120,46){\line(1,0){69}}
\put(140,79){$F$}
\put(160,79){$G$}
\put(180,79){$E$}
\put(180,33){$E$}
\put(140,33){$F\mlc G$}
\put(123,56){\circle*{5}}
\put(123,56){\line(0,1){19}}
\put(143,56){\circle*{5}}
\put(143,56){\line(0,1){18}}
\put(143,56){\line(-1,1){19}}
\put(143,56){\line(2,1){40}}
\put(164,56){\line(1,1){19}}
\put(164,56){\circle*{5}}
\put(164,56){\line(0,1){18}}
\put(191,44){\scriptsize $\mlc$}
\put(123,10){\circle*{5}}
\put(123,10){\line(0,1){19}}
\put(153,10){\circle*{5}}
\put(153,10){\line(0,1){18}}
\put(153,10){\line(-3,2){30}}
\put(153,10){\line(3,2){30}}
\end{picture}
\end{center}

Perhaps this rule is easier to comprehend in the bottom-up (from conclusion to premise) view. To obtain a premise from the conclusion (where $F\mlc G$ is the introduced conjunction), we  
``split" every ogroup $\Gamma$ containing $F\mlc G$ into two adjacent ogroups $\Gamma^F$ and $\Gamma^G$, where $\Gamma^F$ contains $F$ (but not $G$), and $\Gamma^G$ contains $G$ (but not $F$); all other ($\not= F\mlc G$) oformulas of $\Gamma$ --- and only such oformulas --- should be included in either $\Gamma^F$, or $\Gamma^G$, or both. In what we call {\bf conservative $\mlc$-introduction}, {\em all} of the non-$F\mlc G$ oformulas of $\Gamma$ should be included in {\em both} $\Gamma^F$ and $\Gamma^G$.

The following is an example of an application of the $\mlc$-introduction rule in a little bit more complex case where the conclusion has two ogroups containing the introduced conjunction. It is not a conservative one. To make this application conservative, we should add two more arcs to the premise: one connecting ogroup \#3 with $J$, and one connecting ogroup \#4 with $E$.

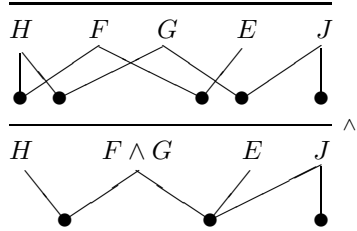
\begin{center} \begin{picture}(107,102)
\put(0,92){\line(1,0){122}}
\put(0,79){$H$}
\put(30,79){$F$}
\put(56,79){$G$}
\put(86,79){$E$}
\put(116,79){$J$}
\put(19,56){\line(-4,5){14}}
\put(19,56){\line(2,1){39}}
\put(19,56){\circle*{5}}
\put(4,56){\circle*{5}}
\put(4,56){\line(0,1){17}}
\put(4,56){\line(3,2){30}}
\put(88,56){\line(-3,2){30}}
\put(73,56){\line(4,5){15}}
\put(73,56){\line(-2,1){39}}
\put(88,56){\line(3,2){30}}
\put(73,56){\circle*{5}}
\put(88,56){\circle*{5}}
\put(118,56){\circle*{5}}
\put(118,56){\line(0,1){19}}

\put(0,46){\line(1,0){122}}
\put(126,44){\scriptsize $\mlc$}

\put(0,33){$H$}
\put(35,33){$F\mlc G$}
\put(88,33){$E$}
\put(115,33){$J$}
\put(21,10){\line(-4,5){15}}
\put(21,10){\line(3,2){28}}
\put(21,10){\circle*{5}}
\put(76,10){\line(-3,2){28}}
\put(76,10){\line(4,5){15}}
\put(76,10){\line(2,1){42}}
\put(76,10){\circle*{5}}
\put(118,10){\circle*{5}}
\put(118,10){\line(0,1){20}}

\end{picture}
\end{center}

\section{Cirquent calculus systems}\label{sccs}

By a {\bf cirquent calculus system} in the present context we mean any subset of the set of the eight rules of the previous section. The one that has the full collection of all eight rules we denote by {\bf CCC} (``Classical Cirquent Calculus''), and the one that has all rules but contraction we denote by {\bf CL5}. Any other system we denote by placing the abbreviated names of the corresponding rules between parentheses.  For instance, (AME) stands for the system that has the axioms, mix and exchange.

Let $S$ be a cirquent calculus system, and $C,A_1,\ldots,A_n$ (possibly $n=0$) any cirquents. A {\bf derivation}  of $C$ from $A_1,\ldots,A_n$ in $S$ is a tree of cirquents with $C$ at its root, 
where each node is a cirquent that either follows from its children by one of the rules of $S$, or is among  
$A_1,\ldots,A_n$ (and has no children). A derivation of $C$ in $S$ from the empty set of cirquents is said to be a {\bf proof} of $C$ in $S$. Of course, if $S$ does not contain axioms, 
then there will be no proofs in it.

Throughout this paper we identify each formula $F$ with the singleton cirquent $(\seq{\{1\}},\seq{F})$, i.e. the cirquent 

\begin{center} \begin{picture}(6,43)
\put(0,33){$F$}
\put(4,10){\line(0,1){19}}
\put(4,10){\circle*{5}}
\end{picture}
\end{center}

Correspondingly, a  proof or derivation of a given formula $F$ in a given system $S$ is a proof or derivation of $(\seq{\{1\}},\seq{F})$. 
The following is an example of a proof of $\gneg F\mld (F\mlc F)$ in (AMEC$\mld\mlc$):

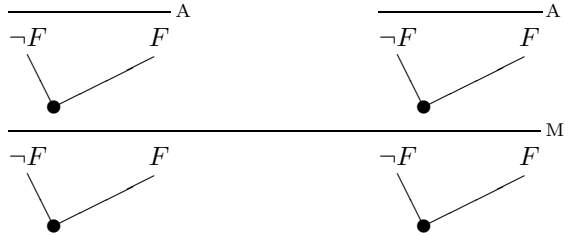
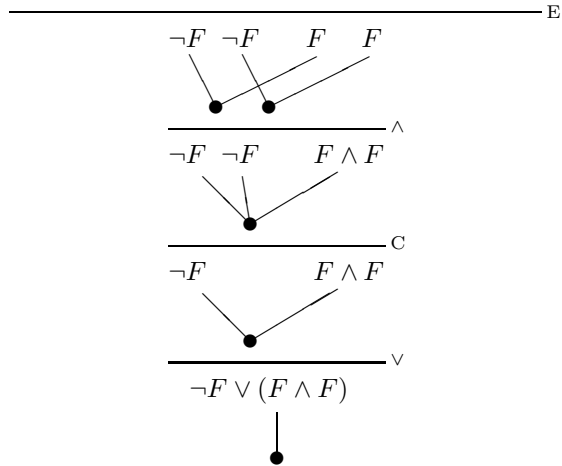
\begin{center}\begin{picture}(199,66)
\put(0,46){\line(1,0){61}}
\put(63,44){\scriptsize A}
\put(0,33){$\gneg F$}
\put(53,33){$F$}
\put(17,10){\line(-1,2){10}}
\put(17,10){\line(2,1){38}}
\put(17,10){\circle*{5}}

\put(140,46){\line(1,0){61}}
\put(203,44){\scriptsize A}
\put(184,46){\line(1,0){18}}
\put(140,33){$\gneg F$}
\put(193,33){$F$}
\put(157,10){\line(-1,2){10}}
\put(157,10){\line(2,1){38}}
\put(157,10){\circle*{5}}
\end{picture}
\end{center}

\begin{center}\begin{picture}(199,34)
\put(0,33){$\gneg F$}
\put(53,33){$F$}
\put(17,10){\line(-1,2){10}}
\put(17,10){\line(2,1){38}}
\put(17,10){\circle*{5}}

\put(0,46){\line(1,0){201}}
\put(203,44){\scriptsize M}
\put(140,33){$\gneg F$}
\put(193,33){$F$}
\put(157,10){\line(-1,2){10}}
\put(157,10){\line(2,1){38}}
\put(157,10){\circle*{5}}
\end{picture}
\end{center}

\begin{center}\begin{picture}(199,34)
\put(0,46){\line(1,0){201}}
\put(203,44){\scriptsize E}
\put(60,33){$\gneg F$}
\put(80,33){$\gneg F$}
\put(112,33){$F$}
\put(133,33){$F$}
\put(98,10){\line(-1,2){10}}
\put(98,10){\circle*{5}}
\put(78,10){\circle*{5}}
\put(78,10){\line(-1,2){10}}
\put(78,10){\line(2,1){38}}
\put(98,10){\line(2,1){38}}

\end{picture}
\end{center}

\begin{center}\begin{picture}(199,34)
\put(60,46){\line(1,0){82}}
\put(144,45){\scriptsize $\mlc$}
\put(60,33){$\gneg F$}
\put(80,33){$\gneg F$}
\put(115,33){$F\mlc F$}
\put(91,10){\line(-1,1){18}}
\put(91,10){\line(-1,6){3}}
\put(91,10){\line(5,3){33}}
\put(91,10){\circle*{5}}
\end{picture}
\end{center}

\begin{center}\begin{picture}(199,34)
\put(60,46){\line(1,0){82}}
\put(144,45){\scriptsize C}
\put(60,33){$\gneg F$}
\put(115,33){$F\mlc F$}
\put(91,10){\line(-1,1){18}}
\put(91,10){\line(5,3){33}}
\put(91,10){\circle*{5}}
\end{picture}
\end{center}

\begin{center}\begin{picture}(199,34)
\put(60,46){\line(1,0){82}}
\put(144,45){\scriptsize $\mld$}
\put(68,33){$\gneg F\mld (F\mlc F)$}
\put(101,10){\line(0,1){17}}
\put(101,10){\circle*{5}}
\end{picture}
\end{center}

It is our convention that if a proof is a proof of a formula $F$, then the last cirquent we simply represent as ``$F$" rather than through a diagram. Just to save space. In a similar space-saving spirit, we will often combine several obvious steps together, labeling the combined application of a ``rule" by the name of the system which contains all of the rules that have been combined. For instance, the above derivation of $\gneg F\mld (F\mlc F)$  we might want to rewrite in a more compact  yet clear way as follows:
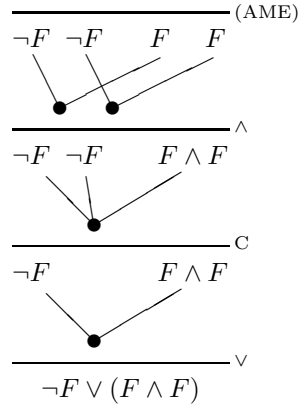
\vspace{10pt}

\begin{center}\begin{picture}(212,56)
\put(60,46){\line(1,0){82}}
\put(144,44){\scriptsize (AME)}
\put(60,33){$\gneg F$}
\put(80,33){$\gneg F$}
\put(112,33){$F$}
\put(133,33){$F$}
\put(98,10){\line(-1,2){10}}
\put(98,10){\circle*{5}}
\put(78,10){\circle*{5}}
\put(78,10){\line(-1,2){10}}
\put(78,10){\line(2,1){38}}
\put(98,10){\line(2,1){38}}
\end{picture}
\end{center}

\begin{center}\begin{picture}(212,34)
\put(60,46){\line(1,0){82}}
\put(144,45){\scriptsize $\mlc$}
\put(60,33){$\gneg F$}
\put(80,33){$\gneg F$}
\put(115,33){$F\mlc F$}
\put(91,10){\line(-1,1){18}}
\put(91,10){\line(-1,6){3}}
\put(91,10){\line(5,3){33}}
\put(91,10){\circle*{5}}
\end{picture}
\end{center}

\begin{center}\begin{picture}(212,34)
\put(60,46){\line(1,0){82}}
\put(144,45){\scriptsize C}
\put(60,33){$\gneg F$}
\put(115,33){$F\mlc F$}
\put(91,10){\line(-1,1){18}}
\put(91,10){\line(5,3){33}}
\put(91,10){\circle*{5}}
\end{picture}
\end{center}

\begin{center}\begin{picture}(212,11)
\put(60,23){\line(1,0){82}}
\put(144,22){\scriptsize $\mld$}
\put(71,10){$\gneg F\mld (F\mlc F)$}
\end{picture}
\end{center}
\ \vspace{-5pt}

Below is an (AME$\mld\mlc$)-proof of {\bf Blass's principle}  mentioned in Section \ref{intr}:
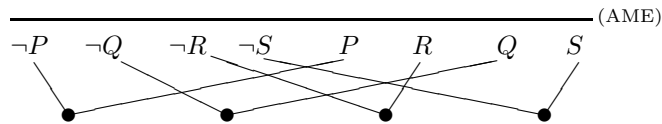
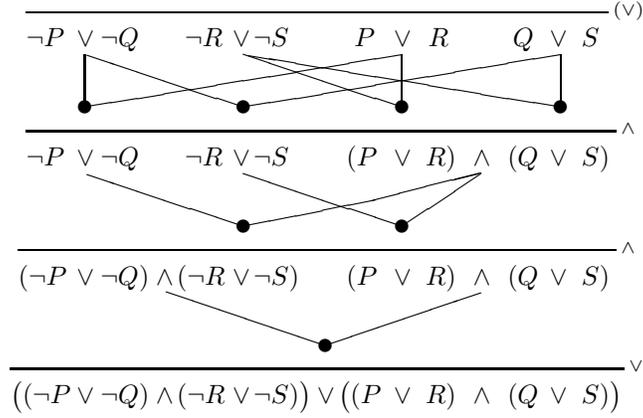
\vspace{20pt}

\begin{center}\begin{picture}(228,44)
\put(0,46){\line(1,0){220}}
\put(222,45){\scriptsize (AME)}
\put(0,33){$\gneg P$}
\put(28,33){$\gneg Q$}
\put(22,10){\line(-2,3){13}}
\put(22,10){\line(5,1){104}}
\put(22,10){\circle*{5}}

\put(60,33){$\gneg R$}
\put(86,33){$\gneg S$}
\put(82,10){\line(-2,1){40}}
\put(82,10){\line(5,1){102}}
\put(82,10){\circle*{5}}

\put(124,33){$P$}
\put(152,33){$R$}
\put(142,10){\line(-3,1){66}}
\put(142,10){\line(2,3){13}}
\put(142,10){\circle*{5}}

\put(184,33){$Q$}
\put(210,33){$S$}
\put(202,10){\line(-5,1){106}}
\put(202,10){\line(2,3){13}}
\put(202,10){\circle*{5}}
\end{picture}
\end{center}

\begin{center}\begin{picture}(228,34)
\put(0,46){\line(1,0){220}}
\put(222,45){\scriptsize ($\mld$)}
\put(0,33){$\gneg P$}
\put(28,33){$\gneg Q$}
\put(19,33){$\mld$}
\put(22,10){\line(0,1){20}}
\put(22,10){\line(6,1){120}}
\put(22,10){\circle*{5}}

\put(60,33){$\gneg R$}
\put(86,33){$\gneg S$}
\put(78,33){$\mld$}
\put(82,10){\line(-3,1){59}}
\put(82,10){\circle*{5}}
\put(82,10){\line(6,1){120}}

\put(124,33){$P$}
\put(152,33){$R$}
\put(139,33){$\mld$}
\put(142,10){\line(-3,1){60}}
\put(142,10){\line(0,1){20}}
\put(142,10){\circle*{5}}

\put(184,33){$Q$}
\put(210,33){$S$}
\put(202,10){\line(-6,1){118}}
\put(202,10){\line(0,1){20}}
\put(202,10){\circle*{5}}
\put(198,33){$\mld$}
\end{picture}
\end{center}

\begin{center}\begin{picture}(228,34)
\put(0,46){\line(1,0){223}}
\put(225,45){\scriptsize $\mlc$}
\put(0,33){$\gneg P$}
\put(28,33){$\gneg Q$}
\put(19,33){$\mld$}

\put(60,33){$\gneg R$}
\put(86,33){$\gneg S$}
\put(78,33){$\mld$}
\put(82,10){\line(-3,1){59}}
\put(82,10){\circle*{5}}
\put(82,10){\line(5,1){50}}
\put(131,20){\line(4,1){41}}

\put(121,33){$(P$}
\put(151,33){$R)$}
\put(139,33){$\mld$}
\put(142,10){\line(-3,1){60}}
\put(142,10){\line(3,2){30}}
\put(142,10){\circle*{5}}

\put(169,33){$\mlc$}

\put(182,33){$(Q$}
\put(210,33){$S)$}

\put(198,33){$\mld$}
\end{picture}
\end{center}

\begin{center}\begin{picture}(228,34)
\put(-3,46){\line(1,0){226}}
\put(225,45){\scriptsize $\mlc$}
\put(-3,33){$(\gneg P$}
\put(28,33){$\gneg Q$)}
\put(19,33){$\mld$}

\put(49,33){$\mlc$}

\put(57,33){$(\gneg R$}
\put(86,33){$\gneg S)$}
\put(78,33){$\mld$}
\put(113,10){\line(-3,1){60}}
\put(113,10){\circle*{5}}
\put(113,10){\line(3,1){59}}

\put(121,33){$(P$}
\put(151,33){$R)$}
\put(139,33){$\mld$}

\put(169,33){$\mlc$}

\put(182,33){$(Q$}
\put(210,33){$S)$}

\put(198,33){$\mld$}
\end{picture}
\end{center}

\begin{center}\begin{picture}(228,11)
\put(-6,23){\line(1,0){232}}
\put(228,22){\scriptsize $\mld$}
\put(-6,10){$\bigl((\gneg P$}
\put(28,10){$\gneg Q$)}
\put(19,10){$\mld$}

\put(49,10){$\mlc$}

\put(57,10){$(\gneg R$}
\put(86,10){$\gneg S)\bigr)$}
\put(78,10){$\mld$}
\put(110,10){$\mld$}

\put(118,10){$\bigl((P$}
\put(151,10){$R)$}
\put(139,10){$\mld$}

\put(169,10){$\mlc$}

\put(182,10){$(Q$}
\put(210,10){$S)\bigr)$}

\put(198,10){$\mld$}
\end{picture}
\end{center}

\section{Classical and affine logics}

In {\bf sequent calculus} (where a {\bf sequent} means a nonempty sequence of formulas), {\bf classical logic} can be axiomatized by the following six rules, where $F,G$ stand for any formulas and $\Gamma,\Delta$ stand for any --- possibly empty --- sequences of formulas:\vspace{10pt}

\begin{center}\begin{picture}(327,68)
\put(2,58){\bf Axiom:} 
\put(6,23){\line(1,0){25}}
\put(33,21){\scriptsize A}
\put(6,10){$\gneg F,F$}

\put(80,58){\bf Exchange:} 
\put(82,42){\line(1,0){42}}
\put(82,30){$\Gamma,F,G,\Delta$}
\put(82,23){\line(1,0){42}}
\put(126,21){\scriptsize E}
\put(82,10){$\Gamma,G,F,\Delta$}

\put(166,58){\bf Weakening:}
\put(181,42){\line(1,0){19}}
\put(181,30){$\Gamma,\Delta$} 
\put(176,23){\line(1,0){29}}
\put(207,21){\scriptsize W}
\put(176,10){$\Gamma,F,\Delta$}

\put(262,58){\bf Contraction:} 
\put(269,42){\line(1,0){42}}
\put(269,30){$\Gamma,F,F,\Delta$}
\put(269,23){\line(1,0){42}}
\put(313,21){\scriptsize C}
\put(276,10){$\Gamma,F,\Delta$}

\end{picture}\end{center}

\begin{center}\begin{picture}(274,78)
\put(20,58){\bf $\mld$-introduction:} 
\put(35,42){\line(1,0){43}}
\put(35,30){$\Gamma,F,G,\Delta$}
\put(32,23){\line(1,0){48}}
\put(82,21){\scriptsize $\mld$}
\put(32,10){$\Gamma,F\mld G,\Delta$}

\put(177,58){\bf $\mlc$-introduction:}
\put(182,42){\line(1,0){19}} 
\put(222,42){\line(1,0){21}}
\put(182,30){$\Gamma,F$}
\put(222,30){$G,\Delta$}
\put(182,23){\line(1,0){59}}
\put(243,21){\scriptsize $\mlc$}
\put(187,10){$\Gamma,F\mlc G,\Delta$}
\end{picture}\end{center}

{\bf Affine logic} is obtained from classical logic by deleting contraction. As noted earlier, the term ``affine logic" in this paper refers to what is called the {\em multiplicative fragment} of this otherwise more expressive logic.
A {\bf sequent calculus system}, in general, is any subset of the above six rules.  
The definition of {\bf provability} of a sequent $\Gamma$ in a sequent calculus system $S$ is standard: this means existence of a tree of sequents --- called a {\bf proof tree} for $\Gamma$ --- with $\Gamma$ at its root, in which every node of the tree follows from its children (where the set of children may be empty in the case of axiom) by one of the rules of $S$. A formula $F$ is considered provable in a sequent calculus system iff $F$, viewed as a one-element sequent, is provable.

At the end of Section \ref{sccs} we saw that cirquent calculus needs neither weakening nor contraction (nor duplication) to prove Blass's principle. Replacing all atoms by $P$ in our proof tree 
for Blass's principle also yields an (AME$\mld\mlc$)-proof of  
\begin{equation}\label{bp} 
\bigl((\gneg P\mld \gneg P)\mlc(\gneg P\mld\gneg P)\bigr)\mld \bigl((P\mld P)\mlc(P\mld P)\bigr).
\end{equation}
The following Fact \ref{f1} establishes that, in contrast, sequent calculus needs {\em both} weakening and contraction to prove (\ref{bp}), let alone the more general Blass's principle.\footnote{That affine logic does not prove (\ref{bp}) was shown by Blass in \cite{Bla92}.}   

\begin{fact}\label{f1}
Any proof of (\ref{bp}) in sequent calculus would have to use both weakening and contraction. 
\end{fact}

\begin{proof} First, let us attempt to construct, in a bottom-up fashion, a proof of (\ref{bp}) in affine logic to see that such a proof does not exist. 
The only rule that can yield (\ref{bp}) is $\mld$-introduction, so 
the premise should be the sequent 

\[(\gneg P\mld \gneg P)\mlc(\gneg P\mld\gneg P),\ (P\mld P)\mlc(P\mld P).\]

Weakening is not applicable to the above sequent, for both of its formulas are non-valid in the classical sense and hence, in view of the known fact that all of the sequent calculus rules preserve classical validity, those formulas, in isolation, are not provable. $\mld$-introduction is not applicable, either, for there is no disjunction on the surface of the sequent. And exchange, of course,  would not take us closer to our goal of finding a proof. This leaves us with $\mlc$-introduction. The sequent is symmetric, so we may assume that the introduced conjunction is, say, the first one. The non-active formula $(P\mld P)\mlc(P\mld P)$ of the conclusion can then only be inherited by one of the premises, meaning that the other premise 
will be just $\gneg P\mld \gneg P$. Now we are stuck with that premise, as it is a non-tautological formula which cannot be proven.

Next, for a contradiction, assume that there is a weakening-free (but not necessarily contraction-free) sequent calculus proof of (\ref{bp}). Consider any branch 
$\Gamma_1,\ldots,\Gamma_n$ of the proof tree, where $\Gamma_1$ should be $\gneg F,F$ for some formula $F$, and $\Gamma_n$ be the sequent consisting just of (\ref{bp}). Notice that once a given sequent $\Gamma_i$ of the above sequence contains a formula $G$, $G$ will be inherited by 
each of the subsequent sequents $\Gamma_{i+1},\ldots,\Gamma_n$ --- either as a formula of the sequent, or as a subformula of such. So, both $F$ and $\gneg F$ should be subformulas of (\ref{bp}). This leaves us only with the possibility 
$\{F,\gneg F\}=\{P,\gneg P\}$, because (\ref{bp}) does not contain any other subformula $F$ together with $\gneg F$.  Let $i$ be the greatest number 
among \ $1,\ldots,n-1$ \ such that $\Gamma_{i+1}$ is neither $\gneg P,P$ nor $P,\gneg P$. $\Gamma_{i+1}$ cannot be derived from $\Gamma_i$ by exchange because then $\Gamma_{i+1}$ would again be $\gneg P,P$ or $P,\gneg P$. Nor can it be derived by contraction which is simply not applicable to $\Gamma_i$. Nor can $\Gamma_{i+1}$ be derived by $\mld$-introduction, because then $\Gamma_{i+1}$ would be $\gneg P\mld P$ or $P\mld \gneg P$, which is not a subformula of (\ref{bp}). Finally, $\Gamma_{i+1}$ cannot be derived from $\Gamma_i$ 
(and an arbitrary other premise) by $\mlc$-introduction either. This is so because an application of this rule would introduce a conjunction where $\gneg P$ or $P$ is a conjunct; but, again, (\ref{bp}) does not have such a subformula.  
\end{proof}

As we just saw, cirquent calculus indeed offers a substantially more flexible machinery for constructing (substructural) deductive systems than sequent calculus does. Sequent calculus can be seen as a simple special case of cirquent calculus that we call ``primitive". Specifically, we say that  
a cirquent is  {\bf primitive} iff all of its ogroups are (pairwise) disjoint. 
The groups of such a cirquent can be thought of as --- and identified with --- sequents: in this section we will not terminologically distinguish between a group $\Gamma$ of a primitive cirquent and the sequent $\Gamma$ consisting exactly of the oformulas that $\Gamma$ contains, arranged in the same order as they appear in the pool of the cirquent.   

For any given cirquent calculus system $S$, we let $S^*$ denote the version of $S$ where the definition of a proof or a derivation has the additional condition that every cirquent in the proof or derivation should be primitive. So, $S^*$ can be called the ``primitive version" of $S$. Of course, $S$ proves or derives everything that $S^*$ does. 

 Strictly speaking, a sequent or cirquent calculus system is the particular collection of its rules, so that even if two systems prove exactly the same formulas or sequents or cirquents, they should count as different systems. Yet, often  we identify a sequent or cirquent calculus system with the set of formulas (or sequents, or cirquents) provable in it, as done in the following Theorem \ref{ll}. 
The equalities in the left column of that theorem, as can be seen from our subsequent proof of it, extend to all other natural pairs of systems obtained by allowing/disallowing various rules, such as affine logic without weakening (i.e. linear logic in the proper sense) vs. the primitive version of {\bf CL5} without weakening, or classical logic without weakening vs. the primitive version of {\bf CCC} without weakening.
It is such equalities that allow us to say that sequent calculus is nothing but the primitive version 
of cirquent calculus. That primitiveness makes cirquent calculus degenerate to sequent calculus is no surprise. The former owes its special power to the ability to express resource sharing, and it is exactly resource sharing that primitive cirquents forbid. 

\begin{theorem}\label{ll} With the following sequent calculus and cirquent calculus systems identified with the sets of formulas that they prove, we have:
\[\begin{array}{lclclcl}
\mbox{1. Affine logic} & = & {\mathbf CL5}^* & \subseteq &  {\mathbf CL5} & \not= & \mbox{Affine logic.}\\
\mbox{2. Classical logic} & = & {\mathbf CCC}^* & \subseteq & {\mathbf CCC} & = & \mbox{Classical logic.}
\end{array}\]
\end{theorem}

\begin{proof} As noted earlier, the inclusions of the type $S^*\subseteq S$ are trivial. The inequality  
{\bf CL5} $\not=$ {\em Affine logic} \ immediately follow from Fact \ref{f1} together with the earlier-established provability of (\ref{bp}) in (AME$\mld\mlc$). The equality  {\bf CCC} $=$ {\em Classical logic}  follows from Theorem \ref{clascomp}, which will be proven in the next section. The latter implies 
that a formula is provable in {\bf CCC} iff it is a tautology in the classical sense, and it just remains to remember that the same is known to be true for 
{\em Classical logic}.  

Our task now is to verify the equalities \ {\em Affine logic} = {\bf CL5}$^*$ \ and \ {\em Classical logic} = {\bf CCC}$^*$. 

The inclusions  {\em Affine logic} $\subseteq$ {\bf CL5}$^*$  \  and \  {\em Classical logic} $\subseteq$ {\bf CCC}$^*$ can be proven by showing that whenever either sequent calculus system proves a sequent $\Gamma$, the corresponding 
primitive cirquent calculus system proves the cirquent whose only group --- as well as pool ---  is $\Gamma$. 
This can be easily done by induction on the heights of proof trees. The steps of such induction are rather straightforward, for every application of a sequent calculus rule --- except weakening and $\mlc$-introduction ---  directly translates into an application of the same-name rule of cirquent calculus as shown below: \vspace{10pt}

\begin{center}\begin{picture}(309,44)
\put(150,41){$\Longleftrightarrow$}
\put(99,43){\line(1,0){26}}
\put(99,30){$\gneg F,F$} 
\put(187,43){\line(1,0){30}}
\put(219,41){\scriptsize A}
\put(127,41){\scriptsize A}
\put(187,30){$\gneg F\hspace{8pt} F$}
\put(203,10){\circle*{5}}
\put(203,10){\line(-1,2){8}}
\put(203,10){\line(1,2){8}} 
\end{picture}\end{center}

\begin{center}\begin{picture}(309,92)
\put(150,41){$\Longleftrightarrow$}
\put(187,82){\line(1,0){115}}
\put(187,69){$E_1 \ldots E_n \ F\hspace{10pt} G \ H_1 \ldots H_m$}
\put(243,49){\circle*{5}}
\put(243,49){\line(-3,1){48}}  
\put(243,49){\line(-6,5){20}}
\put(243,49){\line(-1,2){8}}
\put(243,49){\line(1,2){8}}
\put(243,49){\line(6,5){20}}
\put(243,49){\line(3,1){48}} 

\put(-1,49){$E_1,\ldots,E_n,F,G,H_1,\ldots,H_m$} 
\put(-1,43){\line(1,0){125}}
\put(-1,62){\line(1,0){125}}
\put(-1,30){$E_1,\ldots,E_n,G,F,H_1,\ldots,H_m$} 
\put(187,43){\line(1,0){115}}
\put(306,41){\scriptsize E}
\put(127,41){\scriptsize E}
\put(187,30){$E_1 \ldots E_n \ G\hspace{10pt} F \ H_1 \ldots H_m$}
\put(243,10){\circle*{5}}
\put(243,10){\line(-3,1){48}}  
\put(243,10){\line(-6,5){20}}
\put(243,10){\line(-1,2){8}}
\put(243,10){\line(1,2){8}}
\put(243,10){\line(6,5){20}}
\put(243,10){\line(3,1){48}} 
\end{picture}\end{center}

\begin{center}\begin{picture}(309,92)

\put(150,41){$\Longleftrightarrow$}
\put(190,82){\line(1,0){115}}
\put(190,69){$E_1 \ldots E_n \ F\hspace{10pt} F \ H_1 \ldots H_m$}
\put(246,49){\circle*{5}}
\put(246,49){\line(-3,1){48}}  
\put(246,49){\line(-6,5){20}}
\put(246,49){\line(-1,2){8}}
\put(246,49){\line(1,2){8}}
\put(246,49){\line(6,5){20}}
\put(246,49){\line(3,1){48}} 

\put(2,49){$E_1,\ldots,E_n,F,F,H_1,\ldots,H_m$} 
\put(2,43){\line(1,0){124}}
\put(2,62){\line(1,0){124}}
\put(7,30){$E_1,\ldots,E_n,F,H_1,\ldots,H_m$} 
\put(190,43){\line(1,0){115}}
\put(308,41){\scriptsize C}
\put(128,41){\scriptsize C}
\put(190,30){$E_1 \ldots E_n\hspace{13pt} F\hspace{13pt}  H_1 \ldots H_m$}
\put(246,10){\circle*{5}}
\put(246,10){\line(-3,1){48}}  
\put(246,10){\line(-6,5){20}}
\put(246,10){\line(0,1){17}}
\put(246,10){\line(6,5){20}}
\put(246,10){\line(3,1){48}} 
\end{picture}\end{center}

\begin{center}\begin{picture}(309,92)
\put(150,41){$\Longleftrightarrow$}
\put(190,82){\line(1,0){115}}
\put(190,69){$E_1 \ldots E_n \ F\hspace{10pt} G \ H_1 \ldots H_m$}
\put(246,49){\circle*{5}}
\put(246,49){\line(-3,1){48}}  
\put(246,49){\line(-6,5){20}}
\put(246,49){\line(-1,2){8}}
\put(246,49){\line(1,2){8}}
\put(246,49){\line(6,5){20}}
\put(246,49){\line(3,1){48}} 

\put(-3,49){$E_1,\ldots,E_n,F,G,H_1,\ldots,H_m$} 
\put(-6,43){\line(1,0){132}}
\put(-2,62){\line(1,0){125}}
\put(-6,30){$E_1,\ldots,E_n,F\mld G,H_1,\ldots,H_m$} 
\put(190,43){\line(1,0){115}}
\put(308,41){\scriptsize $\mld$}
\put(128,41){\scriptsize $\mld$}
\put(190,30){$E_1 \ldots E_n \ F\mld G \ H_1 \ldots H_m$}
\put(246,10){\circle*{5}}
\put(246,10){\line(-3,1){48}}  
\put(246,10){\line(-6,5){20}}
\put(246,10){\line(0,1){17}}
\put(246,10){\line(6,5){20}}
\put(246,10){\line(3,1){48}} 
\end{picture}\end{center}

As for weakening and $\mlc$-introduction, their sequent-calculus to cirquent-calculus translations take two steps:

\begin{center}\begin{picture}(309,131)

\put(16,69){$E_1,\ldots,E_n,H_1,\ldots,H_m$} 
\put(11,63){\line(1,0){112}}
\put(16,82){\line(1,0){102}}
\put(10,50){$E_1,\ldots,E_n,F,H_1,\ldots,H_m$} 
\put(125,61){\scriptsize W}
\put(150,61){$\Longrightarrow$}

\put(190,121){\line(1,0){115}}
\put(190,108){$E_1 \ldots E_n \hspace{30pt} \ H_1 \ldots H_m$}
\put(246,88){\circle*{5}}
\put(246,88){\line(-3,1){48}}  
\put(246,88){\line(-6,5){20}}
\put(246,88){\line(6,5){20}}
\put(246,88){\line(3,1){48}}

\put(190,82){\line(1,0){115}}
\put(309,80){\scriptsize W}
\put(190,69){$E_1 \ldots E_n\hspace{13pt} F\hspace{13pt}  H_1 \ldots H_m$}
\put(246,49){\circle*{5}}
\put(246,49){\line(-3,1){48}}  
\put(246,49){\line(-6,5){20}}
\put(246,49){\line(6,5){20}}
\put(246,49){\line(3,1){48}}

\put(190,43){\line(1,0){115}}
\put(309,41){\scriptsize W}
\put(190,30){$E_1 \ldots E_n\hspace{13pt} F\hspace{13pt}  H_1 \ldots H_m$}
\put(246,10){\circle*{5}}
\put(246,10){\line(-3,1){48}}  
\put(246,10){\line(-6,5){20}}
\put(246,10){\line(0,1){17}}
\put(246,10){\line(6,5){20}}
\put(246,10){\line(3,1){48}} 
\end{picture}\end{center}

\begin{center}\begin{picture}(309,125)
\put(190,117){\line(1,0){51}}
\put(253,117){\line(1,0){52}}
\put(190,105){$E_1 \ldots E_n \ F\hspace{11pt} G \ H_1 \ldots H_m$}
\put(237,85){\circle*{5}}
\put(255,85){\circle*{5}}
\put(237,85){\line(-5,2){41}}  
\put(237,85){\line(-4,5){14}}
\put(237,85){\line(0,1){17}}
\put(255,85){\line(5,2){41}}  
\put(255,85){\line(4,5){14}}
\put(255,85){\line(0,1){17}}

\put(150,61){$\Longrightarrow$}
\put(190,79){\line(1,0){115}}
\put(307,77){\scriptsize M}
\put(190,69){$E_1 \ldots E_n \ F\hspace{11pt} G \ H_1 \ldots H_m$}
\put(237,49){\circle*{5}}
\put(255,49){\circle*{5}}
\put(237,49){\line(-5,2){41}}  
\put(237,49){\line(-4,5){14}}
\put(237,49){\line(0,1){17}}
\put(255,49){\line(5,2){41}}  
\put(255,49){\line(4,5){14}}
\put(255,49){\line(0,1){17}}

\put(-20,69){$E_1,\ldots,E_n,F\hspace{21pt}G,H_1,\ldots,H_m$} 
\put(-20,63){\line(1,0){143}}
\put(-20,82){\line(1,0){60}}
\put(60,82){\line(1,0){61}}
\put(126,61){\scriptsize $\mlc$}
\put(-15,50){$E_1,\ldots,E_n,F\mlc G,H_1,\ldots,H_m$} 
\put(190,43){\line(1,0){115}}
\put(309,41){\scriptsize $\mlc$}
\put(190,30){$E_1 \ldots E_n \ F\mlc G \ H_1 \ldots H_m$}
\put(246,10){\circle*{5}}
\put(246,10){\line(-3,1){48}}  
\put(246,10){\line(-6,5){20}}
\put(246,10){\line(0,1){17}}
\put(246,10){\line(6,5){20}}
\put(246,10){\line(3,1){48}} 
\end{picture}\end{center}

The inclusions  {\bf CL5}$^*$ $\subseteq$ {\em Affine logic}  and 
{\bf CCC}$^*$ $\subseteq$ {\em Classical logic} can be verified in a rather similar way. Specifically, this can be done by showing that, whenever the primitive version of either cirquent calculus system 
proves a given cirquent $C$, the corresponding sequent calculus system proves each of the groups  
of $C$ understood as sequents.  
Induction on the heights of proof trees is again the way to proceed. 
The basis of induction is straightforward, taking into account that the translation shown earlier for identity axiom works in either direction, and that the case of the empty cirquent axiom is vacuous as there are no groups in its ``conclusion''.  Duplication does not need to be considered, for either the premise or the conclusion of it has to be non-primitive, which automatically bans this rule in primitive cirquent calculus systems. The inductive steps dealing with mix or pool weakening are trivial, because these rules do not create new groups or affect the contents of the existing groups. The same is true for ogroup exchange, as well as oformula exchange if it is {\em external}, i.e. swaps oformulas that are in different groups, as opposed to {\em internal} exchange that swaps oformulas that are in the same group.  
What internal oformula exchange, ogroup weakening, contraction, $\mld$-introduction and $\mlc$-introduction do in primitive cirquents is that they modify 
one or two of the groups of the premise without affecting any other groups (if there are such). This local behavior allows us to pretend for our present purposes that simply there are no other groups in the cirquent under question. So, in inductive steps dealing with internal oformula exchange, contraction and $\mld$-introduction, we can rely on the fact that the above-illustrated translations between the sequent- and cirquent-calculus versions of these rules work in either direction. As for ogroup weakening and $\mlc$-introduction, their cirquent-calculus to sequent-calculus translations work as follows:\vspace{10pt}

\begin{center}\begin{picture}(309,92)

\put(16,49){$E_1,\ldots,E_n,H_1,\ldots,H_m$} 
\put(11,43){\line(1,0){112}}
\put(16,62){\line(1,0){102}}
\put(10,30){$E_1,\ldots,E_n,F,H_1,\ldots,H_m$} 
\put(125,41){\scriptsize W}
\put(150,41){$\Longleftarrow$}

\put(190,82){\line(1,0){115}}
\put(190,69){$E_1 \ldots E_n\hspace{13pt} F\hspace{13pt}  H_1 \ldots H_m$}
\put(246,49){\circle*{5}}
\put(246,49){\line(-3,1){48}}  
\put(246,49){\line(-6,5){20}}
\put(246,49){\line(6,5){20}}
\put(246,49){\line(3,1){48}}

\put(190,43){\line(1,0){115}}
\put(309,41){\scriptsize W}
\put(190,30){$E_1 \ldots E_n\hspace{13pt} F\hspace{13pt}  H_1 \ldots H_m$}
\put(246,10){\circle*{5}}
\put(246,10){\line(-3,1){48}}  
\put(246,10){\line(-6,5){20}}
\put(246,10){\line(0,1){17}}
\put(246,10){\line(6,5){20}}
\put(246,10){\line(3,1){48}} 
\end{picture}\end{center}

\begin{center}\begin{picture}(309,92)

\put(150,41){$\Longleftarrow$}
\put(190,82){\line(1,0){115}}
\put(190,69){$E_1 \ldots E_n \ F\hspace{11pt} G \ H_1 \ldots H_m$}
\put(237,49){\circle*{5}}
\put(255,49){\circle*{5}}
\put(237,49){\line(-5,2){41}}  
\put(237,49){\line(-4,5){14}}
\put(237,49){\line(0,1){17}}
\put(255,49){\line(5,2){41}}  
\put(255,49){\line(4,5){14}}
\put(255,49){\line(0,1){17}}

\put(-18,62){\line(1,0){58}}
\put(62,62){\line(1,0){61}}
\put(-18,49){$E_1,\ldots,E_n,F\hspace{21pt}G,H_1,\ldots,H_m$} 
\put(-18,43){\line(1,0){142}}
\put(-13,30){$E_1,\ldots,E_n,F\mlc G,H_1,\ldots,H_m$} 
\put(190,43){\line(1,0){115}}
\put(309,41){\scriptsize $\mlc$}
\put(126,41){\scriptsize $\mlc$}
\put(190,30){$E_1 \ldots E_n \ F\mlc G \ H_1 \ldots H_m$}
\put(246,10){\circle*{5}}
\put(246,10){\line(-3,1){48}}  
\put(246,10){\line(-6,5){20}}
\put(246,10){\line(0,1){17}}
\put(246,10){\line(6,5){20}}
\put(246,10){\line(3,1){48}} 
\end{picture}\end{center}

\end{proof}

\section{Tautologies}\label{stau}
By a {\bf classical model}, or simply {\bf model}, we mean a function $M$ that assigns a truth value --- {\em true} ($1$) or {\em false} ($0$) --- to 
each atom, and extends to compound formulas in the standard classical way. The traditional concepts of truth and tautologicity naturally extend from formulas to groups and  
cirquents.  Let $M$ be a model, and $C$ a cirquent. We say that a group $\Gamma$ of $C$ is {\bf true} in $M$ iff at least one of its oformulas is so. And 
$C$ is {\bf true} in $M$ if every group of $C$ is so. 
``{\bf False}", as always, means ``not true". 
Finally, $C$ or a group $\Gamma$ of it is a {\bf tautology} iff it is true in every model. Identifying each formula $F$ with the singleton 
 cirquent $(\langle\{1\}\rangle,\langle F\rangle)$, our concepts of truth and tautologicity of cirquents preserve the standard meaning of these concepts for formulas. Let us mark the evident fact that 
a cirquent is tautological if and only if all of its groups are so. Note also that a cirquent containing the empty group is always false, while a cirquent with no groups, such as 
the empty cirquent $(\seq{},\seq{})$, is always true.

\begin{lemma}\label{tau}
All of the rules of Section \ref{srules} preserve truth in the top-down direction --- that is, whenever the premise(s) of an application of any given rule is (are) true in a given model, so is the conclusion. Taking no premises, (the conclusions of) axioms are thus tautologies.
\end{lemma}

\begin{proof} A routine examination of those rules and our definition of truth for cirquents. \end{proof} 

\begin{lemma}\label{tau2}
The rules of mix, exchange, duplication, contraction, conservative $\mld$-introduction and conservative $\mlc$-introduction preserve truth in the bottom-up direction as well --- that is, whenever the conclusion of an application of such a rule is true in a given model, so is (are) the premise(s). 
\end{lemma}
 
\begin{proof} The above statement for mix, exchange, duplication and contraction is rather obvious. Let us only examine it for 
the conservative versions of $\mld$- and $\mlc$-introduction. 

{\em Conservative $\mld$-introduction}: Assume the disjunction that the rule introduced is $F\mld G$. Notice that the only difference between the conclusion and the premise is that wherever the conclusion has an ogroup $\Gamma$ containing the oformula $F\mld G$, the premise has the ogroup $\Gamma'=(\Gamma-\{F\mld G\})\cup\{F,G\}$ instead.  Since truth-semantically a group is nothing but the disjunction of its oformulas, the truth values of $\Gamma$ and $\Gamma'$ (in whatever model $M$) are the same. Hence so are those of the conclusion and the premise.  

{\em Conservative $\mlc$-introduction}:  
Assume the conjunction that the rule introduced is $F\mlc G$. The only difference between the conclusion and the premise is that wherever the conclusion has an ogroup $\Gamma$ containing $F\mlc G$, the premise has the two ogroups 
$\Gamma_F=(\Gamma-\{F\mlc G\})\cup\{F\}$ and $\Gamma_G=(\Gamma-\{F\mlc G\})\cup\{G\}$ instead. Obviously this implies that if $\Gamma$ is true in a given model, then so are both $\Gamma_F$ and $\Gamma_G$. The above, in turn, implies that if the conclusion is true and hence all of its groups are true, then so are all of the groups of the premise, and hence the premise itself. 
\end{proof}

We say that an oformula $F$ of (the pool of) a cirquent $C$ is {\bf homeless} iff no group of $C$ contains $F$. 

A {\bf literal} means $P$ ({\bf positive literal of type $P$}) or $\gneg P$ ({\bf negative literal of type $P$}) for some atom $P$. The term {\bf oliteral} has the expected meaning: this is a particular occurrence of a literal in a formula or in (the pool or an oformula of) a cirquent. A {\bf literal cirquent} is a cirquent whose pool contains only literals. An {\bf essentially literal} cirquent is one every oformula of whose pool either is an oliteral or is homeless.

\begin{theorem}\label{clascomp}
A cirquent is provable in {\bf CCC} iff it is a tautology.
\end{theorem}

\begin{proof} The soundness part of this theorem is an immediate corollary of Lemma \ref{tau}. For the  completeness part, consider any tautological cirquent $A$. 
 In the bottom-up sense, keep applying to it\footnote{As this can be understood, here and later in similar contexts, ``it'' means $A$ only in the beginning, then ``it'' becomes the premise of $A$, then the premise of the premise of $A$, and so on.}  conservative $\mld$-introduction and conservative $\mlc$-introduction --- in whatever order you like --- until you hit an
essentially literal cirquent $B$, such as the one shown in the following example:

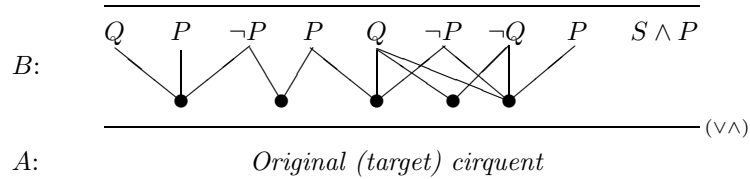
\begin{center}\begin{picture}(270,82)
\put(35,72){\line(1,0){225}}
\put(0,46){$B$:}
\put(35,59){$Q$}
\put(60,59){$P$}
\put(82,59){$\gneg P$}
\put(110,59){$P$}
\put(134,59){$Q$}
\put(156,59){$\gneg P$}
\put(180,59){$\gneg Q$}
\put(210,59){$P$}
\put(234,59){$S\mlc P$}

\put(64,36){\line(-5,4){25}}
\put(64,36){\line(5,4){26}}
\put(64,36){\line(0,1){19}}
\put(64,36){\circle*{5}}

\put(102,36){\line(-3,5){12}}
\put(102,36){\line(3,5){12}}
\put(102,36){\circle*{5}}

\put(138,36){\line(-5,4){25}}
\put(138,36){\line(5,4){26}}
\put(138,36){\line(0,1){19}}
\put(138,36){\circle*{5}}

\put(167,36){\line(-3,2){29}}
\put(167,36){\line(1,1){21}}
\put(167,36){\circle*{5}}

\put(188,36){\line(-5,2){50}}
\put(188,36){\line(0,1){21}}
\put(188,36){\line(-6,5){25}}
\put(188,36){\line(6,5){25}}
\put(188,36){\circle*{5}}

\put(35,26){\line(1,0){225}}
\put(262,24){\scriptsize ($\mld\mlc$)}
\put(90,10){\em Original (target) cirquent}
\put(0,10){$A$:}
\end{picture}\end{center}

The above procedure will indeed always hit an essentially literal cirquent because conservative $\mld$-introduction is obviously always applicable when a given cirquent  (conclusion)  has a non-homeless oformula $F\mld G$, and so is conservative $\mlc$-introduction whenever such a cirquent has a non-homeless oformula $F\mlc G$. $A$ thus follows from $B$ in ($\mld\mlc$).
In view of Lemma \ref{tau2}, $B$ is 
a tautology, and since all of its non-homeless oformulas are literals, the tautologicity of $B$ obviously means that every group of it contains at least one pair of $P,\gneg P$ of same-type positive and negative oliterals. Fix one such pair for each group, and then apply (in the bottom-up sense) to $B$ a series of weakenings
to first delete, in each group,  all arcs but the two arcs pointing to the two chosen oliterals, and next delete all homeless oformulas if any such oformulas are present. This is illustrated below:
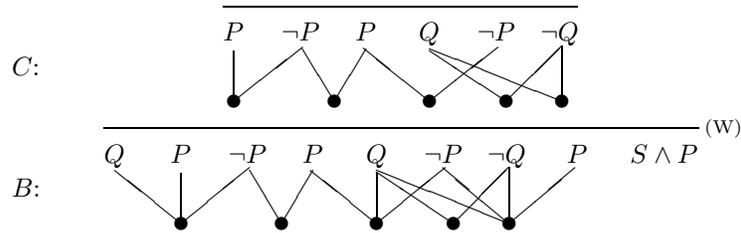
\vspace{10pt}

\begin{center}\begin{picture}(264,102)
\put(80,92){\line(1,0){134}}
\put(0,66){$C$:}
\put(80,79){$P$}
\put(102,79){$\gneg P$}
\put(130,79){$P$}
\put(154,79){$Q$}
\put(176,79){$\gneg P$}
\put(200,79){$\gneg Q$}

\put(84,56){\line(5,4){26}}
\put(84,56){\line(0,1){19}}
\put(84,56){\circle*{5}}

\put(122,56){\line(-3,5){12}}
\put(122,56){\line(3,5){12}}
\put(122,56){\circle*{5}}

\put(158,56){\line(-5,4){25}}
\put(158,56){\line(5,4){26}}
\put(158,56){\circle*{5}}

\put(187,56){\line(-3,2){29}}
\put(187,56){\line(1,1){21}}
\put(187,56){\circle*{5}}

\put(208,56){\line(-5,2){50}}
\put(208,56){\line(0,1){21}}
\put(208,56){\circle*{5}}

\put(35,46){\line(1,0){225}}
\put(262,44){\scriptsize (W)}
\put(0,20){$B$:}

\put(35,33){$Q$}
\put(60,33){$P$}
\put(82,33){$\gneg P$}
\put(110,33){$P$}
\put(134,33){$Q$}
\put(156,33){$\gneg P$}
\put(180,33){$\gneg Q$}
\put(210,33){$P$}
\put(234,33){$S\mlc P$}

\put(64,10){\line(-5,4){25}}
\put(64,10){\line(5,4){26}}
\put(64,10){\line(0,1){19}}
\put(64,10){\circle*{5}}

\put(102,10){\line(-3,5){12}}
\put(102,10){\line(3,5){12}}
\put(102,10){\circle*{5}}

\put(138,10){\line(-5,4){25}}
\put(138,10){\line(5,4){26}}
\put(138,10){\line(0,1){19}}
\put(138,10){\circle*{5}}

\put(167,10){\line(-3,2){29}}
\put(167,10){\line(1,1){21}}
\put(167,10){\circle*{5}}

\put(188,10){\line(-5,2){50}}
\put(188,10){\line(0,1){21}}
\put(188,10){\line(-6,5){25}}
\put(188,10){\line(6,5){25}}
\put(188,10){\circle*{5}}
\end{picture}\end{center}

Every group of the resulting cirquent $C$ will thus have exactly two oformulas: some atom and its negation. 
$B$ follows from $C$ in (W), so that $A$ follows from 
$C$ in (W$\mld\mlc$). Now apply (again in the bottom-up fashion) a series of contractions to $C$ to separate all shared oliterals, as illustrated in the example below, with the resulting cirquent called $D$:
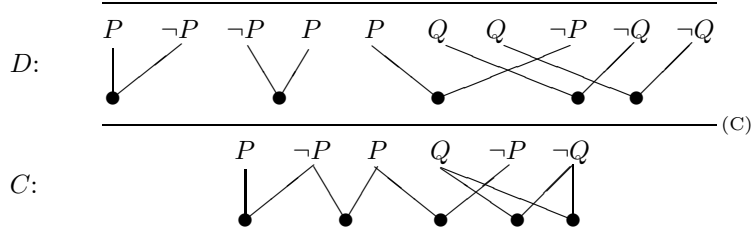
\vspace{10pt}

\begin{center}\begin{picture}(273,102)
\put(35,92){\line(1,0){232}}
\put(0,66){$D$:}
\put(35,79){$P$}
\put(57,79){$\gneg P$}
\put(82,79){$\gneg P$}
\put(110,79){$P$}
\put(134,79){$P$}
\put(158,79){$Q$}
\put(180,79){$Q$}
\put(204,79){$\gneg P$}
\put(228,79){$\gneg Q$}
\put(252,79){$\gneg Q$}

\put(39,56){\line(5,4){26}}
\put(39,56){\line(0,1){19}}
\put(39,56){\circle*{5}}

\put(102,56){\line(-3,5){12}}
\put(102,56){\line(3,5){11}}
\put(102,56){\circle*{5}}

\put(162,56){\line(-5,4){25}}
\put(162,56){\line(5,2){50}}
\put(162,56){\circle*{5}}

\put(215,56){\line(1,1){21}}
\put(215,56){\line(-5,2){50}}
\put(215,56){\circle*{5}}

\put(237,56){\line(1,1){21}}
\put(237,56){\line(-5,2){50}}
\put(237,56){\circle*{5}}

\put(35,46){\line(1,0){232}}
\put(269,44){\scriptsize (C)}
\put(0,20){$C$:}

\put(85,33){$P$}
\put(107,33){$\gneg P$}
\put(135,33){$P$}
\put(159,33){$Q$}
\put(181,33){$\gneg P$}
\put(205,33){$\gneg Q$}

\put(89,10){\line(5,4){26}}
\put(89,10){\line(0,1){19}}
\put(89,10){\circle*{5}}

\put(127,10){\line(-3,5){12}}
\put(127,10){\line(3,5){12}}
\put(127,10){\circle*{5}}

\put(163,10){\line(-5,4){25}}
\put(163,10){\line(5,4){26}}
\put(163,10){\circle*{5}}

\put(192,10){\line(-3,2){29}}
\put(192,10){\line(1,1){21}}
\put(192,10){\circle*{5}}

\put(213,10){\line(-5,2){50}}
\put(213,10){\line(0,1){21}}
\put(213,10){\circle*{5}}

\end{picture}\end{center}
 
So, our original cirquent $A$ is derivable from $D$ in (WC$\mld\mlc$). Every ogroup of $D$ is disjoint from every other ogroup and, as in $C$, each such ogroup contains exactly two oformulas: $P$ and $\gneg P$ for some atom $P$. Therefore $D$ is provable in 
(AME) as illustrated below:

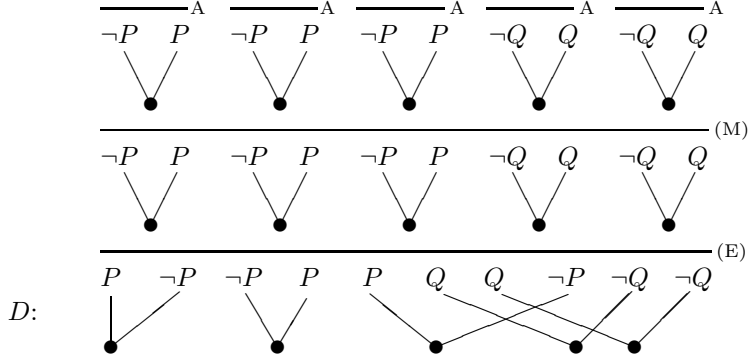
\begin{center}\begin{picture}(273,147)
\put(35,138){\line(1,0){33}}
\put(69,136){\scriptsize A}
\put(84,138){\line(1,0){33}}
\put(118,136){\scriptsize A}
\put(132,138){\line(1,0){33}}
\put(167,136){\scriptsize A}
\put(181,138){\line(1,0){33}}
\put(216,136){\scriptsize A}
\put(230,138){\line(1,0){33}}
\put(265,136){\scriptsize A}
\put(35,125){$\gneg P$}
\put(61,125){$P$}
\put(84,125){$\gneg P$}
\put(110,125){$P$}
\put(133,125){$\gneg P$}
\put(159,125){$P$}
\put(182,125){$\gneg Q$}
\put(208,125){$Q$}
\put(231,125){$\gneg Q$}
\put(257,125){$Q$}

\put(54,102){\line(1,2){10}}
\put(54,102){\line(-1,2){10}}
\put(54,102){\circle*{5}}

\put(103,102){\line(-1,2){10}}
\put(103,102){\line(1,2){10}}
\put(103,102){\circle*{5}}

\put(152,102){\line(-1,2){10}}
\put(152,102){\line(1,2){10}}
\put(152,102){\circle*{5}}

\put(201,102){\line(-1,2){10}}
\put(201,102){\line(1,2){10}}
\put(201,102){\circle*{5}}

\put(250,102){\line(1,2){10}}
\put(250,102){\line(-1,2){10}}
\put(250,102){\circle*{5}}

\put(35,92){\line(1,0){230}}
\put(267,90){\scriptsize (M)}
\put(35,79){$\gneg P$}
\put(61,79){$P$}
\put(84,79){$\gneg P$}
\put(110,79){$P$}
\put(133,79){$\gneg P$}
\put(159,79){$P$}
\put(182,79){$\gneg Q$}
\put(208,79){$Q$}
\put(231,79){$\gneg Q$}
\put(257,79){$Q$}

\put(54,56){\line(1,2){10}}
\put(54,56){\line(-1,2){10}}
\put(54,56){\circle*{5}}

\put(103,56){\line(-1,2){10}}
\put(103,56){\line(1,2){10}}
\put(103,56){\circle*{5}}

\put(152,56){\line(-1,2){10}}
\put(152,56){\line(1,2){10}}
\put(152,56){\circle*{5}}

\put(201,56){\line(-1,2){10}}
\put(201,56){\line(1,2){10}}
\put(201,56){\circle*{5}}

\put(250,56){\line(1,2){10}}
\put(250,56){\line(-1,2){10}}
\put(250,56){\circle*{5}}

\put(35,46){\line(1,0){231}}
\put(268,44){\scriptsize (E)}
\put(0,20){$D$:}

\put(35,33){$P$}
\put(57,33){$\gneg P$}
\put(82,33){$\gneg P$}
\put(110,33){$P$}
\put(134,33){$P$}
\put(158,33){$Q$}
\put(180,33){$Q$}
\put(204,33){$\gneg P$}
\put(228,33){$\gneg Q$}
\put(252,33){$\gneg Q$}

\put(39,10){\line(5,4){26}}
\put(39,10){\line(0,1){19}}
\put(39,10){\circle*{5}}

\put(102,10){\line(-3,5){12}}
\put(102,10){\line(3,5){11}}
\put(102,10){\circle*{5}}

\put(162,10){\line(-5,4){25}}
\put(162,10){\line(5,2){50}}
\put(162,10){\circle*{5}}

\put(215,10){\line(1,1){21}}
\put(215,10){\line(-5,2){50}}
\put(215,10){\circle*{5}}

\put(237,10){\line(1,1){21}}
\put(237,10){\line(-5,2){50}}
\put(237,10){\circle*{5}}

\end{picture}\end{center}

\noindent (In the pathological case of $D$ having no groups at all, it is simply the empty cirquent and hence an axiom.)

We conclude that $A$ is provable in (AMEWC$\mld\mlc$) and hence in {\bf CCC}.
\end{proof}

\begin{remark}\label{rem1}
Note that our completeness proof of Theorem \ref{clascomp} does not appeal to duplication, which means that 
{\bf CCC}=(AMEWC$\mld\mlc$). 
Duplication is thus syntactic sugar for {\bf CCC}. Furthermore, from the completeness proof given in the next section it can be seen that downward (though not upward) duplication does not really add anything to the deductive power of 
{\bf CL5}, either. How sweet is the sugar of duplication? 
It can certainly improve proof sizes, but the possible magnitude of that improvement is unknown at this point. One claim that we are making without a proof, however, is that {\bf CL5} minus duplication has polynomial-size proofs while still remaining strictly stronger than affine logic (remember that Blass's principle is provable in such a system but not in affine logic).  
\end{remark}

\section{Binary tautologies and their instances}

Let $C$ be a cirquent. An {\bf oatom} $P$ of $C$, i.e. an {\bf occurrence} of an atom $P$ in $C$, is an occurrence of $P$ in an oformula of (the diagram of) $C$. Such an oatom is {\bf negative} if it comes with a $\gneg$; otherwise it is {\bf positive}.
When $C$ is a cirquent or formula, by an ``{\bf atom of $C$}" we mean an atom that has at least one occurrence in $C$. 

A {\bf substitution} is a function $\sigma$ that sends every atom $P$ to some formula $\sigma(P)$; if such a $\sigma(P)$ always (for every $P$) is an atom, then $\sigma$ is said to be an {\bf atomic-level substitution}. Function $\sigma$ extends from atoms to all formulas in the expected way: $\sigma(\gneg P)=\gneg \sigma(P)$; $\sigma(F\mld G)=\sigma(F)\mld \sigma(G)$; $\sigma(F\mlc G)=\sigma(F)\mlc \sigma(G)$. 
$\sigma$ also extends to cirquents $C$ by stipulating that $\sigma(C)$ is the result of replacing in $C$ every oformula $F$ by $\sigma(F)$. 
 
Let $A$ and $B$ be cirquents. We say that $B$ is a (substitutional) {\bf instance} of $A$ iff $B=\sigma(A)$ for some substitution $\sigma$; and $B$ is an {\bf atomic-level instance} of $A$ iff $B=\sigma(A)$ for some atomic-level substitution $\sigma$.  
Example: the second cirquent of Figure 3 is an instance --- though not an atomic-level one --- of the first cirquent; the (relevant part of the) substitution $\sigma$ used here is defined by $\sigma(P)=Q\mld P$, $\sigma(Q)=Q$ and $\sigma(R)=P$.

\begin{center}\begin{picture}(365,83)
\put(0,66){\line(1,0){157}}
\put(0,53){$\gneg P$}
\put(43,53){$(Q\mlc R)\mld P$}
\put(117,53){$\gneg Q\mld \gneg R$}

\put(180,66){\line(1,0){187}}
\put(180,53){$\gneg Q\mlc\gneg P$}
\put(233,53){$(Q\mlc P)\mld (Q\mld P)$}
\put(326,53){$\gneg Q\mld \gneg P$}

\put(42,30){\line(-5,3){34}}
\put(42,30){\line(5,3){31}}
\put(42,30){\circle*{5}}
\put(104,30){\line(-5,3){32}}
\put(104,30){\line(5,3){34}}
\put(104,30){\circle*{5}}

\put(237,30){\line(-2,1){37}}
\put(237,30){\line(2,1){37}}
\put(237,30){\circle*{5}}
\put(310,30){\line(-2,1){37}}
\put(310,30){\line(2,1){37}}
\put(310,30){\circle*{5}}
\put(163,10){\bf Figure 3}
\end{picture}\end{center}

\begin{lemma}\label{l1}
If a given cirquent calculus system proves a cirquent $C$, then it also proves every instance of $C$.
\end{lemma}

\begin{proof} Consider a proof tree $T$ of an arbitrary cirquent $C$, and an arbitrary instance $C'$ of $C$. Let $\sigma$ be a substitution with $\sigma(C)=C'$. Replace every oformula $F$ of every cirquent of $T$ by $\sigma(F)$. It is not hard to see that the resulting tree $T'$, which uses exactly the same rules as $T$ does, is a proof of  $C'$. 
\end{proof}

A cirquent is said to be {\bf binary} iff no atom has more than two occurrences in it. A binary cirquent is said to be {\bf normal} iff, whenever it has two occurrences of an atom, one occurrence is negative and the other is positive. A {\bf binary tautology} (resp. {\bf normal binary tautology}) is a binary (resp. normal binary) cirquent that is a tautology in the sense of the previous section. This terminology also extends to formulas understood as cirquents. The left cirquent of Figure 3 is an example of a normal binary tautology. 
 
\begin{lemma}\label{may4}
A cirquent is an instance of some binary tautology iff it is an atomic-level instance of some normal binary tautology. 
\end{lemma}

\begin{proof} The ``if'' part is trivial. For the ``only if'' part, assume $A$ is an instance of a binary tautology $B$. Let $P_1,\ldots,P_n$ be all of the atoms of $B$ that have two positive or two negative occurrences in $B$. Let $Q_1,\ldots,Q_n$ be any pairwise distinct atoms not occurring in $B$. Let $C$ be the result of replacing in $B$ one of the two occurrences of $P_i$ by $Q_i$, for each $i=1,\ldots,n$. Then obviously $C$ is a normal binary cirquent, and $B$ an instance of it. By transitivity, $A$ (as an instance of $B$) is also an instance of $C$. 

We want to see that $C$ is a tautology. Deny this. Then there is a classical model $M$ in which $C$ is false. 
Let $M'$ be the model such that:
\begin{itemize}
\item $M'$ agrees with $M$ on all atoms that are not among  
$P_1,\ldots,P_n,Q_1,\ldots,Q_n$;
\item for each $i\in\{1,\ldots,n\}$,  $M'(P_i)=M'(Q_i)=\left\{\begin{array}{l}
\mbox{{\em false} if $P_i$ and $Q_i$ are positive in $C$}\\
\mbox{{\em true} if $P_i$ and $Q_i$ are negative in $C$}.\end{array}\right. $ 
\end{itemize}
By induction on complexity, it can be easily seen that, for every subformula $F$ of a formula of $C$, whenever $F$ is false in $M$, so is it in $M'$. This extends from (sub)formulas to groups of $C$ and hence $C$ itself. Thus $C$ is false in $M'$ because it is false in $M$. But $M'$ does not distinguish between $P_i$ and $Q_i$ (any $1\leq i\leq n$). This clearly implies that $C$ and $B$ have the same truth value in $M'$. That is, $B$ is false in $M'$, which is however impossible because $B$ is a tautology.
From this contradiction we conclude that $C$ is a (normal binary) tautology. 
  
Let $\sigma$ be a substitution such that $A=\sigma(C)$. Let $\sigma'$ be a substitution such that, 
for each atom $P$ of $C$, $\sigma'(P)$ is the result of replacing in $\sigma(P)$ each occurrence of each atom by a new atom in such a way that: (1) no atom occurs more than once in $\sigma'(P)$, and (2) whenever 
$P\not=Q$, no atom occurs in both $\sigma'(P)$ and $\sigma'(Q)$. As an instance of the tautological $C$,  
$\sigma'(C)$ remains a tautology (this follows from Lemma \ref{l1} and Theorem \ref{clascomp}). $\sigma'(C)$ can also be easily seen to be a normal binary cirquent, because $C$ is so. Finally, with a little thought, $A$ can be seen to be an atomic-level instance of $\sigma'(C)$. 
\end{proof}

\begin{lemma}\label{mar13}
The rules of mix, exchange, duplication, $\mlc$-introduction and $\mld$-introduction preserve binarity and normal binarity in both top-down and bottom-up directions. 
\end{lemma}

\begin{proof} This is so because the above five rules in no way affect what atoms occur in a cirquent and how many times they occur. 
\end{proof}

\begin{lemma}\label{mar13a}
Weakening preserves binarity and normal binarity in the bottom-up direction.  
\end{lemma}

\begin{proof} This is so because, in the bottom-up view, weakening can (delete but) never create any new occurrences of atoms. 
\end{proof}

\begin{theorem}\label{th4}
A cirquent is provable in {\bf CL5} iff it is an instance of a binary tautology.  
\end{theorem}

\begin{proof} ($\Longrightarrow$:) Consider an arbitrary cirquent $A$ provable in {\bf CL5}. By induction on the height of its proof tree, we want to show that $A$ is an
instance of a binary tautology. 

The above is obvious when $A$ is an axiom. 

Suppose now $A$ is derived by exchange from $B$. Let us just consider oformula exchange, with ogroup exchange being similar. By the induction hypothesis, $B$ is an instance of a binary tautology $B'$. Let $A'$ be the result of applying exchange to $B'$ ``at the same place" as it was applied to $B$ when deriving $A$ from it, as illustrated in the following example:

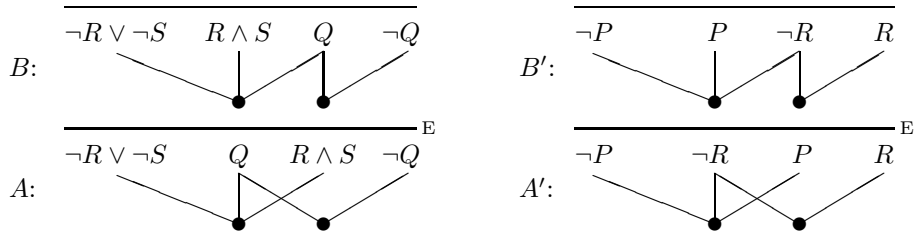
\vspace{10pt}

\begin{center}\begin{picture}(341,102)
\put(0,66){$B$:}
\put(21,92){\line(1,0){133}}
\put(21,79){$\gneg R\mld\gneg S$}
\put(74,79){$R\mlc S$}
\put(115,79){$Q$}
\put(141,79){$\gneg Q$}

\put(87,56){\line(-5,2){46}}
\put(87,56){\line(5,3){32}}
\put(87,56){\line(0,1){19}}
\put(87,56){\circle*{5}}

\put(119,56){\line(0,1){19}}
\put(119,56){\line(5,3){32}}
\put(119,56){\circle*{5}}

\put(21,46){\line(1,0){133}}
\put(156,44){\scriptsize E}
\put(0,20){$A$:}
\put(21,33){$\gneg R\mld\gneg S$}
\put(83,33){$Q$}
\put(106,33){$R\mlc S$}
\put(141,33){$\gneg Q$}

\put(87,10){\line(-5,2){46}}
\put(87,10){\line(5,3){32}}
\put(119,10){\line(-5,3){32}}
\put(87,10){\circle*{5}}

\put(87,10){\line(0,1){19}}
\put(119,10){\line(5,3){32}}
\put(119,10){\circle*{5}}

\put(193,66){$B'$:}
\put(214,92){\line(1,0){121}}
\put(214,79){$\gneg P$}
\put(264,79){$P$}
\put(290,79){$\gneg R$}
\put(327,79){$R$}

\put(267,56){\line(-5,2){46}}
\put(267,56){\line(5,3){32}}
\put(267,56){\line(0,1){19}}
\put(267,56){\circle*{5}}

\put(299,56){\line(0,1){19}}
\put(299,56){\line(5,3){32}}
\put(299,56){\circle*{5}}

\put(214,46){\line(1,0){121}}
\put(337,44){\scriptsize E}
\put(193,20){$A'$:}
\put(214,33){$\gneg P$}
\put(258,33){$\gneg R$}
\put(296,33){$P$}
\put(327,33){$R$}

\put(267,10){\line(-5,2){46}}
\put(267,10){\line(5,3){32}}
\put(299,10){\line(-5,3){32}}
\put(267,10){\circle*{5}}

\put(267,10){\line(0,1){19}}
\put(299,10){\line(5,3){32}}
\put(299,10){\circle*{5}}

\end{picture}\end{center}

\noindent Obviously $A$ will be an instance of $A'$. It remains to note that, by Lemmas \ref{tau} and \ref{mar13}, $A'$ is 
a binary tautology.  

The rules of duplication, $\mld$-introduction and $\mlc$-introduction can be handled in a similar way. 

Next, suppose $A$ is derived from $B$ and $C$ by mix. By the induction hypothesis, $B$ and $C$ are instances of some binary tautologies $B'$ and $C'$, respectively.   We may assume that no atom $P$ occurs in both $B'$ and $C'$, for otherwise, in one of the cirquents, rename $P$ into something different from everything else.  Let $A'$ be the result of applying mix to 
$B'$ and $C'$. By Lemmas \ref{tau} and \ref{mar13}, $A'$ is a binary tautology. And, as in the cases of the other rules, it is evident that $A$ is an instance of $A'$.  

Finally, suppose $A$ is derived from $B$ by weakening. If this is ogroup weakening, the conclusion 
is an instance of a binary tautology for the same reasons as in the case of exchange, duplication, $\mld$-introduction or $\mlc$-introduction. Assume now we are dealing with pool weakening, so that  
$A$ is the result of inserting a new oformula $F$ into $B$. By the induction hypothesis, $B$ is an instance 
of a binary tautology  $B'$.
Let $P$ be an atom not occurring in $B'$.
And let $A'$ be the result of applying weakening to  $B'$ that inserts $P$ in the same place into $B'$ as the above application of weakening inserted $F$ into $B$ when deriving $A$. Obviously $A'$ inherits binarity from $B'$; by 
Lemma \ref{tau}, it inherits from $B'$ tautologicity as well.  
And, for the same reasons as in all previous cases, $A$ is an instance of $A'$.\vspace{8pt}

($\Longleftarrow$:) Consider an arbitrary cirquent $A$ that is an  instance of a binary tautology $A'$. In view of Lemma \ref{l1}, it would suffice to show that {\bf CL5} proves $A'$.  
We construct a proof of $A'$, in the bottom-up fashion, as follows. Starting from $A'$, we keep applying conservative $\mld$-introduction and conservative $\mlc$-introduction until we hit an essentially literal cirquent $B$. As in the proof of Theorem \ref{clascomp}, such a cirquent $B$ is guaranteed to be a tautology, and $A'$ follows from it in ($\mld\mlc$). Furthermore, in view of Lemma 
\ref{mar13}, $B$ is in fact a binary tautology.
Continuing as in the proof of Theorem \ref{clascomp}, we apply to $B$ a series of weakenings and hit a 
tautological cirquent $C$ with no homeless oformulas, where every group only has two oformulas: $P$ and $\gneg P$ for some atom $P$. By Lemma \ref{mar13a}, $C$ remains binary. Our target cirquent $A'$ is thus derivable from $C$ in (W$\mld\mlc$). In the proof of Theorem \ref{clascomp} we next applied a series of contractions to separate shared oformulas. In the present case it suffices to use downward duplication (preceeded with ogroup exchange if necessary) instead of contraction: as it is easy to see, the binarity of $C$ implies that there are no shared oformulas in it except the cases when oformulas are shared by identical-content ogroups. Applying  to $C$ a series of ogroup exchanges and downward duplications, as illustrated below,  yields a cirquent $D$ that no longer has identical-content ogroups and hence no longer has any shared oformulas. 

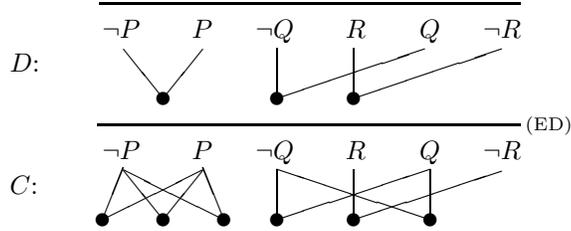
\begin{center}\begin{picture}(205,102)
\put(34,92){\line(1,0){159}}
\put(0,66){$D$:}
\put(35,33){$\gneg P$}
\put(69,33){$P$}
\put(93,33){$\gneg Q$}
\put(127,33){$R$}
\put(155,33){$Q$}
\put(179,33){$\gneg R$}

\put(35,10){\line(2,5){8}}
\put(35,10){\line(2,1){38}}
\put(35,10){\circle*{5}}

\put(58,10){\line(-4,5){15}}
\put(58,10){\line(4,5){15}}
\put(58,10){\circle*{5}}

\put(81,10){\line(-2,5){8}}
\put(81,10){\line(-2,1){38}}
\put(81,10){\circle*{5}}

\put(101,10){\line(0,1){19}}
\put(101,10){\line(3,1){58}}
\put(101,10){\circle*{5}}

\put(130,10){\line(0,1){19}}
\put(130,10){\line(3,1){56}}
\put(130,10){\circle*{5}}

\put(159,10){\line(0,1){19}}
\put(159,10){\line(-3,1){58}}
\put(159,10){\circle*{5}}

\put(33,46){\line(1,0){160}}
\put(195,44){\scriptsize (ED)}
\put(0,20){$C$:}

\put(35,79){$\gneg P$}
\put(69,79){$P$}
\put(93,79){$\gneg Q$}
\put(127,79){$R$}
\put(155,79){$Q$}
\put(179,79){$\gneg R$}

\put(58,56){\line(-4,5){15}}
\put(58,56){\line(4,5){15}}
\put(58,56){\circle*{5}}

\put(101,56){\line(0,1){19}}
\put(101,56){\line(3,1){56}}
\put(101,56){\circle*{5}}

\put(130,56){\line(0,1){19}}
\put(130,56){\line(3,1){56}}
\put(130,56){\circle*{5}}

\end{picture}\end{center}
 
$A'$ is thus derivable from $D$ in (EWD$\mld\mlc$). In turn, as in the proof of Theorem \ref{clascomp}, $D$ is provable in (AME). So, {\bf CL5} proves 
$A'$. \end{proof}

\section{Abstract resource semantics}

\subsection{Elaborating the abstract resource intuitions}\label{s8.1}

Among the basic notions used in our presentation of abstract resource semantics is that of {\em atomic resource}. This is an undefined concept, and we can only point at some examples of what might be intuitively considered atomic resources. These can be a specified amount of money; electric power of a specified voltage and amperage; a specified task performed by a computer, such as providing Internet browsing capabilities; 
a specified number of bits of memory; the standard collection of tasks/duties of an employee in a given enterprise; the choice between a candy and an apple that a vending machine offers to whoever inserts a \$1 bill into it; etc.

From atomic resources we will be building {\em compound resources}. Of course, whether a resource is considered atomic or thought of as a combination of some more basic resources depends on the degree of abstraction or encapsulation we choose in a given treatment. For instance, \$2 can be treated as an atomic resource, but it can as well be understood as a 
combination of \$1 and \$1 --- specifically, the combination $\$1\mlc\$1$, with $\alpha \mlc \beta$ generally being the resource having which intuitively means having both $\alpha$ and $\beta$. Similarly, a multiple-piece software package can be encapsulated and treated as an atomic resource, but in more subtle considerations it can be seen as a combination of the programs, data, etc. of which the package consists. And, had we extended our present approach to choice (``additive'') connectives, the earlier listed atomic resource providing a choice between a candy and an apple could be deatomized and  
understood as the choice conjunction $\adc$ of {\em Candy} and {\em Apple}.

When talking about resources, we always have two parties in mind: the resource provider and the resource user. Correspondingly, every entity that we call a resource comes in two flavors, depending on who is ``responsible'' for providing the resource. Suppose Victor received a salary of \$3,000 in the morning, and paid a \$3,000 mortgage bill in the afternoon. We are talking about the same resource of \$3,000, but in one case it came to Victor as an {\em income} ({\em input}), while in the other case it was an {\em expense} ({\em output}). In the morning Victor was the user (and his employer the provider), while in the afternoon he was the provider (and the mortgage company the user) of the resource \$3,000. 

Or, imagine a car dealer selling a Toyota to his customer for \$20,000. Two atomic resources can be seen involved in this transaction: the Toyota and the \$20,000. The \$20,000 is an income/input from the dealer's perspective while an expense/output from the customer's perspective; on the other hand, the Toyota is an expense/output for the dealer while an income/input for the customer. 

Or, compare the two devices: a mini power generator that produces 100 watts of electric power, and a lamp with a 100-watt light bulb. Power is an output for the generator while an input for the lamp. In turn, the generator does not produce power for free: it takes/consumes certain input such as fuel in a specified quantity. Similarly, the lamp outputs light in exchange for power input. 

The power generator and the lamp are also resources, in our treatment being compound ones unlike the atomic {\em Power}, {\em Fuel} or {\em Light}. Analyzing our intuitive concept of resources, one can notice that we are willing to call a resource anything that can be used --- perhaps in combination with some other resources --- to achieve certain goals. And ``achieving a goal'', in turn, can be understood as nothing but obtaining/generating certain resources. 
\$20,000 is a resource because it  can be used --- in combination with the resource {\em Car dealer} --- to obtain the resource {\em Toyota}. Similarly, a generator and a lamp can help us --- in combination with the resource {\em Fuel} --- obtain {\em Light}. The component atomic resources of {\em Generator} are {\em Fuel} and {\em Power}, the former being an input as already noted, and the latter being an output. We will be using the term {\em port} as a common name for inputs and outputs. To indicate that a given port is an input, the name of the corresponding atomic resource will be prefixed with a ``$-$''; the absence of such a prefix will mean that the port is an output. It should be noted that ``$-$'' is merely an indication of the input/output status of a port, and nothing more; it should not be mistaken for an operation on resources as, say, the later-defined operation $\gneg$ is. 

The sequence of all ports of a given compound resource we call its {\em interface}. Thus, the interface of the resource {\em Generator} is $\seq{\mbox{\em $-$Fuel, Power}}$, and the interface of {\em Lamp} is $\seq{\mbox{\em $-$Power, Light}}$. 
Generally, a compound resource may take any number of inputs and outputs.  For example, Victor possessing both a generator and a lamp can be seen as possessing just one compound resource {\em Generator}\hspace{2pt}$\mlc$\hspace{1pt}{\em Lamp}. 
The interface of this compound resource is then the concatenation of the interfaces of its two components, i.e. $\langle${\em $-$Fuel, Power, $-$Power, Light}$\rangle$. The first and the third ports of this interface stand for the resources that are expected to be provided by the user Victor, so they are inputs
and hence come with a ``$-$''; the second and the fourth ports, on the other hand, stand for the resources that Victor expects to receive, so they are outputs and hence come without a ``$-$''. As for Victor (as opposed to the provider of 
his resource {\em Generator}\hspace{2pt}$\mlc$\hspace{1pt}{\em Lamp}), he sees the same ports, but in negative colors: for him, the first and the third ports are outputs while the second and the fourth ports are inputs.  
 To visualize {\em Generator}\hspace{2pt}$\mlc$\hspace{1pt}{\em Lamp} as a one resource, it may be helpful for  us to imagine a generator and a lamp mounted together on one common board, with the
{\em $-$Fuel} port in the form of a pipe, the {\em Power} port in the form of a socket, the $-${\em Power} port in the form of a plug, and the {\em Light} port in the form of a light bulb. 
Multiple inputs and/or outputs are very common. A TV set takes two inputs: power and cable. And a look at the back panel of a personal computer will show a whole array of what the engineers indeed call ``ports''. This explains the choice of some of our terminology. 

The list of all inputs and outputs is only a half of a full description of a compound resource. The other half is what we call the resource's {\em truth function}. Formally, the latter is a function that returns a truth value --- 0 or 1 --- for each assignment of truth values to the ports of the compound resource. Intuitively, the value $1$ for a given resource --- whether it be atomic or compound --- means that the resource is ``functioning'', or ``doing its job'', or ``keeping its promise''. In such cases we will simply say that the resource is {\em true}. And the value $0$, as expected, means {\em false}, i.e. not true.  For instance, the value $1$ for the resource {\em Power} means that power is indeed generated/supplied, and the value $0$ means 
that there is no power supply. If Victor plugs the plug of his lamp into a (functioning, i.e. true) outlet, then the input port {\em $-$Power} of the resource {\em Lamp} becomes true; otherwise the port will probably remain false. We call assignments of truth values to the ports of a given compound resource {\em situations}. In these terms, the truth function $F$ of the compound resource tells us in which situations $\mathbf s$ the resource is considered true ($F({\mathbf s})=1$) and in which situations $\mathbf s$ it is false ($F({\mathbf s})=0$). Intuitively, such an $F$ can be seen as a description of the job that the resource is ``supposed'' (or ``promises'') to perform. Specifically, the job/promise of the resource is to be true, i.e. ensure that no situations $\mathbf s$ with $F({\mathbf s})=0$ will arise. Jobs are not always done and promises not always kept however. So a resource, whether elementary or compound,  may or may not be true.

Going back to the resource {\em Generator}, its job is to ensure that whenever there is fuel input, there is also power output. That is, whenever the (input) {\em $-$Fuel} port is true, so is the (output) {\em Power} port. 
In more intuitive but less precise terms, this job can be characterized as ``turning fuel into power''. The following two tables contain full descriptions of the two resources {\em Generator} and {\em Lamp}, each table showing both the interface and the truth function (the rightmost column) of the corresponding compound resource:\vspace{10pt}

\begin{center}\begin{picture}(330,138)

\put(0,128){\line(0,-1){98}}
\put(15,121){INTERFACE}
\put(44,111){\line(0,-1){81}}
\put(88,128){\line(0,-1){98}}
\put(92,103){\em Generator}
\put(141,111){\line(0,-1){81}}
\put(8,103){\em $-$Fuel}
\put(52,103){\em Power}
\put(0,94){\line(1,0){141}}
\put(20,83){0}
\put(64,83){0}
\put(113,83){1}
\put(0,78){\line(1,0){141}}
\put(20,67){0}
\put(64,67){1}
\put(113,67){1}

\put(0,62){\line(1,0){141}}
\put(20,51){1}
\put(64,51){0}
\put(113,51){0}
\put(0,46){\line(1,0){141}}
\put(20,35){1}
\put(64,35){1}
\put(113,35){1}
\put(0,30){\line(1,0){141}}

\put(190,128){\line(0,-1){98}}
\put(205,121){INTERFACE}
\put(234,111){\line(0,-1){81}}
\put(278,128){\line(0,-1){98}}
\put(292,103){\em Lamp}
\put(331,111){\line(0,-1){81}}
\put(194,103){\em $-$Power}
\put(244,103){\em Light}
\put(190,94){\line(1,0){141}}
\put(210,83){0}
\put(254,83){0}
\put(303,83){1}
\put(190,78){\line(1,0){141}}
\put(210,67){0}
\put(254,67){1}
\put(303,67){1}

\put(190,62){\line(1,0){141}}
\put(210,51){1}
\put(254,51){0}
\put(303,51){0}
\put(190,46){\line(1,0){141}}
\put(210,35){1}
\put(254,35){1}
\put(303,35){1}
\put(190,30){\line(1,0){141}}

\put(58,10){{\bf Figure 4:} The resources {\em Generator} and {\em Lamp}}
\end{picture}\end{center}

As we see from the above figure, the only situation in which {\em Generator} is false, i.e. considered to have failed to do its job, is $10$ --- the situation in which fuel was supplied to the generator but the latter did not produce power. While the situation $01$ is unlikely to occur with a real generator, {\em Generator} is considered true in it. For this is a situation in which the generator not only did not break its promise, but in fact did even more than promised. A customer who ended up receiving a \$20,000-priced Toyota for something less than \$20,000 (or for free) would hardly be upset and call the generous dealer a deal-breaker. The philosophy ``it never hurts to do more than necessary'' is inherent to our present approach, which is formalized in the requirement that the truth function of a resource should always be {\em monotone}, in the sense that changing $0$ to $1$ in an output or $1$ to $0$ in an input can never turn a true resource into a false one. Since {\em Generator} is true in situation $00$ (no input, no output), so should it be in $01$, for the generator ``produced even more than expected''; an equally good reason
why {\em Generator} is true in situation $01$ is that it is true in $11$, so that, in $01$, the generator ``consumed even less than expected''.   

We remember that the interface of the combination $\alpha\mlc \beta$ of resources is the concatenation of those of $\alpha$ and $\beta$. As for the truth function of $\alpha\mlc \beta$, it 
should account for the intuition that $\alpha\mlc \beta$ is considered to be doing its job iff both $\alpha$ and $\beta$ are doing their jobs. This can be seen from the following table for the resource {\em Generator}\hspace{2pt}$\mlc$\hspace{1pt}{\em Lamp}:

\begin{center}\begin{picture}(244,330)

\put(0,320){\line(0,-1){290}}
\put(40,303){\line(0,-1){273}}
\put(80,303){\line(0,-1){273}}
\put(120,303){\line(0,-1){273}}
\put(160,320){\line(0,-1){290}}
\put(244,303){\line(0,-1){273}}
\put(23,313){I \ N \ T \ E \ R \ F \ A \ C \ E}
\put(82,295){\em $-$Power}
\put(128,295){\em Light}
\put(5,295){\em $-$Fuel}
\put(47,295){\em Power}
\put(163,295){{\em Generator}\hspace{2pt}$\mlc$\hspace{1pt}{\em Lamp}}

\put(0,286){\line(1,0){244}}
\put(18,275){0}
\put(58,275){0}
\put(98,275){0}
\put(138,275){0}
\put(200,275){1}
\put(0,269){\line(1,0){244}}
\put(18,259){0}
\put(58,259){0}
\put(98,259){0}
\put(138,259){1}
\put(200,259){1}

\put(0,254){\line(1,0){244}}
\put(18,243){0}
\put(58,243){0}
\put(98,243){1}
\put(138,243){0}
\put(200,243){0}
\put(0,238){\line(1,0){244}}
\put(18,227){0}
\put(58,227){0}
\put(98,227){1}
\put(138,227){1}
\put(200,227){1}
\put(0,222){\line(1,0){244}}

\put(18,211){0}
\put(58,211){1}
\put(98,211){0}
\put(138,211){0}
\put(200,211){1}
\put(0,205){\line(1,0){244}}
\put(18,195){0}
\put(58,195){1}
\put(98,195){0}
\put(138,195){1}
\put(200,195){1}

\put(0,190){\line(1,0){244}}
\put(18,179){0}
\put(58,179){1}
\put(98,179){1}
\put(138,179){0}
\put(200,179){0}
\put(0,174){\line(1,0){244}}
\put(18,163){0}
\put(58,163){1}
\put(98,163){1}
\put(138,163){1}
\put(200,163){1}
\put(0,158){\line(1,0){244}}

\put(18,147){1}
\put(58,147){0}
\put(98,147){0}
\put(138,147){0}
\put(200,147){0}
\put(0,141){\line(1,0){244}}
\put(18,131){1}
\put(58,131){0}
\put(98,131){0}
\put(138,131){1}
\put(200,131){0}

\put(0,126){\line(1,0){244}}
\put(18,115){1}
\put(58,115){0}
\put(98,115){1}
\put(138,115){0}
\put(200,115){0}
\put(0,110){\line(1,0){244}}
\put(18,99){1}
\put(58,99){0}
\put(98,99){1}
\put(138,99){1}
\put(200,99){0}
\put(0,94){\line(1,0){244}}

\put(18,83){1}
\put(58,83){1}
\put(98,83){0}
\put(138,83){0}
\put(200,83){1}
\put(0,78){\line(1,0){244}}
\put(18,67){1}
\put(58,67){1}
\put(98,67){0}
\put(138,67){1}
\put(200,67){1}

\put(0,62){\line(1,0){244}}
\put(18,51){1}
\put(58,51){1}
\put(98,51){1}
\put(138,51){0}
\put(200,51){0}
\put(0,46){\line(1,0){244}}
\put(18,35){1}
\put(58,35){1}
\put(98,35){1}
\put(138,35){1}
\put(200,35){1}
\put(0,30){\line(1,0){244}}

\put(29,10){{\bf Figure 5:} The resource {\em Generator}\hspace{2pt}$\mlc$\hspace{1pt}{\em Lamp}}
\end{picture}\end{center}
    
The resource {\em Generator}\hspace{2pt}$\mlc$\hspace{1pt}{\em Lamp} is false in situations $1000$, $1001$ and  $1011$ because its {\em Generator} component failed to do its job: there was fuel input but no power output. The reason why  {\em Generator}\hspace{2pt}$\mlc$\hspace{1pt}{\em Lamp} is false in situations $0010$, $0110$ and  $1110$ is that 
{\em Lamp} malfunctioned: there was power input for it but no light output. And in situation $1010$ {\em Generator}\hspace{2pt}$\mlc$\hspace{1pt}{\em Lamp} is false because neither {\em Generator} nor {\em Lamp} kept the promise: there were both fuel input (for the generator) and power input (for the lamp), yet neither power output nor light output was generated. 

The value of {\em Generator}\hspace{2pt}$\mlc$\hspace{1pt}{\em Lamp} in a situation $xyzt$ only depends on the value $u$ of {\em Generator} in situation $xy$ and the value $v$ of {\em Lamp} in situation $zt$. So, such a $xyzt$ can be simply seen as $uv$, and the $\langle${\em $-$Fuel, Power$\rangle$ and  $\langle${\em $-$Power, Light$\rangle$ parts of the interface seen as simply 
{\em Generator} and {\em Lamp}, respectively. That is, the two compound conjuncts of  {\em Generator}\hspace{2pt}$\mlc$\hspace{1pt}{\em Lamp} can be encapsulated and treated as atomic resources, which yields the following, simpler table:

\begin{center}\begin{picture}(193,138)

\put(0,128){\line(0,-1){98}}
\put(10,121){I N T E R F A C E}
\put(53,111){\line(0,-1){81}}
\put(106,128){\line(0,-1){98}}
\put(111,103){{\em Generator}\hspace{2pt}$\mlc$\hspace{1pt}{\em Lamp}}
\put(193,111){\line(0,-1){81}}
\put(5,103){\em Generator}
\put(68,103){\em Lamp}
\put(0,94){\line(1,0){193}}
\put(25,83){0}
\put(78,83){0}
\put(147,83){0}
\put(0,78){\line(1,0){193}}
\put(25,67){0}
\put(78,67){1}
\put(147,67){0}

\put(0,62){\line(1,0){193}}
\put(25,51){1}
\put(78,51){0}
\put(147,51){0}
\put(0,46){\line(1,0){193}}
\put(25,35){1}
\put(78,35){1}
\put(147,35){1}
\put(0,30){\line(1,0){193}}

\put(-31,10){{\bf Figure 6:} The resource {\em Generator}\hspace{2pt}$\mlc$\hspace{1pt}{\em Lamp} in lesser detail}
\end{picture}\end{center}

Ignoring the minor --- at least seemingly so --- technical detail that columns are required to contain an additional bit of information indicating input/output status, 
our tables in the above style bear resemblance with those used in classical logic. Yet there is one crucial difference, which does not show itself in Figure 6 but catches the eye in Figure 5. In our tables, the same atom --- such as {\em Power} in Figure 5 --- may occur more than once, and the rows that assign different truth values to different occurrences of the same atom are as meaningful as any other rows. This is so because the {\em expressions} ``{\em Power}\hspace{1pt}'', ``{\em Fuel}\hspace{1pt}'', etc.,  as such, stand just for {\em resource types}, while particular {\em occurrences} of such expressions in a table (or a formula, cirquent etc.) stand for {\em individual resources} of those types. That is, we (again) deal with the necessity to differentiate between {\em ports} and {\em oports}, {\em inputs} and {\em oinputs}, etc. 
It is possible that the generator is producing power but the lamp is not receiving any: maybe Victor was not smart enough to insert the lamp's input plug into the generator's output socket, or it was a hot and sunny noontime, and he decided to use those 100 watts to feed his fan instead of the lamp.  In the particular case of Figure 5, the two {\em Power}-labeled columns happen to be of different genders: one an output and the other an input.  This, however, would not  always be the case, and generally a table can contain any number of columns with identical labels 
of either gender. For example, the non-simplified table for {\em Generator}\hspace{2pt}$\mlc$\hspace{1pt}{\em Generator} would have two (input) columns labeled {\em $-$Fuel}, and two (output) columns labeled {\em Power}. Such a 16-row table would not be at all the same as the 4-row one for the resource {\em Generator} seen in Figure 4. On the other hand, in classical logic, the truth table for any formula $F$ would be no different from that for $F\mlc F$, because for classical logic, which sees formulas as propositions or Boolean functions, $F$ and $F\mlc F$ are indistinguishable. 
{\em Generator} and {\em Generator}\hspace{2pt}$\mlc$\hspace{1pt}{\em Generator}, on the other hand, are certainly not the same as resources. Provided that Victor has enough fuel to feed two generators, with 
{\em Generator}\hspace{2pt}$\mlc$\hspace{1pt}{\em Generator} he can produce 200 watts of electricity while with just {\em Generator} only 100 watts. It may, however,  happen that Victor 
decides to provide input for the first generator but not for the second one, so that rows with a $1$ in the first {\em $-$Fuel} column and a $0$ in the second {\em $-$Fuel} column cannot be dismissed as meaningless or impossible. And even the fully classical-looking table of Figure 6 would no longer look classical with {\em Generator}\hspace{2pt}$\mlc$\hspace{1pt}{\em Generator}
instead of {\em Generator}\hspace{2pt}$\mlc$\hspace{1pt}{\em Lamp}. Such a table would still have 4 rows, while the classical table for $P\mlc P$ (atomic $P$) would only have 2 rows.

The meaning of the disjunction $\alpha\mld \beta$ of resources must be easy to guess. The interface of $\alpha\mld \beta$ is the same as that of $\alpha\mlc \beta$: just the interfaces of $\alpha$ and $\beta$ put together. As for the truth function, it corresponds to the intuition that the resource $\alpha\mld \beta$ is considered to have failed its job if and only if so did both of its components $\alpha$ and $\beta$. Hence, say, the table for {\em Generator}\hspace{2pt}$\mld$\hspace{1pt}{\em Lamp} would differ from the one of Figure 5 in that the last column would have a $0$ only in one --- the $1010$ --- row. The job of {\em Generator}\hspace{2pt}$\mld$\hspace{1pt}{\em Lamp} is thus to generate either power output or light output (or both) whenever {\em both} fuel and power inputs are received. 

The implicative combination $\alpha\mli \beta$ of resources, in rough intuitive terms, can be characterized as the resource that ``consumes $\alpha$ and produces $\beta$''. To see this, it would suffice to point out that our old friends {\em Generator} and {\em Lamp} (Figure 4) are nothing but {\em Fuel}\hspace{2pt}$\mli$\hspace{1pt}{\em Power} and {\em Power}\hspace{2pt}$\mli$\hspace{1pt}{\em Light}, respectively. Another intuitive characterization of $\mli$ is to say that this is a resource reduction operation. For example, the resource {\em Lamp} = {\em Power}\hspace{2pt}$\mli${\em Light} reduces (the task of generating) {\em Light} to (the task of generating) {\em Power}. Generally, the interface of $\alpha\mli \beta$ is the concatenation of those of $\alpha$ and $\beta$, only, in the $\alpha$ part, the input/output status of each port is reversed. That  is, in the antecedent the roles of the provider and the user are interchanged: the resource provider of $\alpha\mli \beta$ acts as a provider in the $\beta$ part but as a user in the $\alpha$ part. Indeed, {\em Generator} provides {\em Power}, but uses {\em Fuel}. This is why {\em Fuel} --- the atomic resource which, in isolation, is its own ``output'' --- is an input rather than an output of {\em Generator}.  The promise that the resource $\alpha\mli \beta$ carries is to make $\beta$ true as long as $\alpha$ is true. In other words, to guarantee that either $\beta$ true, or $\alpha$ is false (or both). Imagine a car dealer who promised to his customer to sell to him a Toyota for a ``to be negotiated'' price. One way to keep this promise is, of course, to actually sell a Toyota. But what if the dealer has run out of cars by the time the customer arrives? A way out for the dealer is to request an unreasonable price that he believes the customer would never be able or willing to pay.   

The above discussion makes it clear that $\alpha\mli \beta$ is, in fact, a disjunction. Specifically, it is $\gneg \alpha\mld \beta$, where $\gneg \alpha$, intuitively, is ``the opposite of $\alpha$'': the interface of $\gneg \alpha$ is that of $\alpha$ with all inputs turned into outputs and all outputs turned into inputs; and $\gneg \alpha$ is true in exactly the situations in which $\alpha$ is false. Alternatively and equivalently, we can define $\gneg \alpha$ as $\alpha\mli {\mathbf 0}$; here $\mathbf 0$ is the empty-interface resource that is just (``always'') false, intuitively meaning a resource that no one can ever provide.

We close this subsection with a look at an intuitive example where $\mli$ takes a compound antecedent. Let us imagine that Victor has {\em Fuel} and {\em Lamp}. Could he generate light if these two resources are all of his possessions?
 Not really. What he needs is a generator. That is, Victor cannot (successfully) provide the resource {\em Light}, but with his {\em Fuel}, {\em Lamp} and some thought, he can provide the weaker resource {\em Generator}\hspace{2pt}$\mli${\em Light}, i.e.  $(${\em Fuel}\hspace{2pt}$\mli${\em Power}$)$ $\mli$ {\em Light}. Or can he? ``Some thought'' was not originally listed among the resources of Victor, and he might have a hard time exercising it should the example 
be more complex than it is. This is where {\bf CL5} can come to help. As it turns out, {\bf CL5} is exactly the logic that provides a systematic, sound and complete answer to the question on {\em what} and {\em how} Victor can generally do. Back to our present example, with perfect knowledge of {\bf CL5}, Victor has a guarantee of success because {\bf CL5} proves
\[\mbox{\em Fuel}\mlc(\mbox{\em Power}\mli\mbox{\em Light}) \ \ \mli \ \ \bigl((\mbox{\em Fuel}\mli\mbox{\em Power})\mli \mbox{\em Light}\bigr).\]

\subsection{Resources and resource operations defined formally}

Before we move any further,  let us summarize, as formal definitions, the explanations given in the previous subsection. First of all,
we agree that what we have been calling ``atomic resources'' are nothing but propositional letters, i.e. atoms of the language underlying cirquent calculus. This, of course, is some abuse of concepts because, strictly speaking, 
the atoms of the language are {\em variables} ranging over atomic resources rather than atomic resources as such. Similar terminological liberty extends to the 
concepts formally defined below as ``ports'', (compound) ``resources'', etc. 

A {\bf port} is $P$ or $-P$, where $P$ is an atom called the {\bf type} of the port. A port which is just an
atom is said to be an {\bf output}, and a 
port which is a ``$-$''-prefixed atom is said to be an {\bf input}. The input/output status of a port is said to be the {\bf gender} of that port. 
The two genders {\em input}, {\em output} are said to be {\bf opposite}.

An {\bf interface} is a finite sequence of ports. A particular occurrence of a port (input, output) in an interface will be referred to as an {\bf oport} ({\bf oinput}, {\bf ooutput}). As we did with oformulas in the context of cirquents, 
we usually refer to an oport by the name of the corresponding port (as in the phrase ``the oport $P$''), even though different oports may be identical as ports. 

Let $I=\seq{X_1,\ldots,X_n}$ be an interface. A {\bf situation} for $I$ is a function ${\mathbf s}$ of the type $\{1,\ldots,n\}\rightarrow\{0,1\}$. We identify such a function ${\mathbf s}$ with the bit string $a_1\ldots a_n$, where $a_1={\mathbf s}(1),\ldots,a_n={\mathbf s}(n)$; we can also write ${\mathbf s}(X_i)$ instead of ${\mathbf s}(i)$, thinking of $\mathbf s$ as a function assigning truth values to oports. When  ${\mathbf s}(X_i)=1$, we say that 
$X_i$ is {\bf true} in $\mathbf s$, and when ${\mathbf s}(X_i)=0$, we say that 
$X_i$ is {\bf false} in $\mathbf s$.
We define the relation $\leq_I$ on situations for $I$ by stipulating that ${\mathbf s}\leq_I{\mathbf s}'$ iff, for each $i\in\{1,\ldots,n\}$, we have:
\begin{itemize}
\item if $X_i$ is an output, then ${\mathbf s}(X_i)\leq{\mathbf s}'(X_i)$;
\item if $X_i$ is an input, then ${\mathbf s}'(X_i)\leq{\mathbf s}(X_i)$.
\end{itemize}   
The relations $\geq_I, <_I$ and $>_I$ have the expected meanings: ${\mathbf s}\geq_I {\mathbf s}'$ iff ${\mathbf s}'\leq_I {\mathbf s}$; ${\mathbf s}<_I {\mathbf s}'$ iff ${\mathbf s}\leq_I {\mathbf s}'$ and ${\mathbf s}\not= {\mathbf s}'$; and  ${\mathbf s}>_I {\mathbf s}'$ iff ${\mathbf s}' <_I {\mathbf s}$.

\begin{definition}\label{res}
An {\bf abstract resource} --- henceforth simply ``{\bf resource}'' --- is a pair $\alpha=(\inp{\alpha},\fn{\alpha})$, where:

\begin{enumerate}
\item $\inp{\alpha}$, called the {\bf interface of $\alpha$}, is an interface.
\item $\fn{\alpha}$, called the {\bf truth function of $\alpha$}, is a function that sends every situation ${\mathbf s}$ for $\inp{\alpha}$ to $0$ or $1$, such that the following {\bf monotonicity condition} is satisfied:
\begin{itemize}
\item Whenever ${\mathbf s}\leq_{\inp{\alpha}} {\mathbf s}'$, we have $\fn{\alpha}({\mathbf s})\leq \fn{\alpha}({\mathbf s}')$. 
\end{itemize}    
\end{enumerate}
\end{definition}

When $\fn{\alpha}({\mathbf s})=1$ (resp. $=0$), we say that $\alpha$ is {\bf true} (resp. {\bf false}) in situation ${\mathbf s}$. ``{\bf Output} of $\alpha$'', ``{\bf oport of $\alpha$}'', ``{\bf situation for $\alpha$}'', etc. mean ``output of $\inp{\alpha}$'', ``oport of $\inp{\alpha}$'', ``situation for $\inp{\alpha}$'', etc. 

We need to agree on some notation used in the following definition of the basic resource operations. 
Let $I,I_1,I_2$ be interfaces. We write $-I$ to mean the interface that is the same as $I$ only with the genders of all oports reversed; that is, $-I$ is obtained from $I$ by deleting the prefix ``$-$'' wherever it was present and, simultaneously, adding such a prefix wherever it was absent. Next, $I_1I_2$ will stand for the concatenation of $I_1$ and $I_2$, i.e. the result of appending the oports of $I_2$ to those of $I_1$. Note that every situation (understood as a bit string) for $I_1I_2$ will have the form ${\mathbf s}_1{\mathbf s}_2$ (the concatenation of strings ${\mathbf s}_1$ and ${\mathbf s}_2$), where ${\mathbf s}_1$ is a situation for $I_1$ and ${\mathbf s}_2$ is a situation for $I_2$. The {\bf empty interface} will be denoted by 
$\langle\rangle$.  There is only one possible situation for it, and we denote that situation by $\epsilon$. With situations understood as bit strings, $\mathbf \epsilon$ is thus the {\bf empty bit string}.
  
\begin{definition}\label{neg}
Let $\alpha$, $\alpha_1$, $\alpha_2$ be resources. The operations $\gneg$, $\mlc$, $\mld$, $\mli$, ${\mathbf 0}$, $\mathbf 1$ are defined as follows:\vspace{5pt}

\noindent 1. 
{\bf Negation} $\gneg \alpha$:
\begin{itemize}
\item $\inp{\gneg \alpha}=-\inp{\alpha}$; 
\item For any situation ${\mathbf s}$ for $\gneg \alpha$, $\fn{\gneg \alpha}({\mathbf s})=1$ iff $\fn{\alpha}({\mathbf s})=0$. 
\end{itemize}

\noindent 2. 
{\bf Conjunction} $\alpha_1\mlc \alpha_2$:
\begin{itemize}
\item $\inp{\alpha_1\mlc \alpha_2}=\inp{\alpha_1}\inp{\alpha_2}$;
\item For any situation ${\mathbf s}={\mathbf s}_1{\mathbf s}_2$ for $\alpha_1\mlc \alpha_2$, where ${\mathbf s}_1$ is a situation for $\alpha_1$ and ${\mathbf s}_2$ is a situation for $\alpha_2$, we have $\fn{\alpha_1\mlc \alpha_2}({\mathbf s})=1$ iff $\fn{\alpha_1}({\mathbf s}_1)=1$ and $\fn{\alpha_2}({\mathbf s}_2)=1$. 
\end{itemize}

\noindent 3. {\bf Disjunction} $\alpha_1\mld \alpha_2$:
\begin{itemize}
\item $\inp{\alpha_1\mld \alpha_2}=\inp{\alpha_1}\inp{\alpha_2}$;
\item For any situation ${\mathbf s}={\mathbf s}_1{\mathbf s}_2$ for $\alpha_1\mld \alpha_2$, where ${\mathbf s}_1$ is a situation for $\alpha_1$ and ${\mathbf s}_2$ is a situation for $\alpha_2$, we have $\fn{\alpha_1\mld \alpha_2}({\mathbf s})=1$ iff $\fn{\alpha_1}({\mathbf s}_1)=1$ or $\fn{\alpha_2}({\mathbf s}_2)=1$. 
\end{itemize}

\noindent 4. {\bf Implication} $\alpha_1\mli \alpha_2$:
\begin{itemize}
\item $\inp{\alpha_1\mli \alpha_2}=(-\inp{\alpha_1})\inp{\alpha_2}$;
\item For any situation ${\mathbf s}={\mathbf s}_1{\mathbf s}_2$ for $\alpha_1\mli \alpha_2$, where ${\mathbf s}_1$ is a situation for $\alpha_1$ and ${\mathbf s}_2$ is a situation for $\alpha_2$, we have $\fn{\alpha_1\mli \alpha_2}({\mathbf s})=1$ iff $\fn{\alpha_1}({\mathbf s}_1)=0$ or $\fn{\alpha_2}({\mathbf s}_2)=1$. 
\end{itemize}

\noindent 5. {\bf Empty-interface} ({\bf constant}) resources ${\mathbf 0}$ and ${\mathbf 1}$:
\begin{itemize}
\item $\inp{\mathbf 0}=\inp{\mathbf 1}=\langle\rangle$;
\item $\fn{{\mathbf 0}}({\mathbf \epsilon})=0$;   \ $\fn{{\mathbf 1}}({\mathbf \epsilon})=1$.
\end{itemize}
\end{definition}

For safety, we need to verify that the above operations (with ${\mathbf 0}$ and ${\mathbf 1}$ seen as $0$-ary operations) are indeed operations on resources --- that is, that they do not violate the monotonicity condition of Definition \ref{res}. The truth function of ${\mathbf 0}$ (and similarly for ${\mathbf 1}$) is trivially monotone, because there is only one situation for the empty interface. Of the other operations, it would be sufficient to only consider $\mli$ for, in view of Theorem \ref{interdef}, just as this is the case in classical logic, $\gneg$, $\mlc$ and $\mld$ are definable in terms of $\mli$ and $\mathbf 0$.

Consider any two situations ${\mathbf s},{\mathbf s}'$ for $\alpha_1\mli \alpha_2$ such that ${\mathbf s}\leq_{\inp{\alpha_1\mli \alpha_2}}{\mathbf s}'$. Our goal is to show that if $\alpha_1\mli \alpha_2$ is false in ${\mathbf s}'$, then so is it in ${\mathbf s}$. Let ${\mathbf s}_1,{\mathbf s}'_1$ be the situations for $\alpha_1$, and   ${\mathbf s}_2,{\mathbf s}'_2$ the situations for $\alpha_2$, such that 
${\mathbf s}={\mathbf s}_1{\mathbf s}_2$ and ${\mathbf s}'={\mathbf s}'_1{\mathbf s}'_2$. Then the condition ${\mathbf s}\leq_{\inp{\alpha_1\mli \alpha_2}}{\mathbf s}'$ implies that \ 
${\mathbf s}_1\geq_{\inp{\alpha_1}}{\mathbf s}'_1$ and ${\mathbf s}_2\leq_{\inp{\alpha_2}}{\mathbf s}'_2$. Assume $\alpha_1\mli \alpha_2$ is false in ${\mathbf s}'$. This means that $\fn{\alpha_1}({\mathbf s}'_1)=1$ and $\fn{\alpha_2}({\mathbf s}'_2)=0$. From  $\fn{\alpha_1}({\mathbf s}'_1)=1$ and ${\mathbf s}_1\geq_{\inp{\alpha_1}}{\mathbf s}'_1$ the monotonicity of $\fn{\alpha_1}$ implies 
 $\fn{\alpha_1}({\mathbf s}_1)=1$. Similarly, from $\fn{\alpha_2}({\mathbf s}'_2)=0$ and ${\mathbf s}_2\leq_{\inp{\alpha_2}}{\mathbf s}'_2$ the monotonicity of $\fn{\alpha_2}$ implies 
 $\fn{\alpha_2}({\mathbf s}_2)=0$. This means that $\alpha_1\mli \alpha_2$ is false in ${\mathbf s}={\mathbf s}_1{\mathbf s}_2$. Done.
  
The following theorem can be verified by a routine examination of the relevant definitions, left as an exercise for the reader:

\begin{theorem}\label{interdef}
For any resources $\alpha,\alpha_1,\alpha_2$, the following equalities hold:
\[\begin{array}{rcl}
\gneg \alpha & = & \alpha\mli {\mathbf 0};\\
\gneg\gneg \alpha & = & \alpha;\\
\gneg(\alpha_1\mlc \alpha_2) & = & \gneg \alpha_1\mld \gneg \alpha_2;\\
\gneg(\alpha_1\mld \alpha_2) & = & \gneg \alpha_1\mlc \gneg \alpha_2;\\
\alpha_1\mli \alpha_2 & = & \gneg \alpha_1\mld \alpha_2.\\
\end{array}\]
\end{theorem}

\subsection{Formulas as resources}

By a {\bf situation} for a formula $F$ we mean an assignment $\mathbf s$ of truth values ($0$ or $1$) to its oatoms. Such an $\mathbf s$ will as well be understood as a situation for any osubformula $G$ of $F$ by mechanically restricting its domain to the oatoms of $G$.  Note the difference between a situation and a classical model. The latter is an assignment of 
truth values to {\em atoms} rather than {\em oatoms}. So, say, (the relevant part of) a classical model for $\gneg P\mld (P\mlc P)$ would only have to assign a value to $P$, while a situation for this formula would have to list three --- not necessarily identical --- values for the three oatoms of the formula. Classical models can be viewed as special cases of situations that assign identical truth values to oatoms that are identical as atoms. Disregarding this difference, the {\em truth} status of a formula $F$ in a situation $\mathbf s$ for $F$ is determined in the 
``standard'' way by recursion on (the occurrences of) its subformulas. That is:
\begin{itemize}
\item $\gneg P$ is true in $\mathbf s$ iff $P$ is false in $\mathbf s$; here and later ``false'', as always, means ``not true''. 
\item $G\mlc H$ is true in $\mathbf s$ iff so are both $G$ and $H$;
\item $G\mld H$ is true in $\mathbf s$ iff so is $G$ or $H$ (or both).
\end{itemize}

While the above definition looks like exactly the classical definition of truth, once again we emphasize the implicit and important difference: in the present definition, $\gneg P$, $G\mlc H$, $G\mld H$ and their subformulas are particular occurrences of (sub)formulas rather than formulas as such. This makes it possible for, say, $P\mlc\gneg P$ to be true, which will be the case when the first oatom $P$ is true and the second oatom $P$ is false.   

Notice that, graphically, literals and ports are almost the same, with the difference that ports take the prefix ``$-$'' where literals take the prefix ``$\gneg$''. Since we often convert between literals and ports, here we introduce two functions $\por$ and $\lit$, with $\por$ for literal-to-port conversion and $\lit$ for port-to-literal conversion. Specifically,  for any atom $P$, we have $\por(P)=\lit(P)=P$, $\por(\gneg P)=-P$ and $\lit(-P)=\gneg P$. 

As in the previous subsection, a situation $\mathbf s$ for a formula $F$ can and will be understood as a bit string --- specifically, a bit string of length $n$, where $n$ is the number of oatoms of $F$. Then the same $\mathbf s$ can also be considered a situation for any resource whose interface has $n$ ports. The second clause of the following definition implicitly relies on this seeing no difference between situations for formulas and situations for resources. 

\begin{definition}\label{june1}
With each formula $F$ we associate the resource $F^\clubsuit$, called the {\bf resource represented by $F$}, defined as follows:
\begin{itemize}
\item $\inp{F^\clubsuit}=\seq{\por(L_1),\ldots,\por(L_n)}$, where $L_1,\ldots,L_n$ is the sequence of the oliterals of $F$ listed in the order of their appearance in $F$. 
\item For every situation $\mathbf s$ for $F^\clubsuit$, $\fn{F^\clubsuit}({\mathbf s})=1$ iff $F$ is true in $\mathbf s$. 
\end{itemize}
\end{definition}

According to the following theorem, $^\clubsuit$ respects the meanings of $\gneg,\mlc,\mld$ (and, as we may guess, also of $\mli$ if the latter was officially allowed in formulas)  as operations on resources. Hence, with propositional letters understood as representing atomic resources,
$F^\clubsuit$ is indeed ``the resource represented by $F$''. The 
theorem can be verified by a straightforward analysis of the relevant definitions, so we state it without a proof:

\begin{theorem}\label{may31}
For any formulas $F,G$, we have $(\gneg F)^\clubsuit=\gneg (F^{\clubsuit})$, $(F\mlc G)^\clubsuit=F^\clubsuit\mlc G^\clubsuit$ and  $(F\mld G)^\clubsuit=F^\clubsuit\mld G^\clubsuit$.  
\end{theorem}

While every formula corresponds to a resource, vice versa does not hold. For example, it is obvious that no formula in our present sense does represent the resource $\mathbf 0$ or $\mathbf 1$, for there 
are no atomless formulas. This is not a serious issue of course, for we could painlessly add $\mathbf 0$ and $\mathbf 1$ (perhaps using the symbols $\tlg$ and $\twg$ instead) to our formal language, laziness being the only reason for not having done so.  

The resource whose interface is $\seq{P,Q,R}$ and which is true when so is $Q$ and at 
least one of $P,R$ cannot be expressed by a formula, either. For example, $Q\mlc(P\mld R)$ would not fit the bill because of the wrong order of its atoms. Again, this is not a ``serious'' problem, for we 
may be willing to not distinguish between {\em permutationally equivalent} resources --- resources that, informally speaking, only differ in the order in which their interfaces list the ports. 

There are, however, really serious reasons that make it impossible to capture all resources with formulas, reasons that essentially call for switching to non-traditional means of expression such as cirquents. The closure of atomic resources under $\gneg$, $\mlc$, $\mld$, $\mli$ and any other operations in a similar style does not yield the class of all resources, and this is so because, vaguely speaking, such operations do not allow us to account for the possibility of resource sharing. An analysis of the following example can make this point clear.

\begin{example}\label{may25}
Let $\beta$ be the resource defined by:
\begin{itemize}
\item $\inp{\beta}=\seq{P,Q,R}$;
\item $\fn{\beta}({\mathbf s})=1$ iff at least two of the three bits of ${\mathbf s}$ are $1$s.  
\end{itemize}

The ``two out of three'' Boolean function in classical logic would be expressed by the formula $F=(P\mlc Q)\mld (P\mlc R)\mld (Q\mlc R)$. $F^\clubsuit$, however, is not $\beta$, for the former has six oports rather than three. We need yet do not have means to indicate that, say, the two occurrences of $P$ in $F$ stand for the same individual resource rather than two different resources of the same type. In other words, we need yet do not have means to indicate that $P$ is {\em shared} between the two subresources $P\mlc Q$ and $P\mlc R$. It should not be hard for the reader to convince himself or herself that generally for no $\gneg,\mlc,\mld$-combination $F$ of $P,Q,R$ do we have $F^\clubsuit=\beta$.  
\end{example}

Thus, while the collection $(\gneg,\mlc,\mld)$ enjoys what is called {\em functional completeness} in classical logic, in no reasonable sense is the expressive power of this collection complete when it comes to resource semantics. Nor is there any easy remedy such as adding some extra connectives to the language. What we need is to go substantially beyond the traditional formalisms of logic. As we are going to see in the next subsection (Theorem \ref{may26}), the formalism of cirquent calculus turns out to be sufficient.

\subsection{Cirquents as resources}

Let $C$ be a cirquent.
By a {\bf situation} for $C$ we mean an assignment $\mathbf s$ of truth values ($0$ or $1$) to the oatoms of $C$. Understood as a bit string, such an ${\mathbf s}$ is the concatenation of situations for the oformulas of $C$. Then an oformula of the cirquent is considered true or false in $\mathbf s$ if it is so in the corresponding substring of $\mathbf s$. For instance, if the 
pool of $C$ is $\seq{G,H}$ where $G$ has 3 oatoms and $H$ has 2 oatoms,  
then $G$ is true in situation $10111$ iff it is true (in the sense of the previous subsection) in $101$, and $H$ is true in $10111$ iff it is true in $11$.  
Next, we consider a group $\Gamma$ of $C$ true in a given situation $\mathbf s$ for $C$ iff at least one of the oformulas of $\Gamma$ is true. Notice that, unlike oformulas, the truth 
values of two different ogroups that are identical as groups would always be the same. So, truth can be considered a property of groups rather than ogroups. Finally, $C$ is true in $\mathbf s$ iff all of its groups are so.

\begin{definition}\label{june1a}
With each cirquent $C$ we associate the resource $C^\clubsuit$, called the {\bf resource represented by $C$}, defined as follows:
\begin{enumerate}
\item $\inp{C^\clubsuit}=\seq{\por(L_1),\ldots,\por(L_n)}$, where $L_1,\ldots,L_n$ is the sequence of the oliterals of $C$ listed in the order of their appearance in $C$.
\item For every situation $\mathbf s$ for $C^\clubsuit$, \ $\fn{C^\clubsuit}({\mathbf s})=1$ iff $C$ is true in $\mathbf s$. 
\end{enumerate}
\end{definition}

\begin{theorem}\label{may26}
For every resource $\alpha$ there is a cirquent $C$ --- in fact a literal one --- which represents $\alpha$, i.e. such that $C^\clubsuit=\alpha$. 
\end{theorem}
 
\begin{proof} Consider an arbitrary resource $\alpha$ with $\inp{\alpha}=\seq{X_1,\ldots,X_n}$. We construct a corresponding cirquent $C$ as follows. The pool of $C$ is 
$\seq{L_1,\ldots,L_n}\ = \langle \lit(X_1),$ 
$\ldots,\lit(X_n)\rangle$.  In view of condition 1 of Definition \ref{june1a}, it is already clear that $\inp{C^\clubsuit}=\inp{\alpha}$. We now need to define the structure of $C$ and show that we also have $\fn{C^\clubsuit}=\fn{\alpha}$.   

Let us call a situation  $\mathbf s$  for $\alpha$ \ {\bf critical} \ iff $\alpha$ is false in $\mathbf s$ and,  for every situation ${\mathbf s}'$ for $\alpha$ with \ ${\mathbf s}<_{\inp{\alpha}}{\mathbf s}'$,  \mbox{$\alpha$ is true} in ${\mathbf s}'$. 
Let ${\mathbf s}_1,\ldots,{\mathbf s}_k$ be a list of all critical situations for $\alpha$. Seeing them as bit strings, these situations, of course, are also situations for 
$C$ (no matter what the structure of $C$ is). For each $i\in\{1,\ldots,k\}$, we define the group $\Gamma_i$ by
\[\Gamma_i=\{L_j \ |\ 1\leq j\leq n,\ L_j \mbox{\em \ is false in } {\mathbf s}_i\}.\] 

Now, we define the structure of $C$ to be $\seq{\Gamma_1,\ldots,\Gamma_k}$. To see that $\fn{C^\clubsuit}=\fn{\alpha}$,   
consider any situation $\mathbf s$ for $C$ (and hence for $C^\clubsuit$ and $\alpha$). 

Suppose $\fn{C^\clubsuit}({\mathbf s})=0$, by condition 2 of Definition \ref{june1a} meaning that $C$ is false in $\mathbf s$. Then there is a group $\Gamma_i$ ($1\leq i\leq k$) which is false in $\mathbf s$, and
this, in turn, means that every $L_j$ which is false in ${\mathbf s}_i$ (i.e. every oformula of $\Gamma_i$) is false in $\mathbf s$. In other words, ${\mathbf s}\leq_{\inp{\alpha}} {\mathbf s}_i$. But $\alpha$ is false in ${\mathbf s}_i$ because ${\mathbf s}_i$ is critical. Hence, by the monotonicity of $\fn{\alpha}$, we have $\fn{\alpha}({\mathbf s})=0$.  

Now suppose $\fn{\alpha}({\mathbf s})=0$. We may assume that $\mathbf s$ is critical, for otherwise replace it with a critical situation ${\mathbf s}'$ such that ${\mathbf s}\leq_{\inp{\alpha}} {\mathbf s}'$ (by monotonicity, such an ${\mathbf s}'$ is guaranteed to exist). Since $\mathbf s$ is critical, ${\mathbf s}={\mathbf s}_i$ for some $1\leq i\leq k$. Remembering now how the group $\Gamma_i$ was chosen, it is clear that $\Gamma_i$ is false in $\mathbf s$, whence $C$ is false in $\mathbf s$, i.e. $\fn{C^\clubsuit}({\mathbf s})=0$.
\end{proof}

\begin{example}\label{23}
Applying the construction from our proof of Theorem \ref{may26} to the resource $\beta$ of Example \ref{may25}, $\beta$ is represented by the following cirquent:
 
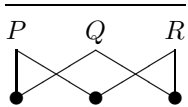
\begin{center} \begin{picture}(68,57)

\put(0,45){\line(1,0){68}}
\put(0,32){$P$}
\put(30,32){$Q$}
\put(60,32){$R$}

\put(4,10){\line(0,1){18}}
\put(64,10){\line(0,1){18}}
\put(4,10){\line(5,3){30}}
\put(64,10){\line(-5,3){30}}
\put(34,10){\line(-5,3){30}}
\put(34,10){\line(5,3){30}}
\put(4,10){\circle*{5}}
\put(34,10){\circle*{5}}
\put(64,10){\circle*{5}}
\end{picture}
\end{center}

\end{example}
\subsection{Resource-semantical validity}

Every logical semantics has a concept of validity, and so does our abstract resource semantics. In classical semantics valid formulas are called {\em tautologies}, and in abstract resource semantics they will be called {\em trivialities}. Trivialities and tautologies are similar in many respects. Asserting a tautology means ``asserting nothing''. Likewise, possessing (or providing) a triviality means ``possessing (or providing) nothing''. Of course, the word ``nothing'' has a negative flavor. But there is a positive side as well. If possessing $\alpha$ amounts to possessing nothing,
this means that, in fact, everyone possesses $\alpha$. Tautologicity is a guarantee of truth, and can be eventually used in finding true (and non-tautological) statements. Similarly, triviality 
is a guarantee of success in providing resources, and can be eventually used in finding what (nontrivial) resources can be generated. The formula given at the end of Subsection \ref{s8.1} is a triviality. And it is exactly this fact that allowed us to be confident that Victor, who possesses the nontrivial resources represented by the antecedent of the formula, can successfully generate the nontrivial resource represented by the consequent. 

The simplest examples of trivial resources are those of the form $\alpha\mli \alpha$. What makes such a resource trivial is that it merely returns back what it takes. Even the poorest person in the world would be able to provide the resource \$3,000$\hspace{1pt}\mli\hspace{1pt}$\$3,000, i.e. pay \$3,000 if he or she receives \$3,000. And it would not take a power plant to 
support the resource {\em Power}\hspace{2pt}$\mli${\em Power}\hspace{1pt}: assuming here that the input port comes to the resource provider in the form of a socket and the output port in the form of a plug, enough to just insert the plug into the socket. A similar trick works with {\em Power}$\hspace{1pt}\mlc${\em Power}\hspace{2pt}$\mli${\em Power},
providing which not only does not require any spending, but in fact can be even done with a benefit.  While $\alpha\mlc \alpha\mli \alpha$ is a triviality, $\alpha\mli \alpha\mlc \alpha$ is not. Only someone having \$3,000 of his or her own would be able to pay two \$3,000 bills while only receiving an income of \$3,000. And a duplex power adapter that seemingly turns {\em Power} into {\em Power}$\hspace{1pt}\mlc${\em Power}, does so only seemingly: while the voltage being accurately reproduced, the amperage in the two output (o)ports of the adapter would be inevitably lower than in its input port, as 100w cannot be converted into 200w ``for free''.     
  
What makes $P\mli P$, i.e. $\gneg P\mld P$ valid in classical logic is that the truth table for this formula has a $1$ in both of its two rows. The resource-semantical table for $\gneg P\mld P$, however, has a $0$ in one --- $10$ --- of its four rows. Generally, no formula in our present sense would have only $1$s in its resource-semantical table. 
Yet Victor, acting as a provider of $\gneg P\mld P$, has a way to make sure that the falsifying situation $10$ never occurs. This way is to {\em allocate} the input oport/resource $-P$ to the output oport/resource $P$. If here $P$ is 
{\em Power} with oport $-P$ given in the form of a socket and oport $P$ in the form of a plug, the physical meaning of allocation, as noted above, could be inserting the plug into the socket, which guarantees that whenever the socket has power supply, so does the plug. If $P$ is computer memory, then allocation may literally mean allocation. If $P$ is \$3,000, then allocation probably means redirection or mutual cancellation. If $P$ is {\em Light}, allocating probably means just letting the user utilize the light generated by himself or herself. For yet more diversity, assume $P$ is the resource {\em Chess} whose promise is to play the game of chess white and win. Then providing the resource $\gneg \mbox{\em Chess}\mld \mbox{\em Chess}$ means to play chess on two boards and win on at least one of them; specifically, to play white on the right (``output'') board, and play black on the left board, which is considered an ``input'' and hence the roles --- colors --- of the two players are interchanged. In this case, a way to ``allocate''  $-\mbox{\em Chess}$ and $\mbox{\em Chess}$ to each other   
would be to mimic on one board the moves made by the adversary on the other board, and vice versa. Obviously this strategy, amounting to having the user/adversary play against himself, guarantees that whenever the game is won by 
white (i.e. the promise of $\mbox{\em Chess}$ kept) on the left board, so is it on the right board.  

The concept of allocation thus allows various, rather different particular interpretations. Attempting to formalize those meanings (if so, which one?) is beyond the scope and dignity of our resource semantics which is meant to be an abstract, general-purpose formal framework. Essentially we treat allocation as a basic, undefined concept, just as we treat the concepts of atomic resources, inputs, outputs, ports, or truth (for atomic resources). Our semantics only focuses on the conditions that all intended interpretations of its basic concepts would satisfy. There are three such conditions pertaining to allocation: 
\begin{description}
\item[(i)] Only same-type and opposite-gender oports can be  allocated to each other.
\item[(ii)] Allocations should be ``monogamous'', in the sense that no oport can be involved in more than one allocation.
\item[(iii)] Each allocation guarantees that whenever its input oport is true, so is the output oport. 
\end{description}

Notice the egalitarian view of the genders implied by our phrases such as ``allocating to {\em each other}'': if oinput $X$ is allocated to ooutput $Y$, we also say that ooutput $Y$ is allocated to oinput $X$. This is in concordance with the fact that the intuitive distinction between input and output is often blurred and, when modeling real-life situations, initial decisions regarding whether a certain port should be listed as an input or an output can be arbitrary. For instance, when representing a computer as a compound resource, the (resource provided by the) monitor would most likely  become an  output and the keyboard an input; the fates of the modem or the floppy disk drive, however, are not just as clear. And,   in our earlier example with $\gneg\mbox{\em Chess}\mld\mbox{\em Chess}$, there were hardly any reasons in favor of treating the left rather than the right board as an ``input''.   In each case, however, we still deal with two opposite --- provider's and user's --- perspectives of the same resource, clearly dictating whether any given pair of same-type oports should have the {\em same} or {\em opposite} genders, even if there is flexibility in choosing {\em what} particular genders they have. So, under any choice, the resource modeling what was originally modeled as $\gneg\mbox{\em Chess}\mld\mbox{\em Chess}$ would have two same-type and opposite-gender oports. Note, however, 
that if genders are reversed, so should be the meanings of the truth values and all aspects of the provider/user perspectives of the corresponding atomic resources. For example, in the case of {\em Chess},
such a reversal would mean considering {\em Chess} (or whatever new name we use for it) true when black rather than white wins the game,
with the role of the provider of this resource now being that of the black rather than the white player. This  interchangeability and symmetry between input and output explains why many-inputs-to-one-output allocations would 
generally be just as inadmissible or impossible 
as many-outputs-to-one-input allocations. The fact that to-be-allocated resources more often than not 
 would involve elements of both ``income'' and ``expense'' (encapsulated {\em Generator} being a simplest example), 
with some thought, can also be seen to speak against the admissibility of any dual standards for inputs and outputs when it comes to allocations.

Back to the topic of validity, triviality of a resource is understood as a possibility to establish allocations that rule out all situations in which the resource is false. Assuming that making allocations only requires intellect and no other, ``external'' resources, and that Victor does possess intellect, the triviality of a resource indeed means a guarantee that Victor can always successfully provide the resource. The following formal definition summarizes these intuitions:
  
\begin{definition}\label{june2}
Let $\alpha$ be a resource. 

1. An {\bf allocation} for $\alpha$ is a pair $(X,Y)$, where $X$ is an oinput of $\alpha$, and $Y$ is an ooutput of $\alpha$ of the same type as $X$. We say that the allocation $(X,Y)$ {\bf utilizes} $X$ and $Y$.

2. An {\bf arrangement} for $\alpha$ is a set of allocations for $\alpha$. 

3. An arrangement ${\cal A}$ for $\alpha$ is said to be {\bf monogamous} iff no oport of $\alpha$ is utilized by more than one allocation of ${\cal A}$. 

4. 
Let ${\mathbf s}$ be a situation for $\alpha$, and ${\cal A}$ an arrangement for $\alpha$. We say that $\mathbf s$ is {\bf consistent with} ${\cal A}$ iff, whenever $(X,Y)\in {\cal A}$, 
${\mathbf s}(X)\leq {\mathbf s}(Y)$. 

5. We say that an arrangement ${\cal A}$ for $\alpha$ is {\bf trivializing} (for $\alpha$) iff $\alpha$ is true in every situation for $\alpha$ consistent with ${\cal A}$.

6. We say that resource $\alpha$ is {\bf trivial}, or is a {\bf triviality}, iff there is a monogamous trivializing arrangement for it.

7. We say that a cirquent (or formula) $C$ is {\bf trivial}, or is a {\bf triviality},  iff the resource represented by such a cirquent (formula) is trivial. That is, $C$ is trivial iff there is a monogamous trivializing arrangement for $C^\clubsuit$. 

\end{definition}

\begin{exercise}\label{june2b}
Show that dropping the monogamicity requirement in the above definition of triviality yields the ordinary concept of tautologicity for cirquents and formulas. That is, a cirquent or formula $C$ is a tautology in the sense of Section \ref{stau} iff there is a (not necessarily monogamous) trivializing arrangement for $C^\clubsuit$.

{\em Hint:} Let us say that an arrangement ${\cal A}$ for a given resource is {\bf greedy} iff ${\cal A}$ contains every possible allocation for that resource. Exploit the fact that, in a sense, classical models for $C$  are nothing but situations consistent with the greedy arrangement for $C^\clubsuit$. \end{exercise}

\begin{lemma}\label{may26a}
A cirquent (or formula) is trivial iff it is an instance of a binary tautology. 
\end{lemma}
 
\begin{proof} Formulas are special cases of cirquents, so let us limit our attention to cirquents in general. Consider an arbitrary cirquent $C$. For readability, here we will be identifying the literals of $C$ with the corresponding ports of $C^\clubsuit$, so that any arrangement for $C^\clubsuit$ can be seen as a set of pairs $(\gneg P,P)$ of oliterals of $C$. 

($\Longrightarrow$): Assume $C$ is trivial. Let ${\cal A}$ be a monogamous trivializing arrangement for $C^\clubsuit$.  Whenever two oatoms of $C$ are in literals utilized by the same allocation of ${\cal A}$, we call such oatoms {\bf coupled}. We replace oatoms of $C$ by new atoms in such a way that coupled oatoms are replaced by the same atom, and any pair of non-coupled atoms end up being replaced by different atoms. Call the resulting cirquent $D$. Of course, $D$ is a normal binary cirquent and $C$ is an (atomic-level) instance of it. What remains to show is that $D$ is a tautology. Suppose it is not. Then there must be a classical model 
$M$ in which $D$ is false. Let $\mathbf s$ be the situation for $C$ that sends every oatom $P$ of $C$ to the same truth value as $M$ sends the atom that replaced $P$ when obtaining $D$ from $C$. Clearly $\mathbf s$ is consistent with ${\cal A}$. It is also obvious that every 
osubformula of every oformula of $C$ has the same truth value in $\mathbf s$ as the corresponding osubformula of $D$ in $M$. We conclude that $C$ --- and hence $C^\clubsuit$ --- is false in $\mathbf s$. This contradicts our assumption that 
${\cal A}$ is trivializing for $C^\clubsuit$. 

($\Longleftarrow$): Assume $C$ is an instance of a binary tautology $D$. In view of Lemma \ref{may4}, we may assume that $D$ is normal and $C$ is an atomic-level instance of it. Let ${\cal A}$ be the arrangement for $C^\clubsuit$ that contains an allocation $(\gneg P,P)$ if and only if the oatoms of the corresponding two oliterals of $D$ are identical as atoms (``corresponding oliterals of $D$'' means the ones that were replaced by $\gneg P$ and $P$ when obtaining $C$ from $D$). Clearly ${\cal A}$ is monogamous. We want to verify that ${\cal A}$ is trivializing for $C^\clubsuit$. Suppose it is not. Then there is a situation ${\mathbf s}$ for $C^\clubsuit$ consistent with ${\cal A}$ such that $C^\clubsuit$ is false in $\mathbf s$. In view of the monotonicity of $\fn{C^\clubsuit}$, we may assume that whenever two oatoms of $C$ are coupled (``coupled'' in the same sense as in the previous paragraph), $\mathbf s$ assigns identical truth values to them. Let then $M$ be a classical model such that, for any atom $P$ of $D$, the truth value of $P$ in $M$ is the same as the truth value of the corresponding oatom(s) of $C$ in $\mathbf s$ (again, ``corresponding oatom(s)'' mean(s) the one(s) that replaced $P$ when obtaining $C$ from $D$). It is not hard to see that then $D$ has the same truth value --- particularly, the value {\em false} --- in $M$ as $C$ has in $\mathbf s$. So, $D$ is not a tautology, which is a contradiction. 
\end{proof}

The following theorem, establishing the soundness and completeness of {\bf CL5} with respect to our abstract resource semantics, is an immediate corollary of Theorem \ref{th4} and Lemma \ref{may26a}:

\begin{theorem}\label{bbb}
A cirquent (or formula) is provable in {\bf CL5} iff it is trivial. 
\end{theorem}

\begin{remark}\label{rem2}
At the end of Subsection \ref{s8.1} we promised that {\bf CL5} would provide an ultimate answer to the question what and how
 Victor can achieve purely by means of smart resource management. Theorem \ref{bbb} directly pertains only to the {\em what}
 part of this question: Victor can successfully provide the resource represented by a formula or cirquent $C$ (i.e. $C$
 is trivial) iff {\bf CL5}$\vdash C$. Specifically, he can achieve his success by just setting up a monogamous trivializing arrangement for 
$C^\clubsuit$ (making allocations is exactly what ``resource management'' means). So, the {\em how} part of the above question reduces to finding such an arrangement. 
From our proof of the $(\Longleftarrow$) part of Lemma \ref{may26a} one can see that finding a monogamous trivializing arrangement for $C^\clubsuit$, in turn, essentially means nothing but finding a normal binary cirquent $D$ such that $C$ is an atomic-level instance of $D$.
Our proof of the soundness part of Theorem \ref{th4} implicitly provides an easy way to turn a {\bf CL5}-proof of $C$ into a {\bf CL5}-proof of a binary (and, in fact, normal) tautological cirquent $C'$ such that $C$ is an instance of $C'$. $C$ is not necessarily an atomic-level instance of such a $C'$ though. No problem: our constructive proof of Lemma \ref{may4} shows how to turn $C'$ into a(nother) normal binary tautology $D$ such that $C$ is an atomic-level instance of $D$. 
Putting the above steps and observations together yields a polynomial-time algorithm that turns a {\bf CL5}-proof of an arbitrary cirquent $C$ into a monogamous trivializing arrangement for $C^\clubsuit$.  
\end{remark}

\section{Cirquent calculus and computability logic}\label{slast}

The claim that the semantics of computability logic is a semantics of resources has been explicitly or implicitly present in every work on CL. 
The abstract resource semantics introduced in this paper makes the same claim. While the two semantics are far from being the same, there is no contradiction or competition here. This is so not only due to the forthcoming Theorem \ref{june11} which, in conjunction with Theorem \ref{bbb}, implies that the two semantics validate the same principles. Abstract resource semantics, as noted, is a general-purpose framework, with its basic notions such as (atomic) resources, their truth values, or allocations being open to various specific interpretations. Computability logic offers one of such interpretations, and its semantics can be seen as a materialization of abstract resource semantics. A ``resource'' in CL has a very specific meaning as briefly explained in Section \ref{intr}. It is a derived concept defined in terms of some other, more basic and subtler-level entities (specifically, games). The same can be said about the other basic semantical concepts: in CL, 
to what we call ``successfully providing the resource'' (truth value $1$) corresponds winning the game; to our ``allocations'' correspond copy-cat subroutines in game-playing strategies; etc. While being more special, however, the semantics of computability logic still turns out to be general enough to invalidate anything that abstract resource semantics does.   
    
The earlier-promised soundness and completeness of {\bf CL5} with respect to the semantics of computability logic can be proven purely syntactically, based on the fact of the soundness and completeness of system {\bf CL2} known from \cite{CL2}. The propositional language of the latter is considerably more expressive than the one in which the formulas of {\bf CL5} are written. One difference is that the language of {\bf CL2} includes the {\em choice operators} $\adc$ and $\add$. The syntactic behavior of these operators, as noted in Section \ref{intr}, is somewhat reminiscent of that of the additive operators of linear logic, just as $\mlc,\mld,\mli$, called {\em parallel operators} in CL, are relatives of the multiplicative operators of linear logic. Another significant difference is that the language of {\bf CL2} has two sorts of atoms: {\bf general} and {\bf elementary}.\footnote{In computability logic, general atoms stand for computational problems (games) of arbitrary degrees of interactivity, while the interpretations of elementary atoms are limited to a special, ``moveless'' sort of games, i.e. computational problems with no interaction at all. Such problems are in fact nothing but propositions or predicates in the standard sense.} General atoms are the same as in the formulas of {\bf CL5}, and we continue using the uppercase $P,Q,R,S$ as metavariables for them. As for elementary atoms, they are foreign to {\bf CL5}. Elementary atoms are divided into {\bf logical} and {\bf non-logical}. There are two logical elementary atoms $\twg$ and $\tlg$, and infinitely many non-logical elementary atoms, for which we will be using the lowercase $p,q,r,s$ as metavariables. There is also a minor difference: in the official version of {\bf CL2} given in \cite{CL2}, formulas are allowed to contain $\mli$, and also the scope of $\gneg$ is not limited to atomic formulas. We may safely pretend that this difference does not exist: whether it be classical, linear or computability logic, $F\mli G$ is virtually the same as $\gneg F\mld G$ and hence can be understood just as an abbreviation. For similar reasons, it is irrelevant whether $\gneg$ is allowed to be applied only to atoms or compound formulas as well. Another insignificant difference is that the language of {\bf CL2} officially treats $\mlc$ and $\mld$ (as well as $\adc$ and $\add$) as variable-arity operators, as opposed to the strictly binary treatment chosen in the present paper. In view of the associativity of these operators in classical, linear and computability logics, this difference, too, can be safely ignored. So, in this presentation we assume that the official language of {\bf CL2} does not include $\mli$, that $\gneg$ is only allowed to be applied to atoms, and that $\mlc,\mld,\adc,\add$ are strictly binary. 

What we have been referring to as ``formulas'' in the previous sections we will now call {\bf CL5-formulas}, and formulas of the language of {\bf CL2} we will call {\bf CL2-formulas}. {\bf CL5}-formulas are thus nothing but {\bf CL2}-formulas that do not contain $\adc,\add$ and elementary atoms. A {\bf CL2}-formula is said to be {\bf elementary} iff it 
contains neither general atoms nor $\adc$,$\add$. Ignoring any differences (in their status) between the general and  elementary sorts of atoms, in this section $\adc,\add$-free {\bf CL2}-formulas --- including elementary formulas --- will be seen as formulas of classical propositional logic, with $\twg$ and $\tlg$ having their standard meanings (``truth'' and ``falsity'', respectively). The terms ``positive occurrence'' and ``negative occurrence'' will also be used with their standard meanings as we did before. A {\bf surface occurrence} of a subformula of a {\bf  CL2}-formula is an occurrence which is not in the scope of $\adc,\add$. The {\bf elementarization} of a {\bf CL2}-formula $F$ is the result of replacing in $F$ every positive surface occurrence of each general atom by $\tlg$, every negative surface occurrence of each general atom by $\twg$, every surface occurrence of each $\add$-subformula by $\tlg$, and every surface occurrence of each $\adc$-subformula by $\twg$. A {\bf CL2}-formula is said to be {\bf stable} iff its elementarization is a tautology of classical logic; otherwise the formula is said to be {\bf instable}. In these terms, and with ${\mathcal P}\mapsto C$ meaning ``from premise(s) $\mathcal P$ conclude $C$'', system {\bf CL2} is given in \cite{CL2} by the following three rules of inference:        

\begin{description}
\item[Rule (a):]  $\vec{H}\mapsto F$, where $F$ is stable and $\vec{H}$ is the smallest set of formulas such that, 
whenever $F$ has a surface occurrence of a subformula $G_1\adc G_2$, for both  
$i\in\{1,2\}$, $\vec{H}$ contains the result of replacing that occurrence in $F$ by $G_i$.
\item[Rule (b):]  $H\mapsto F$, where $H$ is the result of replacing in $F$ a surface occurrence of a subformula $G_1\add G_2$ by $G_1$ or $G_2$. 
\item[Rule (c):] $H\mapsto F$,\label{z18} where $H$ is the result of replacing in $F$ two --- one positive and one negative ---
surface occurrences of some general atom by a nonlogical elementary atom that does not occur in $F$.
\end{description}

Axioms are not explicitly stated, but notice that the set of the premises of Rule {\bf (a)} can be empty, in which case the conclusion acts as an axiom. Let us look at a couple of examples before moving any further. {\bf CL2} proves $P\mlc P\mli P$, i.e. $(\gneg P\mld\gneg P)\mld P$. This formula follows from $(\gneg P\mld \gneg p)\mld p$ by Rule {\bf (c)}. In turn, the stable 
formula $(\gneg P\mld \gneg p)\mld p$ is derivable from the empty set of premises by Rule {\bf (a)}. On the other hand, {\bf CL2} does not prove $P\mli P\mlc P$, i.e. $\gneg P\mld (P\mlc P)$. Indeed, this formula is instable and does not contain choice operators, so the only rule by which it could be derived is {\bf (c)}. The premise should be $\gneg p\mld (p\mlc P)$ or 
$\gneg p\mld (P\mlc p)$ for some nonlogical elementary atom $p$. In either case we deal with an instable formula without choice operators and with only one occurrence of a general atom. Hence it cannot be the conclusion of any of the three rules of {\bf CL2}. 

\begin{lemma}\label{june9}
A {\bf CL5}-formula is provable in {\bf CL2} iff it is an instance of a binary tautology.
\end{lemma}

\begin{proof} In this proof, as agreed a while ago, we treat both general and elementary atoms as atoms of classical propositional logic. With this minor and irrelevant difference, the terms ``instance'', ``binary tautology'' etc. have the same meanings as before.\vspace{7pt}
 
($\Longrightarrow$:) Consider an arbitrary {\bf CL5}-formula $F$ provable in {\bf CL2}. Fix a {\bf CL2}-proof of $F$ in the form of a sequence (rather than tree) $\seq{F_n,F_{n-1},\ldots,F_1}$ of formulas, with $F_1=F$. We may assume that this sequence has no repetitions or other redundancies. We claim that, for each $i$ with $1\leq i\leq n$, the following conditions are satisfied:

\begin{description}
\item[Condition 1:] $F$ is an atomic-level instance of $F_i$, and (hence) $F_i$ does not contain $\adc,\add$.
\item[Condition 2:] Whenever $F_i$ contains an elementary atom, that atom is non-logical, and has exactly two --- one positive and one negative --- occurrences in $F_i$.
\item[Condition 3:] If $i<n$, then $F_i$ is derived from $F_{i+1}$ by Rule {\bf (c)}.
\item[Condition 4:] $F_n$ is derived (from the empty set of premises) by Rule {\bf (a)}.
\end{description}

Condition 4 is obvious, because it is only Rule {\bf (a)} that may take no premises. That Conditions 1-3 are also satisfied can be verified by induction on $i$. For the basis case of $i=1$, Conditions 1 and 2 are immediate because $F$ ($=F_1$) is its own atomic-level instance and, as a {\bf CL5}-formula, it contains neither $\adc,\add$ nor elementary atoms. $F_1$ cannot be derived by Rule {\bf (b)} because, by Condition 1,  $F_1$  does not contain any $\add$. Nor can it be derived by Rule {\bf (a)} unless $n=1$, for otherwise either $F_1$ would have to contain a $\adc$ (which is not the case according to Condition 1), or the proof of $F$ would have redundancies as $F_1$ would not really need any premises. Thus, if $1<n$, the only possibility for $F_1$ is to be derived from $F_2$ by Rule {\bf (c)}.    
For the induction step, assume $i<n$ and the above conditions are satisfied for $F_i$. According to Condition 3, $F_i$ is derived by Rule {\bf (c)} from $F_{i+1}$. This obviously implies that $F_{i+1}$ inherits Conditions 1 and 2 from $F_i$. And that Condition 3 also holds for $F_{i+1}$ can be shown in the same way as we did for $F_1$. 

Condition 1 implies that $F$ is an instance of $F_n$. Therefore, in order to complete our proof of the $(\Longrightarrow)$ part of the lemma, it would now suffice to show that $F_n$ is an instance
 of a binary tautology. As the conclusion of Rule {\bf (a)} (Condition 4), $F_n$ is stable. Let $G$ be the elementarization of $F_n$. The stability of $F_n$ means that $G$ is a tautology.  Let $H$ be the result of replacing in $G$ every
 occurrence of $\twg$ and $\tlg$ by an elementary nonlogical atom not occurring in $G$, in such a way that different occurrences of $\twg,\tlg$ are replaced by different atoms. In view of Condition 2 (applied to
 $F_n$), it is obvious that $H$ is binary. By the same condition, $F_n$ did not contain $\twg,\tlg$. This means that every occurrence of $\twg$ in $G$ comes from replacing a negative occurrence of a general atom in $F_n$, and every occurrence of $\tlg$ in $G$ comes from replacing a positive occurrence of a general atom in $F_n$. It is not hard to see that, for this reason, 
$F_n$ is an instance of $H$. What remains to show is that $H$ is a tautology. But this is indeed so because $H$ results
 from the tautological $G$ by replacing positive occurrences of $\tlg$ and negative occurrences of $\twg$. It is known from classical logic that such replacements do not destroy truth and hence tautologicity of formulas.\vspace{7pt}  

($\Longleftarrow$:) Assume $F$ is a {\bf CL5}-formula which is an instance of a binary tautology $T$. In view of Lemma \ref{may4}, we may assume that $T$ is normal and $F$ is an atomic-level instance of it. For simplicity, we may also assume that all atoms of $T$ are elementary (elementary and non-logical, that is). Let us call the atoms that only have one occurrence in $T$ {\bf single}, and the atoms that have two occurrences {\bf married}. Let $\sigma$ be the substitution with $\sigma(T)=F$. Let $G$ be the formula resulting from $T$ by substituting each single atom $q$ by $\sigma(q)$. It is clear that then $F$ can be derived from $G$ by a series of applications of Rule {\bf (c)}, with each such application replacing two --- a positive and a negative --- occurrences of some married atom $p$ by $\sigma(p)$. So, in order to show that 
{\bf CL2} proves $F$, it would suffice to verify that $G$ is stable and hence it can be derived from the empty set of premises by Rule {\bf (a)}. But $G$ is indeed stable. To see this, consider the elementarization $G'$ of $G$. It results from $T$ by replacing its single atoms by $\twg$ (if negative) or $\tlg$ (if positive). It is known from classical logic that replacing a single atom by whatever formula does not destroy tautologicity of formulas. Hence, as $T$ is a tautology, so is $G'$, meaning that $G$ is stable. 
\end{proof}

\begin{theorem}\label{june11}
A formula is provable in {\bf CL5} iff it is valid in computability logic.
\end{theorem}
\begin{proof} This theorem is an immediate corollary of Lemma \ref{june9}, Theorem \ref{th4} and the known fact (proven in \cite{CL2}) that 
{\bf CL2} is sound and complete with respect to the semantics of computability logic. \end{proof}

Theorems \ref{bbb} and \ref{june11} are similar in that both establish a soundness and completeness of {\bf CL5}. A difference that may catch the eye is that Theorem \ref{bbb} talks about 
all cirquents while the statement of Theorem \ref{june11} is restricted only to formulas. This is so because cirquents
 first emerged in the present paper, and hence the semantics of CL elaborated earlier did not extend to them. Here we 
briefly outline how to fill this gap. This explanation is very informal, and is only meant for readers already familiar
 with the semantics of CL. Every $n$-ary cirquentstructure (see Section \ref{mul}) can be seen as an $n$-ary operation on
 games, and every cirquent $C$ with that structure then seen as the result of applying that operation to the games 
$A_1,\ldots,A_n$ represented by the oformulas of the pool of the cirquent. Specifically, $C$ is a game playing which
 means playing the $n$ games in parallel, i.e. it is a game on $n$ boards, where $A_1$ is played on board \#1, $A_2$ is played on board \#2, etc. The machine (player $\twg$) is considered the winner iff, for each group $\Gamma$ of the cirquent, 
it wins at least one of the games (represented by the oformulas) of the group. This is the whole story. As we see, cirquents can be understood as parallel combinations of games in the same style as $\mlc,\mld$-combinations. Not every combination represented by cirquents can be expressed using the ordinary (such as $\mlc,\mld$) operators of 
CL. Among such combinations is the one expressed by the cirquent of Example \ref{23}. We claim without a proof that,  
with the just-outlined extended semantics of CL, Theorem \ref{june11} generalizes to all cirquents. 

We further claim that Theorem \ref{june11} can be strengthened by replacing {\bf CL5} with the more expressive cirquent calculus system {\bf CL6}. What makes the language of {\bf CL6} stronger than that of {\bf CL5} is that, along with general atoms, {\bf CL6}-formulas may also contain elementary atoms including the logical atoms $\twg$ and $\tlg$, with $\gneg \twg$ now always written as $\tlg$ and $\gneg\tlg$ as $\twg$.  In accordance with the earlier-established meaning of the term ``elementary'', we call a {\bf CL6}-formula {\bf elementary} iff it does not contain general atoms. With ``formulas'' now meaning any {\bf CL6}-formulas, the set of the rules of {\bf CL6} is obtained from that of {\bf CL5} by adding to it $\twg$ (understood as a singleton cirquent) as an additional axiom, plus the rule of contraction limited only to elementary formulas (that is, the contracted formula $F$ is required to be elementary). 
One could show that Theorem \ref{june11} remains true with {\bf CL6} instead of {\bf CL5}. So does it with ``formula'' further replaced by ``cirquent'', where computability-logic validity of cirquents is understood in the sense 
of the previous paragraph.  

\section{What is next?}\label{slastest}

The results of the previous section can be seen as first steps within the program of 
developing cirquent-calculus deductive systems for ever more expressive fragments of computability logic. 
Probably the same applies to our present version of abstract resource semantics.
 Linear-logic literature abounds with 
examples illustrating the intended resource intuitions for additive and exponential operators. Yet, 
just as in the case of multiplicatives, such intuitions have never found a good formal explication,
the reason for which, again, being the inherent incompleteness of linear logic. If and when abstract resource semantics is successfully and naturally extended to additive- and exponential-style operators, one should probably expect it to validate the same class of formulas as computability logic does, thus yielding a logic properly stronger than  linear (affine) logic in its full language. 

It would also be interesting to see if abstract resource semantics can be modified so as to make a meaningful semantics for {\em relevance logic}, specifically, for the obviously relevance-logic-style system 
which results from {\bf CCC} by deleting weakening (but not contraction) and perhaps 
duplication as well. 

Back to proof theory, besides adopting more expressive underlying languages for formulas, a promising direction for cirquent calculus to grow might be relaxing its concept of a cirquentstructure. As illustrated in Section \ref{mul}, cirquentstructures in the sense of the present paper are special sorts of circuits --- ones 
of depth 2, with all level-1 gates being $\mld$-gates and the single level-2 gate being an 
$\mlc$-gate. Allowing circuits of arbitrary depths and forms (including circuits with non-traditional types of gates) as underlying cirquentstructures would yield more general concepts of cirquents. It is appropriate to still use the same term ``cirquent'' in these cases. After all, the gist of the cirquent calculus approach is its ability to syntactically capture resource-sharing, and resource-sharing is exactly  
what circuits are all about. 

As long as obtaining soundness and completeness is the only goal, the semantics of computability logic --- of its $(\gneg,\mlc,\mld)$-fragment for sure --- does not call for 
generalizations in the above style, just as classical logic needs nothing beyond ordinary sequent calculus. But even if,
in some good sense, the present form of cirquent calculus is natural and sufficiently powerful as a ``universal syntax'' 
(not that the author really wants to make such a claim),  
studying generalizations of it could be still worthwhile. Switching to the new syntactic vision known as the {\em calculus of structures} (\cite{Gug01}) introduced by Alessio Guglielmi has proved very beneficial, among the benefits being exponential speedups of proofs, and the possibility of richer combinatorial analysis of proofs for diverse logical systems, including the old and well-axiomatized classical, linear or modal logics. Achieving or enhancing similar effects could be some of the possible motivations for tackling general forms of cirquent calculus.  
A central idea in 
the calculus of structures is what is called {\em deep inference}.  It means deducing inside a formula at any depth, as opposed to the shallow inference of sequent calculus where rules only see the roots of formula trees. Notice 
that all three rules of {\bf CL2} modify subformulas at any ($\mlc,\mld$)-depths,  thus essentially being deep inference rules.  Some time ago this observation  gave rise to the hope that the calculus of structures could be 
a well-suited deductive framework for CL.\footnote{See http://alessio.guglielmi.name/res/cos/crt.html\#ComLog} 
However, efforts to axiomatize even the simplest fragments of CL in it have not been successful so far.
While CL obviously does need deep inference, what is  making the calculus of structures apparently insufficient is the absence of resource-sharing capabilities.  
Just like the rules of {\bf CL2}, our cirquent calculus rules from Section \ref{srules} can be seen as limited sorts of deep inference, specifically, inferences that always take place at depth 2 (when they modify oformulas) or 1 (when they modify ogroups). With arbitrary-depth underlying cirquentstructures, generalized cirquent calculus would naturally invite inferences at any structural levels, thus enjoying the full generality, flexibility and power offered by both deep inference and resource-sharing.  

Last but maybe not least, it would be interesting to understand the impact of sharing on {\em cut 
elimination}, the idol worshiped by all proof-theoreticians. Of course, before attempting to even ask this question,
one needs to figure out what a cut rule should exactly mean in the context of cirquent calculus.

\end{document}